\documentclass[12equipped]{article}
\usepackage{amssymb,amsmath}
\usepackage{hyperref}
\usepackage{graphicx}
\input epsf.tex
\usepackage[utf8]{inputenc}

\title{Notes on $A_{\infty}$-algebras, $A_{\infty}$-categories and non-commutative geometry. I}

\author {Maxim Kontsevich and Yan Soibelman}

\begin{document}
\maketitle
\newtheorem{thm}{Theorem}[subsection]
\newtheorem{defn}[thm]{Definition}
\newtheorem{lmm}[thm]{Lemma}
\newtheorem{rmk}[thm]{Remark}
\newtheorem{prp}[thm]{Proposition}
\newtheorem{conj}[thm]{Conjecture}
\newtheorem{exa}[thm]{Example}
\newtheorem{cor}[thm]{Corollary}
\newtheorem{que}[thm]{Question}
\newtheorem{ack}{Acknowledgements}
\newcommand{\C}{{\bf C}}
\newcommand{\K}{{\bf k}}
\newcommand{\R}{{\bf R}}
\newcommand{\N}{{\bf N}}
\newcommand{\Z}{{\bf Z}}
\newcommand{\Q}{{\bf Q}}
\newcommand{\G}{\Gamma}
\newcommand{\A}{A_{\infty}}
\newcommand{\g}{\bf g}
\newcommand{\ihom}{\underline{Hom}}
\newcommand{\epi}{\twoheadrightarrow}
\newcommand{\mono}{\hookrightarrow}
\newcommand\ra{\rightarrow}
\renewcommand\O{{\cal O}}
\newcommand\nca{nc{\bf A}^{0|1}}
\newcommand{\epp}{\varepsilon}

\tableofcontents

\section{Introduction}

\subsection{$\A$-algebras as spaces}\label{spaces}
The notion of $\A$-algebra introduced by Stasheff  as well as the  notion of $\A$-category introduced
by Fukaya has two different versions. First one is operadic:
an $\A$-algebra is an algebra over the $\A$-operad\footnote{The $\A$-operad can be realized as the operad of singular chains of the operad of intervals
in the real line.} Second one is geometric: an $\A$-algebra is the
same as a non-commutative formal graded manifold $X$ over a field
$\K$, having a marked $\K$-point $pt$, and  equipped with a vector
field $d$ of degree $+1$ such that $d|_{pt}=0$ and  $[d,d]=0$ (such
vector fields are called {\it homological}). By definition the
algebra of functions on the non-commutative formal pointed graded
manifold is isomorphic to the algebra of formal series $\sum_{n\ge
0}\sum_{i_1,i_2,...,i_n\in
I}a_{i_1...i_n}x_{i_1}...x_{i_n}:=\sum_Ma_Mx^M$ of free graded
variables $x_i,i\in I$ (the set $I$ can be infinite). Here
$M=(i_1,...,i_n), n\ge 0$ is a non-commutative multi-index, i.e. an
element of the free monoid generated by $I$. Homological vector
field makes the above graded algebra into a complex of vector spaces. The triple
$(X,pt,d)$ is called a {\it non-commutative formal pointed
differential-graded (or simply dg-) manifold}.

It is an interesting problem to make a dictionary from the pure
algebraic language of $\A$-algebras and $\A$-categories to the
language of non-commutative geometry \footnote {We use ``formal"
non-commutative geometry in tensor categories, which is different
from the non-commutative geometry in the sense of  Alain Connes.}.
One purpose of these notes is to make few steps in this
direction.

From the point of view of Grothendieck's approach to the notion of ``space",
our formal pointed manifolds are given by functors on graded associative Artin algebras
commuting with finite projective limits. It is easy to see that such functors
are represented by graded coalgebras. These coalgebras can be thought of as coalgebras
of distributions on formal pointed manifolds. The above-mentioned algebras of formal
power series are dual algebras to the coalgebras of distributions.

In the case of small $\A$-categories one should modify the above definitions. Instead of one
marked point one considers a closed subscheme of disjoint points
 in a formal graded manifold (they corresponds to objects of the category), and the  homological vector
field $d$ must be compatible with the embedding of this subscheme as
well as with the projection onto it.

\subsection{Some applications of the geometric language}\label{geometric language}
Geometric approach to $\A$-algebras and
$\A$-categories gives answers to several questions which are not transparent from the pure algebraic point of view. In
particular one can obtain an explicit description of the
$\A$-structure on $\A$-functors.  In geometric language
$\A$-functors are interpreted as maps between non-commutative formal
dg-manifolds commuting with homological vector fields. One can
introduce a non-commutative formal dg-manifold of maps between two
such spaces. Functors are just ``commutative" points of the latter.\footnote{By commutative points
we mean points with values in commutative algebras.} 
As a result of the above
considerations one can describe explicitly the $\A$-structure on functors in terms
of sums over  certain sets of trees.
In this paper we provide the answer in the case of $\A$-categories with one object, i.e. in the case of $\A$-algebras. 
The general case reflects the difference between quivers with
one vertex and quivers with many vertices (vertices correspond to objects). \footnote
{Another, purely algebraic approach to the $\A$-structure on
functors was suggested in [Lyu].}  Among other applications of our
geometric language  we mention an interpretation of the Hochschild
chain complex of an $\A$-algebra in terms of cyclic differential forms
on the corresponding formal pointed dg-manifold (Section \ref{Hochschild cochain complex}).

Geometric language simplifies some proofs as well.
For example, Hochschild cohomology of an $\A$-category $\mathcal{C}$
is isomorphic to $Ext^{\bullet} ( Id_{\mathcal{C}}, Id_{\mathcal{C}} )$ taken
in the $\A$-category of endofunctors $\mathcal{C} \to \mathcal{C}$. This
result admits an easy proof, if one interprets Hochschild cochains as vector
fields and functors as maps\footnote{The idea to treat $Ext^{\bullet} ( Id_{\mathcal{C}}, Id_{\mathcal{C}} )$
as the tangent space to deformations of the derived category $D^b({\cal C})$ goes back
to A.Bondal.}

\subsection{Contents of the paper}\label{contents}
Present paper contains two parts out of the planned three (the
last one is devoted to $\A$-categories).
In this paper  we discuss
$\A$-algebras (=non-commutative formal pointed dg-manifolds with fixed affine coordinates). We have tried to
develop in detail the dictionary between algebra and geometry and provide  proofs in the cases when they are not straightforward.

Part \ref{dg-manifolds} is devoted to the geometric description of $\A$-algebras. We start with basics on formal
graded affine schemes, then add  homological vector fields to the story,
thus arriving to the geometric definition of $\A$-algebras as formal
pointed dg-manifolds.
Most of the material is well-known in  algebraic language.
We cannot completely avoid $\A$-categories (subject of  another planned
paper). They appear in the form of categories of $\A$-modules
and $\A$-bimodules, which we define using our geometric language.

Since in the $\A$-world many notions are defined ``up to quasi-isomorphism", their geometric meaning is not obvious.
As an example we mention the notion of {\it weak unit}.
Basically, the latter means that the unit exists at the level of
cohomology only. In Section \ref{weak unit} we discuss the relationship
of weak units with the ``differential-graded" version of the affine line.

We start Part \ref{smoothness and compactness} with the definition of the Hochschild complexes
of $\A$-algebras. As we have already mentioned,
Hochschild cochain complex is interpreted in terms of
graded vector fields on the non-commutative formal affine space.
Dualizing, Hochschild
chain complex is interpreted in terms of degree one cyclic differential forms.
This interpretation is motivated by [Ko2]. It differs from the traditional
picture (see e.g. [Co], [CST]) where one assigns to a
Hochschild chain $a_0 \otimes a_1 \otimes ... \otimes a_n$ the differential
form $a_0 da_1 ...da_n$. In our approach we interepret $a_i$ as the dual to an
affine coordinate $x_i$, and the above expression is dual to the cyclic differential $1$-form $x_1 ...x_ndx_0$.
We also discuss graphical
description of Hochschild chains, the differential, etc.

After that we discuss homologically smooth compact $\A$-algebras.
Those are analogs of smooth projective varieties in algebraic
geometry. Indeed, the derived category $D^b(X)$ of coherent sheaves
on a smooth projective variety $X$ is $\A$-equivalent to the category of
perfect modules over a homologically smooth compact $\A$-algebra
(this can be obtained using the results of [BvB]). The latter
contains as much information about the geometry of $X$ as the
category $D^b(X)$ does. A good illustration of this idea is given by
the  non-commutative version of the Hodge theory presented in Section \ref{Hodge-to-de Rham}. It
is a largely conjectural topic, which eventually should be
incorporated in the theory of ``non-commutative motives". Encoding
smooth proper varieties by homologically smooth compact
$\A$-algebras one can forget about the underlying commutative
geometry and develop a theory of ``non-commutative smooth
projective varieties" in an abstract form. 

Let us briefly explain how all that looks in the case of  the Hodge theory. Let $( C_{\bullet} ( A, A ), b )$
be the Hochschild chain complex of a weakly unital homologically
smooth compact $\A$-algebra $A$. The corresponding negative cyclic
complex $( C_{\bullet} ( A, A ) [ [ u ] ], b + uB )$ gives rise to a
family of complexes over the formal affine line ${\bf A}^1_{form} [
+2 ]$ (shift of the grading reflects the fact that the variable $u$
has degree $+ 2$, cf. [Co], [CST]). We conjecture that the
corresponding family of cohomology groups gives rise to a vector
bundle over the formal line. The generic fiber of this vector bundle
is isomorphic to the periodic cyclic homology of $A$, while the fiber over $u =
0$ is isomorphic to the Hochschild homology of $A$. If a compact
homologically smooth $\A$-algebra $A$ comes from a smooth
projective variety as explained above, then the generic fiber is
just the algebraic de Rham cohomology of the variety, while the
fiber over $u = 0$ is the Hodge cohomology. In this case our conjecture
becomes the classical theorem which claims degeneration of  the
spectral sequence Hodge-to-de Rham.\footnote{In a recent preprint
[Kal], D.Kaledin claims the proof of our conjecture. He uses a
different approach to the cyclic homology theory than we do.}.

Last section of Part  \ref{smoothness and compactness} is devoted to the relationship between moduli
spaces of points on a cylinder and algebraic structures on the
Hochschild complexes.
In Section \ref{generalized Deligne\rq{}s conjecture} we formulate a generalization of Deligne's conjecture.
Recall that Deligne's conjecture
says (see e.g. [KoSo1], [McS], [T]) that the Hochschild cochain complex of an
$\A$-algebra is  an algebra over the operad of chains on the
topological operad of little discs.
In the conventional approach
to non-commutative geometry Hochschild cochains correspond to polyvector
fields, while Hochschild chains correspond to de Rham differential forms. One
can contract a form with a polyvector field or take a Lie
derivative of a form with respect to a polyvector field. This geometric point of view
leads to a generalization of Deligne's conjecture which includes Hochschild chains equipped
with the structure of (homotopy) module over cochains, and to the ``Cartan
type" calculus which involves both chains and cochains (cf. [CST],
[TaT1]). We unify both approaches under one roof formulating a theorem which
says that the pair consisting of the Hochschild chain and
Hochschild cochain complexes of the same $\A$-algebra is
an algebra over the colored operad of singular chains on
configurations of discs on a cylinder with marked points on each of the
boundary circles.{\footnote {After our paper was finished we received
the paper [TaT2] where the authors proved an equivalent result.}}

Sections \ref{scalar product} and \ref{PROP} are devoted to $\A$-algebras with scalar product. Geometrically this is the same
as non-commutative formal symplectic manifolds. In Section \ref{scalar product} we also discuss a homological
version of this notion and explain that it corresponds to the notion of Calabi-Yau structure on a manifold.
In Section \ref{PROP} we define an action of the PROP of
singular chains of the topological PROP of smooth oriented $2$-dimensional surfaces with
boundaries on the Hochschild chain complex of  an $\A$-algebra with scalar product. 
If in addition the algebra $A$ is homologically smooth 
and the spectral sequence Hodge-to-de Rham degenerates, 
then the above action extends to the action of the PROP 
of singular chains on the topological PROP of stable $2$-dimensional surfaces. This is
essentially equivalent to a structure of $2$-dimensional 
Cohomological TFT \footnote{A different approach to this topic was developed by Kevin Costello, see [Cos]}. 
More details and an application of this approach to Fukaya categories
is given in [KaKoP]. We also mention that in Sections \ref{GM connection}, 
\ref{flat connection and operad} we propose an approach to the Gauss-Manin connection on periodic cyclic homology which is different
from the one in [Get]. In particular we conjecture the existence of flat connection which is derived from a certain colored operad.

\subsection{Generalization to $\A$-categories}\label{generalization to categories}
Let us say few words about the  planned
paper  devoted to $\A$-categories.

The formalism of present paper admits a straightforward
generalization to the case of $\A$-categories.
The latter should be viewed
as non-commutative formal dg-manifolds with a closed marked subscheme
of objects. One can hope for a  translation of the theory of $\A$-categories into our geometric
language of non-commutative geometry. This goal is not fully achieved at present.
Although some parts of the theory of $\A$-categories
admit nice interpretation in terms of non-commutative geometry, some other
still wait for it. The latter includes e.g. triangulated $\A$-categories. 
As a first step to the future geometric description of triangulated $\A$-categories
we propose to describe 
their theory  from the point of view of
$\A$-functors from ``elementary " categories to a given $\A$-category (see a summary in [Ko5], [So1], [So2]).
Those ``elementary" categories are, roughly speaking, derived
categories of representations of quivers with small number of vertices.
Our approach has certain advantages over the traditional one. For example the complicated ``octahedron
axiom" admits a natural interpretation in terms of functors from
the $\A$-category associated with the quiver of the Dynkin diagram $A_2$ (there are six
indecomposible objects in the category $D^b ( A_2 - mod )$ corresponding to
six vertices of the octahedron). Interpreting such functors as morphisms of non-commutative formal dg-schemes
one can hope to describe some parts of the theory of $\A$-categories geometrically.

At the same time we have not been able to provide pure geometric proofs of all the  results, thus
relying on less flexible approach which uses differential-graded categories
(see [Dr]). 

In present and subsequent paper on the subject we mostly consider $\A$-algebras and
categories over a field of characteristic zero. This assumption
simplifies many results, but also makes some other less general. We
refer the reader to [Lyu], [LyuOv] for a theory over a ground ring
instead of ground field (the approach of [Lyu],[LyuOv] is pure
algebraic and different from ours). Most of the results of present
paper are valid for an $\A$-algebra $A$ over the unital commutative
associative  ring $\K$, as long as the graded module $A$ is flat over
$\K$. More precisely, the results of Part \ref{dg-manifolds} remain true except of the
results of Section \ref{minimal models} (the minimal model theorem). In
these two cases we assume that $\K$ is a field of characteristic
zero. {\it Constructions} of Part \ref{smoothness and compactness} work over a commutative ring
$\K$. The results of Section \ref{scalar product} are valid (and the conjectures are
expected to be true) over a field of characteristic zero. Algebraic
version of Hodge theory from Section \ref{Hodge-to-de Rham} and the results of Section
\ref{Hochschild cochain and chain complexes} are formulated for an $\A$-algebra over the field of
characteristic zero, although e.g. the Conjecture \ref{integer version} is expected to be
true for any $\Z$-flat $\A$-algebra.

{\it Acknowledgments}. We thank to Vladimir Drinfeld for useful
discussions and to Victor Ginzburg and Kevin Costello for comments on the manuscript. Y.S. thanks to Clay Mathematics Institute for supporting him as a
Fellow during a part of the work. His work was also partially supported by an NSF grant.
He is especially grateful to IHES for
excellent research and living conditions.

\part{$\A$-algebras and non-commutative dg-manifolds}\label{dg-manifolds}

\section{Coalgebras and non-commutative schemes}\label{coalgebras}

Geometric description of $\A$-algebras will be given in terms of geometry of
non-commutative ind-affine schemes in the tensor category of graded vector spaces
(we will use $\Z$-grading or $\Z/2$-grading). In
this section we are going to describe these ind-schemes as functors from
finite-dimensional algebras to sets (cf. with the description of formal
schemes in [Gr]). More precisely, such functors are represented by counital
coalgebras. Corresponding geometric objects are called {\it non-commutative
thin schemes}.

\subsection{Coalgebras as functors}

Let $\K$ be a field, and ${\cal C}$ be a $\K$-linear abelian
symmetric monoidal category (we will call such categories {\it tensor}),
which admits infinite sums and products (we refer to [DM] about
all necessary terminology of tensor categories). Then we can do simple linear algebra
in ${\cal C}$, in particular, speak about associative algebras or coassociative coalgebras.
For the rest of the paper, unless we say otherwise, we will assume that either
${\cal C} = Vect_\K^{\Z}$, which is the tensor category of $\Z$-graded
vector spaces $V = \oplus_{n \in \Z} V_n$,
or ${\cal C}= Vect_\K^{\Z/2}$, which is the tensor category
of $\Z/2$-graded vector spaces (then $V=V_0\oplus V_1$), or  ${\cal C} = Vect_\K$, which is
the tensor category of $\K$-vector spaces. There is a natural embedding (faithful tensor functor) from $Vect_\K$ to the other
two categories such that  objects  and morphisms are concentrated in degree zero. The unit object ${\bf 1}$ in all three examples can be identified
with the field $\K$ placed in degree zero.  Since most of our definitions hold for more general tensor categories, we will abuse the notation and
use ${\bf 1}$ and $\K$ interchangeably. 

Associativity morphisms in
all these categories are identity maps, and commutativity morphisms are given by the
Koszul rule of signs: $c ( v_i \otimes v_j ) = ( - 1 )^{ij} v_j \otimes v_i$,
where $v_n$ denotes an element of degree $n$.

We will denote by ${\cal C}^f$ the artinian category of
finite-dimensional objects in ${\cal C}$ (i.e. objects of finite
length). The category $Alg_{{\cal C}^f}$ of unital
finite-dimensional algebras is closed with respect to finite
projective limits. In particular, finite products and finite fiber
products exist in $Alg_{{\cal C}^f}$. One has also the categories
$Coalg_{{\cal C}}$ (resp. $Coalg_{{\cal C}^f}$) of coassociative
counital (resp. coassociative counital finite-dimensional)
coalgebras. In the case ${\cal C} = Vect_\K$ we will also use the
notation $Alg_\K$, $Alg_\K^f$, $Coalg_\K$ and $Coalg_\K^f$ for these
categories. The category $Coalg_{{\cal C}^f} = Alg_{{\cal
C}^f}^{op}$ admits finite inductive limits.

We will need few basic facts about coalgebras.  The proofs are given
in the Appendix.

\begin{thm}\label{coalgebras as functors}
Let $F : Alg_{{\cal C}^f} \to Sets$ be a covariant functor commuting with
finite projective limits. Then it is isomorphic to a functor of the type $A
\mapsto Hom_{Coalg_{{\cal C}}} ( A^{\ast}, B )$ for some counital coalgebra
$B$. Moreover, the category of such functors is equivalent to the category of
counital coalgebras.

\end{thm}

\begin{prp}\label{coalgebra as union}
If $B \in Ob ( Coalg_{{\cal C}} )$, then $B$ is a union of
finite-dimensional counital coalgebras.

\end{prp}

Objects of the category $Coalg_{{\cal C}^f}=Alg_{{\cal C}^f}^{op}$
can be interpreted as ``very thin" non-commutative affine schemes
(cf. with finite schemes in algebraic geometry).
Proposition \ref{coalgebra as union} implies that the category $Coalg_{\cal C}$ is
naturally equivalent to the category of ind-objects in $Coalg_{{\cal
C}^f}$.

For a counital coalgebra $B$ we denote by $Spc(B)$ (and call the ``spectrum" of
the coalgebra $B$) the corresponding functor on the category of
finite-dimensional algebras. A functor {\it isomorphic} to $Spc(B)$ for
some $B$ is called a {\it non-commutative thin scheme}. The category
of non-commutative thin schemes is equivalent to the category of
counital coalgebras. For a non-commutative scheme $X$ we denote by
$B_X$ the corresponding coalgebra. We will call it the coalgebra of
{\it distributions} on $X$. The {\it algebra of functions} on $X$ is
by definition the dual algebra ${\cal O}(X)=B_X^{\ast}:=Hom_{\cal C}(B, {\bf 1})$.

Non-commutative thin schemes form a full monoidal subcategory
$NAff_{{\cal C}}^{th} \subset Ind ( NAff_{{\cal C}} )$ of the
category of non-commutative ind-affine schemes (see Appendix). Tensor product in this monoidal category
corresponds to the tensor product of coalgebras.

\begin{exa} \label{free coalgebra} Let $V \in Ob ( {\cal C} )$. Then $T ( V ) = \oplus_{n \ge 0}
V^{\otimes n}$ carries a structure of counital cofree coalgebra in
${\cal C}$ with the coproduct $\Delta(v_0\otimes...\otimes
v_n)=\sum_{0\le i\le n}(v_0\otimes...\otimes
v_i)\otimes(v_{i+1}\otimes...\otimes v_n)$.  Here $V^{\otimes 0}={\bf 1}=\K$. The corresponding to $T(V)$
non-commutative thin scheme is called the non-commutative formal affine
space $V_{form}$ (or formal neighborhood of zero in $V$).
\end{exa}

\begin{defn} A non-commutative formal  manifold $X$ is a non-commutative thin scheme
isomorphic to some $Spc ( T ( V ) )$ from the example above. The
dimension of $X$ is defined as $dim_\K V$.

\end{defn}

The algebra $\mathcal{O} ( X )$ of functions on a non-commutative
formal manifold $X$ of dimension $n$ is by definition the topological dual to the corresponding colagebra. Hence it isomorphic to the
topological algebra $\K \langle \langle x_1, ..., x_n \rangle
\rangle$ of formal power series in free graded variables $x_1, ...,
x_n$.

Let $X$ be a non-commutative formal manifold, and $pt:\K\to B_X$ a
$\K$-point in $X$.

\begin{defn} \label{formal pointed manifold} The pair $(X,pt)$ is called a non-commutative formal pointed manifold.
If ${\cal C}=Vect_\K^{\Z}$ it will be called non-commutative formal
pointed graded manifold. If ${\cal C}=Vect_\K^{\Z/2}$ it will be
called non-commutative formal pointed supermanifold.

\end{defn}

The notion of tangent space $T_{pt}X$ at the point $pt$ follows  from the agreement than for $Spc(T(V))$ and the natural map $pt: \K\to V^{\otimes 0}=\K$ 
we set $T_{pt}(Spc(T(V)))=V$.  The above-mentioned coordinates $x_1,...,x_n$ can be interpreted as ``affine coordinates\rq\rq{} on $V$. For more details
see Section \ref{A-infty algebras and dg-manifolds} and Appendix.

The following example is a generalization of the Example \ref{free coalgebra} from the quiver with one vertex to an arbitrary quiver.

\begin{exa}\label{generalization to arbitrary quiver} Let $I$ be a set and $B_I = \oplus_{i \in I} \textbf{1}_i$ be the
direct sum of trivial coalgebras. We denote by ${\cal O}(I)$ the
dual topological algebra. The latter can be thought of as the algebra of
functions on a discrete non-commutative thin scheme $I$.

A quiver $Q$ in ${\it C}$ with the set of vertices $I$ is given by a
collection of objects $E_{ij} \in {\cal C}, i, j \in I$ called
spaces of arrows from $i$ to $j$. The coalgebra of $Q$ is the
coalgebra $B_Q$ generated by the ${\cal O}(I) - {\cal
O}(I)$-bimodule $E_Q = \oplus_{i, j \in I} E_{ij}$, i.e. $B_Q \simeq
\oplus_{n \ge 0} \oplus_{i_0, i_1, ..., i_n \in I} E_{i_0 i_1}
\otimes ... \otimes E_{i_{n - 1} i_n} : = \oplus_{n \ge 0} B_Q^n$,
$B_Q^0 : = B_I$. Elements of $B_Q^0$ are called {\it trivial paths}.
Elements of $B_Q^n$ are called paths of the length $n$. Coproduct is
given by the formula

$$\Delta ( e_{i_0 i_1} \otimes ... \otimes
e_{i_{n - 1} i_n} ) = \oplus_{0 \le m \le n} ( e_{i_0 i_1} \otimes
... \otimes e_{i_{m - 1} i_m} ) \otimes ( e_{i_m i_{m + 1}} ...
\otimes ... \otimes e_{i_{n-1} i_{n}} ),$$
where for $m=0$ (resp. $m=n$) we set  $e_{i_{- 1} i_0} = 1_{i_0}$
(resp. $ e_{i_n i_{n + 1}} = 1_{i_n}$).

In particular, $\Delta ( 1_i ) = 1_i \otimes 1_i, i \in I$ and $\Delta (
e_{ij} ) = 1_i \otimes e_{ij} + e_{ij} \otimes 1_j$, where $e_{ij} \in
E_{ij}$, and $1_m \in B_I$ corresponds to the image of $1 \in {\bf 1}$
under the natural embedding into $\oplus_{m \in I} {\bf 1}$.

The coalgebra $B_Q$ has a counit $\varepsilon$ such that $\varepsilon ( {\bf 1}_i )
= {\bf 1}_i$, and $\varepsilon ( x ) = 0$ for $x \in B_Q^n, n \ge 1$.

\end{exa}

\begin{exa}\label{generalized quivers} (Generalized quivers). Here we replace ${\bf 1}_i$ by a unital simple algebra $A_i$
(e.g. $A_i=Mat(n_i,D_i)$, where $D_i$ is a division algebra). Then $E_{ij}$ are
$A_i-mod-A_j$-bimodules. We leave as an exercise to the reader to write down the
coproduct (for that one uses the tensor product of bimodules) and to check that we indeed
obtain a coalgebra.

\end{exa}

\begin{exa} \label{special quiver}Let $I$ be a set. Then the coalgebra $B_I = \oplus_{i \in I}
\textbf{1}_i$ is a direct sum of trivial coalgebras, isomorphic to
the unit object in $\mathcal{C}$. This is a special case of Example \ref{generalization to arbitrary quiver}
Notice that in general $B_Q$ is a ${\cal O}(I)-{\cal
O}(I)$-bimodule.
\end{exa}

\begin{exa} \label{disjoint union of formal neighborhoods}Let $A$ be an associative unital algebra in $\mathcal{C}$. It gives rise to the functor
$F_A : Coalg_{\mathcal{C}^f} \to Sets$ such that $F_A ( B ) =
Hom_{Alg_{\mathcal{C}}} ( A, B^{\ast} )$. This functor describes
finite-dimensional representations of $A$. It commutes with finite
direct limits, hence it is representable by a coalgebra. If $A =
{\cal O} ( X )$ is the algebra of regular functions on the affine
scheme $X$, then in the case of algebraically closed field $\K$ the
coalgebra representing $F_A$ is isomorphic to $\oplus_{x \in X ( \K
)} {\cal O}_{x, X}^{\ast}$, where ${\cal O}_{x, X}^{\ast}$ denotes
the topological dual to the completion of the local ring ${\cal
O}_{x, X}$. If $X$ is smooth of dimension $n$, then each summand is
isomorphic to the topological dual to the algebra of formal power series
$\K [ [t_1, ..., t_n ] ]$. In other words, this coalgebra corresponds to
the disjoint union of formal neighborhoods of all points of $X$.

\end{exa}

\begin{rmk} \label{thin schemes via structure theorems}One can describe non-commutative thin schemes more precisely by using
structure theorems about finite-dimensional algebras in ${\cal C}$. For example,
in the case ${\cal C}=Vect_\K$ any finite-dimensional algebra $A$ is isomorphic to a sum
$A_0\oplus r$, where $A_0$ is a finite sum of matrix algebras $\oplus_iMat(n_i,D_i)$,
$D_i$ are division algebras, and $r$ is the radical. In $\Z$-graded case a similar
decomposition holds, with $A_0$ being a sum of algebras of the type
$End(V_i)\otimes D_i$,where $V_i$ are some graded vector spaces and $D_i$ are
division algebras of degree zero. In $\Z/2$-graded case the description is slightly more
complicated. In particular $A_0$ can contain summands isomorphic to
$(End(V_i)\otimes D_i)\otimes D_{\lambda}$, where $V_i$ and $D_i$ are $\Z/2$-graded
analogs of the above-described objects, and $D_{\lambda}$ is a $1|1$-dimensional
superalgebra isomorphic to $\K[\xi]/(\xi^2=\lambda)$, $deg\,\xi=1, \lambda\in \K^{\ast}/(\K^{\ast})^2$.

\end{rmk}

\subsection{Smooth thin schemes}

Recall that the notion of an ideal has meaning in any abelian tensor category. A $2$-sided
ideal $J$ is called \textit{nilpotent} if the multiplication map $J^{\otimes
n} \to J$ has zero image for a sufficiently large $n$.

\begin{defn}\label{smooth thin schemes}
Counital coalgebra $B$ in a tensor category $\mathcal{C}$ is called
smooth if the corresponding functor $F_B : Alg_{{\cal C}^f} \to Sets, F_B ( A
) = Hom_{Coalg_{{\cal C}}} ( A^{\ast}, B )$ satisfies the following lifting
property: for any $2$-sided nilpotent ideal $J \subset A$ the map $F_B ( A ) \to
F_B ( A / J )$ induced by the natural projection $A \to A / J$ is surjective.
Non-commutative thin scheme $X$ is called smooth if the corresponding counital
coalgebra $B = B_X$ is smooth.
\end{defn}

\begin{prp}\label{quivers give smooth coalgebras}
For any quiver $Q$ in ${\cal C}$ the corresponding coalgebra $B_Q$ is
smooth.

\end{prp}
{\it Proof.} First let us assume that the result holds for all
finite quivers. We remark that if $A$ is finite-dimensional, and $Q$
is an infinite quiver then for any morphism $f : A^{\ast} \to B_Q$
we have: $f ( A^{\ast} )$ belongs to the coalgebra of a finite
sub-quiver of $Q$. Since the lifting property holds for the latter,
the result follows. 

In order vto apply the above considerations we need to prove the Proposition for a
finite quiver $Q$ . For that let us choose a basis $\{e_{ij,\alpha}\}$ of
each space of arrows $E_{ij}$. Then for a finite-dimensional algebra
$A$ the set $F_{B_Q}(A)$ is isomorphic to the set
$\{((\pi_i),x_{ij,\alpha})_{i,j\in I}\}$, where $\pi_i\in
A,\pi_i^2=\pi_i,\pi_i\pi_j=\pi_j\pi_i$, if $i\ne j$,$ \sum_{i\in
I}\pi_i=1_A$, and $x_{ij,\alpha}\in \pi_iA\pi_j$ satisfy the
condition: there exists $N\ge 1$ such that
$x_{i_1j_1,\alpha_1}...x_{i_mj_m,\alpha_m}=0$ for all $m\ge N$. Let
now $J\subset A$ be the nilpotent ideal from the definition of
smooth coalgebra and $(\pi_i^{\prime},x_{ij,\alpha}^{\prime})$ be
elements of $A/J$ satisfying the above constraints. Our goal is to
lift them to $A$. We can lift  them to the projectors $\pi_i$ and
elements $x_{ij,\alpha}$ for $A$ in such a way that  the above
constraints are satisfied except of the last one, which becomes an
inclusion $x_{i_1j_1,\alpha_1}...x_{i_mj_m,\alpha_m}\in J$ for $m\ge
N$. Since $J^n=0$ in $A$ for some $n$ we see that
$x_{i_1j_1,\alpha_1}...x_{i_mj_m,\alpha_m}=0$ in $A$ for $m\ge nN$.
This proves the result. $\blacksquare$

\begin{rmk}\label{remarks about smooth schemes}
a) According to Cuntz and Quillen (see [CQ2]) a non-commutative {\it
algebra} $R$ in $Vect_\K$ is called {\it smooth} if the functor
$Alg_\K \to  Sets, F_R ( A ) = Hom_{Alg_\K} ( R, A )$ satisfies the
lifting property from the Definition \ref{smooth thin schemes} applied to all (not only
finite-dimensional) algebras. We remark that if $R$ is smooth in the
sense of Cuntz and Quillen then the coalgebra $R_{dual}$
representing the functor $Coalg_k^f\to Sets, B\mapsto
Hom_{Alg_k^f}(R,B^{\ast})$  is smooth. One can prove that any smooth
coalgebra in $Vect_\K$ is isomorphic to a coalgebra of a generalized
quiver (see Example \ref{generalized quivers}).

b) Almost all examples of non-commutative smooth thin schemes
considered in this paper are formal pointed manifolds, i.e. they are
isomorphic to $Spc(T(V))$ for some $V\in Ob({\cal C})$. It is
natural to try to ``globalize" our results to the case of
non-commutative ``smooth" schemes $X$ which satisfy the property
that the completion of $X$ at a ``commutative" point gives rise to a
formal pointed manifold in our sense. An example of the space of
maps is considered in the next subsection.

c) The tensor product of non-commutative smooth thin
schemes is typically non-smooth, since it corresponds to the {\it tensor product} of
coalgebras. Notice that the tensor product of colagebras  is not a categorical product.
\end{rmk}

Let now $x$ be a $\K$-point of a non-commutative smooth thin scheme
$X$. By definition $x$ is a homomorphism of counital coalgebras
$x:\K\to B_X$ (here $\K={\bf 1}$ is the trivial coalgebra
corresponding to the unit object). The completion $\widehat{X}_x$ of
$X$ at $x$ is a formal pointed manifold which can be described such as follows. As a functor $F_{\widehat{X}_x}:Alg_{\cal
C}^f\to Sets$ it assigns to a finite-dimensional algebra $A$ the set
of such homomorphisms of counital colagebras $f:A^{\ast}\to B_X$
which are compositions $A^{\ast}\to A_1^{\ast}\to B_X$, where
$A_1^{\ast}\subset B_X$ is a conilpotent extension of $x$ (i.e.
$A_1$ is a finite-dimensional unital nilpotent algebra such that the
natural embedding $\K\to A_1^{\ast}\to B_X$ coincides with $x:\K\to
B_X$).

Description of the coalgebra $B_{\widehat{X}_x}$ is given in the following Proposition.

\begin{prp} \label{formal neighborhood}The formal neighborhood $\widehat{X}_{x}$  corresponds to
the counital sub-coalgebra $B_{\widehat{X}_x} \subset B_X$ which is the preimage under the
natural projection $B_X \to B_X / x ( \K )$ of the sub-coalgebra
consisting of conilpotent elements in the non-counital coalgebra $B / x ( \K
)$. Moreover, $\widehat{X}_{x}$ is universal for all morphisms from nilpotent
extensions of $x$ to $X$.

\end{prp}

We discuss in Appendix a more general construction of the completion along a non-commutative
thin subscheme.

We leave as an exercise to the reader to prove the following result.

\begin{prp} \label{formal neighborhood of a point} Let $Q$ be a quiver and $pt_i\in X=X_{B_Q}$ corresponds to a vertex
$i\in I$. Then the formal neighborhood $\widehat{X}_{pt_i}$ is a
formal pointed manifold corresponding to the tensor coalgebra
$T(E_{ii})=\oplus_{n\ge 0}E_{ii}^{\otimes n}$, where $E_{ii}$ is the
space of loops at $i$.

\end{prp}

\subsection{Inner Hom}\label{inner Hom}

Let $X, Y$ be non-commutative thin schemes, and $B_X, B_Y$ the
corresponding coalgebras.

\begin{thm}\label{representability of Maps}
The functor $Alg_{{\cal C}^f} \to Sets$ such that
$$ A \mapsto Hom_{Coalg_{\cal C}} ( A^{\ast} \otimes B_X, B_Y ) $$
is representable. 

The corresponding non-commutative thin scheme will be
denoted by $Maps ( X, Y )$.
\end{thm}
{\it Proof.} It is easy to see that the functor under
consideration commutes with finite projective limits. Hence it is of the type
$A \mapsto Hom_{Coalg_{\cal C}} ( A^{\ast}, B )$, where $B$ is a counital coalgebra (Theorem
\ref{coalgebras as functors}). The corresponding non-commutative thin scheme is the desired $Maps ( X, Y
)$. $\blacksquare$

It follows from the definition that $Maps ( X, Y ) = \ihom ( X, Y )$, where
the inner Hom is taken in the symmetric monoidal category of non-commutative
thin schemes. By definition $\ihom ( X, Y )$ is a non-commutative thin scheme,
which satisfies the following functorial isomorphism for any $Z \in Ob (
NAff_{{\cal C}}^{th} )$:
$$ Hom_{NAff_{{\cal C}}^{th}} ( Z, \ihom ( X, Y ) ) \simeq
   Hom_{NAff_{{\cal C}}^{th}} ( Z \otimes X, Y ) . $$
Notice that the monoidal category $NAff_{{\cal C}}$ of all
non-commutative affine schemes does not have inner $Hom's$ even in
the case ${\cal C} = Vect_\K$. If ${\cal C} = Vect_\K$ then one can
define $\ihom ( X, Y )$ for $X = Spec ( A )$, where $A$ is a {\it
finite-dimensional} unital algebra and $Y$ is arbitrary. The
situation is similar to the case of  the conventional (i.e. \lq{}commutative\rq{}) algebraic
geometry, where one can define an affine scheme of maps from a
scheme of finite length to an arbitrary affine scheme. On the other
hand, one can show that the category of non-commutative ind-affine
schemes admits inner Hom's (the corresponding result for commutative
ind-affine schemes is known.

\begin{rmk}\label{Maps and quivers}
The non-commutative thin scheme $Maps(X,Y)$ gives rise to a quiver,
such that its vertices are $\K$-points of $Maps(X,Y)$. In other words,
vertices correspond to homomorphisms $B_X\to B_Y$ of the coalgebras
of distributions. Taking the completion at  a $\K$-point we obtain
a formal pointed manifold. More generally, one can take a completion along
a subscheme of $\K$-points, thus arriving to a non-commutative
formal manifold with a marked closed subscheme (rather than one point). This construction can be used
for the description of the $\A$-structure on $\A$-functors.
We also remark that the space of arrows $E_{ij}$ of a quiver is an example of the geometric
notion of bitangent space at a pair of $\K$-points $i,j$.
We plan to discuss this topic in another  paper.
\end{rmk}

\begin{exa} Let $Q_1=\{i_1\}$ and $Q_2=\{i_2\}$ be quivers with one vertex
such that $E_{i_1i_1}=V_1, E_{i_2i_2}=V_2$, $dim\,V_i<\infty, i=1,2$.
Then $B_{Q_i}=T(V_i), i=1,2$ and $Maps(X_{B_{Q_1}},X_{B_{Q_2}})$ corresponds
to the quiver $Q$ such that the set of vertices $I_Q=Hom_{Coalg_{\cal C}}(B_{Q_1},B_{Q_2})=
\prod_{n\ge 1}\ihom(V_1^{\otimes n},V_2)$ and
for any two vertices $f,g\in I_Q$ the space of arrows is isomorphic
to $E_{f,g}=\prod_{n\ge 0}\ihom(V_1^{\otimes n},V_2)$.

\end{exa}

\begin{defn} Homomorphism $f:B_1\to B_2$ of counital coalgebras is
called a minimal conilpotent extension if it is an inclusion and the
induced coproduct on the non-counital coalgebra $B_2/f(B_1)$ is
trivial.

\end{defn}

Composition of minimal conilpotent extensions is simply called a
conilpotent extension. Definition \ref{smooth thin schemes} can be reformulated in terms of
finite-dimensional coalgebras. Namely, the coalgebra $B$ is smooth if the
functor $C\mapsto Hom_{Coalg_{\cal C}}(C,B)$ satisfies the lifting
property with respect to  conilpotent extensions of
finite-dimensional counital coalgebras. The following result
shows that we can drop the condition of finite-dimensionality.

\begin{prp} \label{smoothness without finite-dimensionality}If $B$ is a smooth coalgebra then the functor
$Coalg_{\cal C}\to Sets$ such that $C\mapsto Hom_{Coalg_{\cal
C}}(C,B)$ satisfies the lifting property for conilpotent extensions.

\end{prp}

{\it Proof.} Let $f:B_1\to B_2$ be a conilpotent extension, and
$g:B_1\to B$ an arbitrary homomorphism of counital coalgebras. It
can be thought of as homomorphism of $f(B_1)\to B$. We need to show
that $g$ can be extended to $B_2$. Let us consider the set of pairs
$(C,g_C)$ such $f(B_1)\subset C\subset B_2$ and $g_C: C\to B$
defines  an extension of counital coalgebras, which coincides with
$g$ on $f(B_1)$. We apply Zorn lemma to the partially ordered set of
such pairs and see that there exists a maximal element
$(B_{max},g_{max})$ in this set. We claim that $B_{max}=B_2$.
Indeed, let $x\in B_2\setminus B_{max}$. Then there exists a
finite-dimensional coalgebra $B_x\subset B_2$ which contains $x$.
Clearly $B_x$ is a conilpotent extension of $f(B_1)\cap B_x$. Since
$B$ is smooth we can extend $g_{max}: f(B_1)\cap B_x\to B$ to
$g_x:B_x\to B$ and,finally to $g_{x,max}: B_x+B_{max}\to B$. This
contradicts to maximality of $(B_{max},g_{max})$. Proposition is
proved. $\blacksquare$

\begin{prp}\label{smothness of Maps}
If $X, Y$ are non-commutative thin schemes and $Y$ is smooth then $Maps ( X,
Y )$ is also smooth.
\end{prp}

{\it Proof.} Let $A\to A/J$ be a nilpotent extension of
finite-dimensional unital algebras. Then $(A/J)^{\ast}\otimes B_X\to
A^{\ast}\otimes B_X$ is a conilpotent extension of counital
coalgebras. Since $B_Y$ is smooth then the  Proposition \ref{smoothness without finite-dimensionality}
implies that the induced  map $Hom_{Coalg_{\cal C}}(A^{\ast}\otimes
B_X,B_Y)\to Hom_{Coalg_{\cal C}}((A/J)^{\ast}\otimes B_X,B_Y)$ is
surjective. This concludes the proof. $\blacksquare$

Let us consider the case when $(X,pt_X)$ and $(Y,pt_Y)$ are
non-commutative formal pointed manifolds in the category ${\cal
C}=Vect_\K^{\Z}$. One can describe ``in coordinates" the
non-commutative formal pointed manifold, which is the formal
neighborhood of a $\K$-point of $Maps(X,Y)$. Namely, let $X = Spc ( B
)$ and $Y = Spc ( C )$, and let $f \in Hom_{NAff_{{\cal C}}^{th}} (
X, Y )$ be a morphism preserving marked points. Then $f$ gives rise
to a $\K$-point of $Z = Maps ( X, Y )$. Since ${\cal O} ( X )$ and
${\cal O} ( Y )$ are isomorphic to the topological algebras of
formal power series in  free graded variables, we can choose  sets
of free topological generators $(x_i )_{i\in I}$ and $( y_j )_{j\in
J}$ for these algebras. Then we can write for the corresponding
homomorphism of algebras $f^{\ast} : {\cal O} ( Y ) \to {\cal O} ( X
)$:
$$ f^{\ast} ( y_j ) = \sum_M c_{j, M}^0 x^M, $$
where $c_{j, M}^0 \in \K$ and $M=(i_1,...,i_n), i_{s}\in I$ is a
non-commutative multi-index (all the coefficients depend on $f$,
hence a better notation should be $c_{j, M}^{f, 0}$). Notice that
for $M = 0$ one gets $c_{j, 0}^0 = 0$ since $f$ is a morphism of
pointed schemes. Then we can consider an ``infinitesimal
deformation" $f_{def}$ of $f$

$$ f^{\ast}_{def} ( y_j ) = \sum_M (c_{j, M}^0+\delta c_{j, M}^0) x^M, $$
where $\delta c_{j, M}^0$ are {\it new variables} commuting with all $x_i$.
Then $\delta c_{j, M}^0$ can be thought of as coordinates in the formal
neighborhood of $f$.
More pedantically it can be spelled out such as follows.
Let $A = \K \oplus m$ be a finite-dimensional graded unital algebra, where $m$
is a graded nilpotent ideal of $A$. Then an $A$-point of the formal
neighborhood $U_f$ of $ f$ is a morphism $\phi \in
Hom_{NAff_{{\cal C}}^{th}} ( Spec ( A ) \otimes X, Y )$, such that it
reduces to $f$ modulo the nilpotent ideal $m$. We have for the corresponding
homomorphism of algebras:
$$ \phi^{\ast} ( y_j ) = \sum_M c_{j, M} x^M, $$
where $M$ is a non-commutative multi-index, $c_{j, M} \in A$, and
$c_{j, M} \mapsto c_{j, M}^0$ under the natural homomorphism $A \to
\K = A / m$. In particular $c_{j, 0} \in m$. We can treat
coefficients $c_{j, M}$ as $A$-points of the formal neighborhood
$U_f$ of $ f  \in Maps ( X, Y )$.

\begin{rmk} The above definitions  play an important role in the case when the non-commutative smooth thin scheme
$Spc ( B_Q )$ is assigned to a (small) $\A$-category. Then
the non-commutative smooth thin scheme
$Maps ( Spc ( B_{Q_1} ), Spc (B_{Q_2} ) )$ can  be used for the
description of the category of $\A$-functors between $\A$-categories, and the formal
neighborhood of a point in the  space
$Maps ( Spc ( B_{Q_1} ), Spc (B_{Q_2} ) )$   corresponds to natural
transformations between $\A$-functors.
\end{rmk}

\section{$\A$-algebras}\label{A-infty algebras}

\subsection{Main definitions}\label{A-infty algebras and dg-manifolds}

From now on assume that ${\cal C}=Vect_\K^{\Z}$ unless we say otherwise.
If $X$ is a thin scheme then a vector field on $X$ is, by definition, a derivation of the coalgebra $B_X$.
Since derivations can have different degrees, the vector fields form a {\it graded} Lie algebra $Vect(X)$.

\begin{defn} \label{nc dg scheme}A non-commutative thin differential-graded (dg for short) scheme is a
pair $( X, d )$ where $X$ is a non-commutative thin scheme, and $d$ is a
vector field on $X$ of degree $+ 1$ such that $[ d, d ] = 0$.

\end{defn}
We will call the vector field $d$
{\it homological vector field}.

Let $X$ be a formal pointed manifold and $x_0$ be its unique $\K$-point. Such a
point corresponds to a homomorphism of counital coalgebras $\K \to B_X$. We say
that the vector field $d$ vanishes at $x_0$ if the corresponding derivation
kills the image of $\K$.

\begin{defn} \label{formal pointed dg manifold}A non-commutative formal pointed dg-manifold is a pair $((X,x_0),d)$
such that $(X,x_0)$ is a non-commutative formal pointed graded manifold, and  $d=d_X$ is
a homological vector field  on $X$ such
that $d|_{ x_0}  = 0$.
\end{defn}

Homological vector field $d$ has an infinite Taylor decomposition at $x_0$. More
precisely, let $T_{x_0} X$ be the tangent space at $x_0$. It is canonically isomorphic to the graded vector space of primitive elements
of the coalgebra $B_X$, i.e. the set of $a\in B_X$ such that
$\Delta(a)=1\otimes a+a\otimes 1$ where $1\in B_X$ is the image
of $1\in \K$ under the homomorphism of coalgebras $x_0:\K\to B_X$
(see Appendix for the general definition
of the tangent space). Then
$d:=d_X$ gives rise to a  (non-canonically defined) collection of linear maps $d_X^{(n)}:=m_n : T_{x_0} X^{\otimes n} \to
T_{x_0} X [1], n \ge 1$ called {\it Taylor coefficients of $d$} which satisfy a system of quadratic relations
arising from the condition $[ d, d ] = 0$. Indeed, our
non-commutative formal pointed manifold is isomorphic to the formal
neighborhood of zero in $T_{x_0} X$, hence the corresponding non-commutative
thin scheme is isomorphic to the cofree tensor coalgebra $T ( T_{x_0} X )$
generated by $T_{x_0} X$. Homological vector field $d$ is a derivation of
a cofree coalgebra,
hence it is uniquely determined by a sequence
of linear maps $m_n$.

\begin{defn}\label{non-unital A-infty algebra} Non-unital $\A$-algebra over $\K$ is given by a non-commutative formal
pointed dg-manifold $( X, x_0, d )$ together with an isomorphism of counital
coalgebras $B_X \simeq T ( T_{x_0} X )$.

\end{defn}

As we already pointed out, a choice of  isomorphism with the tensor coalgebra generated by the tangent
space is a non-commutative analog of a choice of affine structure in the
formal neighborhood of $x_0$.

From the above definitions one can recover the conventional definition of an $\A$-algebra.
We present it below for convenience of the reader.

\begin{defn} \label{standard definition of A-infty algebra}A structure of an $\A$-algebra on $V \in Ob ( Vect_\K^{\Z} )$ is
given by a derivation $d$ of degree $+ 1$ of the non-counital
cofree coalgebra $T_+ ( V [ 1
] ) = \oplus_{n \ge 1} V^{\otimes n}$ such that $[ d, d ] = 0$ in the
differential-graded Lie algebra of coalgebra derivations.

\end{defn}

Traditionally the Taylor coefficients of $d=m_1+m_2+...$
are called (higher) multiplications for $V$.
The pair $( V, m_1 )$ is a complex of $\K$-vector spaces called the
{\it tangent complex}. If $X=Spc(T(V))$ then $V[1]=T_{0}X$ and $m_1=d_X^{(1)}$
is the first Taylor coefficient of the homological vector field $d_X$. The tangent cohomology groups $H^i(V,m_1)$ will be
denoted by $H^i ( V )$. Clearly $H^{\bullet} ( V ) = \oplus_{i\in \Z}
H^i ( V )$ is an associative (non-unital) algebra with the product induced by $m_2$.

An important class of
$\A$-algebras consists of {\it unital} (or strictly unital)
and  {\it weakly
unital} (or homologically unital) ones. We are going to discuss the definition
and the geometric meaning of unitality later.

Homomorphism of $\A$-algebras can be described geometrically as a morphism of the corresponding non-commutative formal pointed dg-manifolds. 
In the algebraic form one recovers the following traditional definition.

\begin{defn} \label{morphism of A-infty algebras}A homomorphism of non-unital $\A$-algebras ($\A$-morphism for short)
$( V, d_V ) \to ( W, d_W )$ is a homomorphism of differential-graded
coalgebras $T_+ ( V [ 1 ] ) \to T_+ ( W [ 1 ] )$.
\end{defn}

A homomorphism $f$ of non-unital $\A$-algebras is determined by its
Taylor coefficients $f_n : V^{\otimes n} \to W [ 1 - n ], n \ge 1$ satisfying
the system of equations

$\sum_{1 \le l_1 < ..., < l_i = n} ( - 1 )^{\gamma_i} m_i^W ( f_{l_1} ( a_1,
..., a_{l_1} ),$\\
$f_{l_2 - l_1} ( a_{l_1 + 1}, ..., a_{l_2} ), ..., f_{n - l_{i - 1}} ( a_{n -
l_{i - 1} + 1}, ..., a_n ) ) =$

$\sum_{s + r = n + 1} \sum_{1 \le j \le s} ( - 1 )^{\epsilon_s} f_s ( a_1,
..., a_{j - 1}, m_r^V ( a_j, ..., a_{j + r - 1} ), a_{j + r}, ..., a_n ) .$

Here $\epsilon_s = r \sum_{1 \le p \le j - 1} deg ( a_p ) + j - 1 + r ( s - j
)$, $\gamma_i = \sum_{1 \le p \le i - 1} ( i - p ) ( l_p - l_{p - 1} - 1 ) +
\sum_{1 \le p \le i - 1} \nu ( l_p ) \sum_{l_{p - 1} + 1 \le q \le l_p} deg (
a_q )$, where we use the notation $\nu ( l_p ) = \sum_{p + 1 \le m \le i} ( 1
- l_m + l_{m - 1} )$, and set $l_0 = 0$.

\begin{rmk} \label{case of supermanifolds}All the above definitions and results are valid
for $\Z/2$-graded $\A$-algebras as well. In this case we consider
formal manifolds in the category $Vect_\K^{\Z/2}$ of $\Z/2$-graded
vector spaces. We will use the correspodning results without
further comments. In this case one traditionally denotes
by $\Pi A$ the $\Z/2$-graded vector space $A[1]$.

\end{rmk}

\subsection{Minimal models of $\A$-algebras}\label{minimal models}

One can do simple differential geometry in the symmetric monoidal category of non-commutative
formal pointed dg-manifolds. New phenomenon
is the possibility to define some structures
up to a quasi-isomorphism.

 \begin{defn} \label{qis}Let $f : ( X, d_X, x_0 ) \to ( Y, d_Y, y_0 )$ be a morphism of
non-commutative formal pointed dg-manifolds. We say that $f$ is a
quasi-isomorphism if the induced morphism of the tangent complexes $f_1 : (
T_{x_0}X, d_X^{( 1 )} ) \to ( T_{y_0}Y, d_Y^{( 1 )} )$ is a quasi-isomorphism.
We will use the same terminology for the corresponding $\A$-algebras.
\end{defn}

\begin{defn} \label{minimal A-infty algebra}An $\A$-algebra $A$ (or the corresponding non-commutative formal
pointed dg-manifold) is called minimal if $m_1 = 0$. It is called contractible
if $m_n = 0$ for all $n \ge 2$ and $H^{\bullet} ( A, m_1 ) = 0$.
\end{defn}

The notion of minimality is coordinate independent, while the
notion of contractibility is not.

It is  not difficult to prove that any $\A$-algebra $A$ has a \textit{minimal model
$M_A$}.  Here we call  $M_A$ a minimal model if  $M_A$ is minimal and there exists a quasi-isomorphism $M_A \to A$.
The proof of the existence of a minimal model is similar to the one from [Ko1], [KoSo2].
The minimal model is unique up to an $\A$-isomorphism.
We are going to use the same
terminology when talking about  non-commutative formal pointed dg-manifolds. In particular, in geometric
language a non-commutative formal pointed dg-manifold $X$ is isomorphic to a {\it categorical} product
(the latter corresponds  to the completed  free product of algebras of functions)
$X_m \times X_{lc}$, where $X_m$ is minimal and $X_{lc}$ is linear
contractible. The above-mentioned quasi-isomorphism corresponds in the geometric language to the
projection $X \to X_m$. In what follows we will freely move from algebraic to geometric language and back at our convenience.

The following result (homological inverse function theorem)
can be easily deduced from the above product
decomposition.

\begin{prp}\label{inverse function theorem} If $f : A \to B$ is a quasi-isomorphism of $\A$-algebras then there is a (non-canonical) quasi-isomorphism $g : B \to A$ such that $fg$ and $gf$ induce identity maps
on zero cohomologies $H^0 ( B )$ and $H^0 ( A )$ respectively.

\end{prp}

\subsection{Centralizer of an $\A$-morphism}\label{centralizer}

Let $A$ and $B$ be two $\A$-algebras, and $( X, d_X, x_0 )$ and $( Y, d_Y, y_0 )$ be the
corresponding non-commutative formal pointed dg-manifolds.
Let $f : A \to B$ be a morphism of $\A$-algebras. Then the corresponding
$\K$-point $f \in Maps ( Spc ( A ), Spc ( B ) )$ gives rise to the formal
pointed manifold $U_{ f }=\widehat{Maps}(X,Y)_f$
(completion at the point $f$).
Functoriality of the construction of $Maps ( X, Y
)$ gives rise to a homomorphism of graded Lie algebras of vector fields $Vect ( X ) \oplus
Vect ( Y ) \to Vect ( Maps ( X, Y ) )$. Since $[ d_X, d_Y ] = 0$ on $X \otimes
Y$, we have a well-defined homological vector field $d_Z$ on $Z = Maps ( X, Y )$. It
corresponds to $d_X\otimes 1_Y-1_X\otimes d_Y$ under the above homomorphism. It is easy to see that
$d_Z|_{f}=0$ and in fact morphisms $f:A\to B$ of $\A$-algebras
are exactly zeros of $d_Z$.
We are going to describe below the $\A$-algebra $Centr ( f )$ (centralizer of $f$) which
corresponds to the formal neighborhood $U_f$ of the point $ f \in Maps (
X, Y )$.
We can write (see  Section \ref{centralizer} for the notation)
$$ c_{j, M} = c_{j, M}^0 + r_{j, M}, $$
where $c_{j, M}^0 \in \K$ and $r_{j, M}$ are formal non-commutative coordinates in the
neighborhood of $f$. Then the $\A$-algebra $Centr ( f )$ has a 
basis $(r_{j, M})_{j,M}$ and the $\A$-structure
is defined by the restriction of the homological vector $d_Z$ to $U_f$.

As a $\Z$-graded vector space $Centr ( f ) = \prod_{n \ge 0} Hom_{Vect^{\Z}_\K}
( A^{\otimes n}, B ) [ - n ]$. Let $\phi_1, ..., \phi_n \in Centr ( f )$, and
$a_1, ..., a_N \in A$. Then we have $m_n ( \phi_1, ..., \phi_n ) ( a_1, ..., a_N ) =
I + R$. Here $I$ corresponds to the term $=1_X\otimes d_Y$ and is given by the following expression:

\vspace{3mm}


\includegraphics{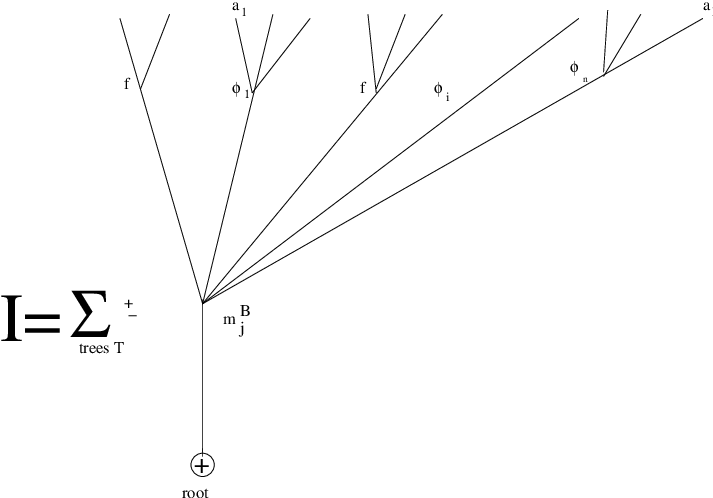}

{\vspace{3mm}}

Similarly $R$ corresponds to the term $d_X\otimes 1_Y$ and it is described by the following picture:

\vspace{3mm}


\includegraphics{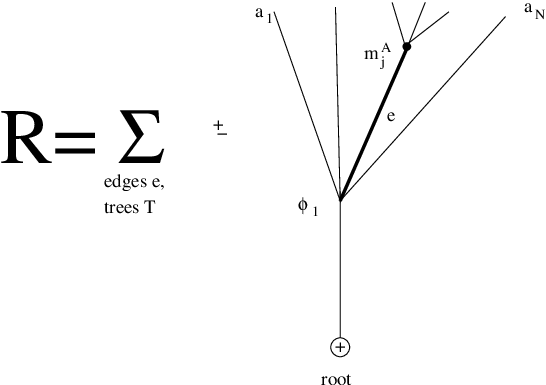}

{\vspace{3mm}}

Here are some comments on the figure describing $I$:

1) We partition a sequence $(a_1, ..., a_N)$ into $l \ge n$ non-empty subsequences.

2) We mark $n$ of these subsequences counting from the left (the set can be empty).

3) We apply multilinear map $\phi_i, 1 \le i \le n$ to the $i$th marked group of elements
$a_l$.

4) We apply Taylor coefficients of $f$ to the remaining subsequences.

Notice that the term $R$ appear only for $m_1$ (i.e. $n=1$).
For all subsequences we have $n\ge 1$.

From geometric point of view the term $I$ corresponds to  the
vector field $d_Y$, while the term $R$ corresponds to the  vector
field $d_X$. 

\begin{prp}\label{commutator of vector fields}
Let $d_{Centr ( f )}$ be the derivation corresponding to the image of $d_X
\oplus d_Y$ in $Maps ( X, Y )$.
Then one  has:
 $$[ d_{Centr ( f )}, d_{Centr ( f )} ] = 0.$$

\end{prp}
{\it Proof.} Straightforward computation .
$\blacksquare$

\begin{rmk}\label{natural transformations}
The $\A$-algebra $Centr(f)$ and its generalization to the
case of $\A$-categories  provide geometric
description of the notion of natural transformation in the $\A$-case (see [Lyu],[LyuOv]
for a pure algebraic approach to this notion).
\end{rmk}

\section{Non-commutative dg-line ${\bf L}$ and weak unit}\label{weak unit}

\subsection{Main definition}\label{main definition of weak unit}

\begin{defn}\label{strictly unital algebras}
An $\A$-algebra is called {\it unital} (or strictly unital) if there exists an element $1 \in V$ of degree zero,
such that $m_2 ( 1, v ) = m_2 ( v, 1 )$ and $m_n ( v_1, ..., 1, ..., v_n ) =
0$ for all $n \ne 2$ and $v, v_1, ..., v_n \in V$. It is called {\it weakly
unital} (or homologically unital) if  the  graded
associative unital algebra $H^{\bullet} ( V )$ has a unit $1 \in H^0 ( V )$.
\end{defn}

The notion of strict unit depends on a choice of
affine coordinates on $Spc(T(V))$,
while the notion of weak unit is ``coordinate free".
Moreover, one can show that a weakly unital $\A$-algebra becomes strictly
unital after an appropriate change of coordinates.

The category of
unital or weakly unital $\A$-algebras are defined in the natural way by the requirement that
morphisms preserve the unit (or weak unit) structure.

In this section we are going to discuss a non-commutative
dg-version of  the odd
$1$-dimensional supervector space ${\bf A}^{0|1}$ and its relationship to weakly unital
$\A$-algebras. The results are valid for both  $\Z$-graded and $\Z/2$-graded
$\A$-algebras.

\begin{defn}\label{dg line}
Non-commutative formal dg-line ${\bf L}$ is a non-commutative formal
pointed dg-manifold corresponding to the one-dimensional $\A$-algebra
$A\simeq \K$ such that $m_2=id, m_{n\ne 2}=0$.

\end{defn}
The algebra of functions ${\cal O}({\bf L})$ is isomorphic to the
topological algebra of formal series
$\K \langle \langle\xi \rangle \rangle$, where $deg\,\xi=1$. The differential
is given by $\partial(\xi)=\xi^2$.

\subsection{Adding a weak unit}\label{adding weak unit}

Let $(X,d_X,x_0)$ be a non-commutative formal pointed dg-manifold
correspodning to a non-unital $\A$-algebra $A$. We would like to describe
geometrically the procedure of adding a weak unit to $A$.

Let us consider the non-commutative formal pointed graded manifold
$X_1={\bf L}\times X$  corresponding to the free product of the corresponding coalgebras, i.e. to 
$B_{\bf L}\ast B_X$. Clearly  one can lift vector fields $d_X$ and
$d_{\bf L}:=\partial/\partial\, \xi$ to $X_1$.

\begin{lmm} \label{commutator of vector field with itself}The vector field
$$d:=d_{X_1}=d_X+ad(\xi)-\xi^2\partial/\partial\,\xi$$
satisfies the condition $[d,d]=0$.

\end{lmm}

{\it Proof.} Straightforward check. $\blacksquare$

\begin{prp}\label{existence of strict unit} The $\A$-algebra $A_1$ has a strict unit.

\end{prp}

{\it Proof.} It follows from the formulas  that $\xi$ appears
in the expansion of $d_{X_1}$ in quadratic expressions only. Let $A_1$ be the $\A$-algebras
corresponding to $X_1$ and  $1\in T_{pt}X_1=A_1[1]$ be the element of $A_1[1]$
dual to $\xi$ (it corresponds to the tangent vector $\partial/\partial\,\xi$).
Thus we see that $m_2^{A_1}(1,a)=m_2^{A_1}(a,1)=a, m_2^{A_1}(1,1)=1$
for any $a\in A$ and $m_n^{A_1}(a_1,...,1,...,a_n)=0$
for all $n\ge 2, a_1,...,a_n\in A$. This completes the proof. $\blacksquare$

Notice that we have a canonical morphism of non-commutative formal pointed
dg-manifolds $e:X\to X_1$ such that
$e^{\ast}|_X=id, e^{\ast}(\xi)=0$.

\begin{defn} \label{weak unit as morphism}Weak unit in $X$ is given by a
morphism of non-commutative formal pointed dg-manifolds
$p:X_1\to X$ such that $p\circ e=id$.

\end{defn}
It follows from the definition that if $X$ has a weak unit then the
associative algebra $H^{\bullet}(A,m_1^{A})$ is unital. Hence our
geometric definition agrees with the pure algebraic one
(explicit algebraic description of the notion of weak unit
can be found e.g. in [FOOO], Section 20 \footnote {V. Lyubashenko has informed
us that the equivalence of two descriptions  also follows from his results with Yu. Bespalov and O. Manzyuk.}).

\section{Modules and bimodules}\label{modules and bimodules}

\subsection{Modules and vector bundles}\label{modules and vector bundles}

Recall that a topological vector space is called
linearly compact if it is a projective limit of finite-dimensional
vector spaces. The duality functor $V\mapsto V^{\ast}$ establishes
an anti-equivalence between the category of vector spaces equipped
with the discrete topology and the category of linearly compact
vector spaces. All that can be extended in the obvious way to the
category of graded vector spaces.

Let $X$ be a non-commutative thin scheme in $Vect_\K^{\Z}$.

\begin{defn}\label{vector bundles}
Linearly compact vector bundle ${\cal E}$ over $X$ is given by a linearly compact topologically free ${\cal O}(X)$-module $\Gamma({\cal E})$, where
${\cal O}(X)$ is the algebra of function on $X$.
Module  $\Gamma({\cal E})$ is called the module of {\it sections } of the linearly compact vector bundle ${\cal E}$.

\end{defn}

Suppose that $(X,x_0)$ is formal pointed graded manifold.
By definition, {\it the fiber}
of ${\cal E}$ over $x_0$ is given by the quotient space
${\cal E}_{x_0}=\Gamma({\cal E})/\overline{m_{x_0}\Gamma({\cal E})}$ where $m_{x_0}\subset {\cal O}(X)$ is the
$2$-sided maximal ideal of functions vanishing at $x_0$, and the bar means the closure.

\begin{defn}\label{dg vector bundle}

A dg-vector bundle over
a formal pointed dg-manifold $(X,d_X,x_0)$
is given by a linearly compact vector bundle
${\cal E}$ over $(X,x_0)$ such that the corresponding module $\Gamma({\cal E})$
carries a differential $d_{\cal E}: \Gamma({\cal E})\to \Gamma({\cal E})[1] , d_{\cal E}^2=0$ so that $(\Gamma({\cal E}),d_{\cal E})$ becomes a dg-module over the dg-algebra
$({\cal O}(X),d_X)$ and $d_{\cal E}$ vanishes on ${\cal E}_{x_0}$.
\end{defn}

\begin{defn} \label{modules} Let $A$ be a non-unital $\A$-algebra. A left $A$-module
$M$ is given by a dg-bundle $E$
over the formal pointed dg-manifold $X=Spc(T(A[1]))$ together with
an isomorphism of vector bundles $\Gamma({\cal E})\simeq{\cal O}(X)\widehat{\otimes} M^{\ast}$ called a trivialization of ${\cal E}$.

\end{defn}
Here $\widehat{\otimes}$ denotes the topologically completed tensor product.
Passing to dual spaces we obtain the following algebraic definition.

\begin{defn} \label{equivalent definition of modules} Let $A$ be an $\A$-algebra and $M$ be a $\Z$-graded vector space. A structure of a left $\A$-module on $M$ over $A$
(or simply a structure of a left $A$-module on $M$) is given by a
differential $d_M$ of degree $+1$ on $T(A[1])\otimes M$ which makes it into
a dg-comodule over the dg-coalgebra $T(A[1])$.

\end{defn}

The notion of {\it  right} $\A$-module
is similar. Right $A$-module is the same as  left $A^{op}$-module.
Here $A^{op}$ is the {\it opposite} $\A$-algebra, which
coincides with $A$ as a $\Z$-graded vector space,
and for the higher multiplications one has
$m_n^{op}(a_1,...,a_n)=(-1)^{n(n-1)/2}m_n(a_n,...,a_1)$.
The $\A$-algebra $A$ carries the natural structures of
the left and right $A$-modules.
If we simply say ``$A$-module" it will always mean ``left $A$-module".

Taking the Taylor series of $d_M$ we obtain a collection of $\K$-linear maps (higher action morphisms) for any $n \ge 1$
$$ m_n^M : A ^{\otimes (n-1)} \otimes M \to M [2-n], $$
satisfying the compatibility conditions which
can be written in exactly the same form as compatibility conditions for the higher products $m_n^A$ (see e.g. [Ke1]).
All those conditions can be derived from just one property
that the cofree $T_+(A[1])$-comodule
$T_+(A[1],M)=\oplus_{n\ge 0}A[1]^{\otimes n}\otimes M$
carries a derivation $m^M=(m_n^M)_{n\ge 0}$ such that
$[m^M,m^M]=0$. In particular $(M,m_1^M)$ is a complex of
vector spaces.

\begin{defn} \label{weakly unital module}Let $A$ be a weakly unital $\A$-algebra. An $A$-module
$M$ is called weakly unital if the cohomology $H^{\bullet}(M,m_1^M)$
is a unital $H^{\bullet}(A)$-module.

\end{defn}

It is easy to see that left $\A$-modules over $A$ form a dg-category $A-mod$ with morphisms
being homomorphisms of the corresponding comodules.
As a graded vector space
$$ Hom_{A - mod} ( M, N ) = \oplus_{n \ge 0} \ihom_{Vect_k^{\Z}} ( A [ 1 ]^{\otimes n} \otimes
   M, N ) . $$
It easy to see that $Hom_{A - mod} ( M, N )$
is a complex.

If $M$ is a right $A$-module and $N$ is a left $A$-module
then one has a naturally defined structure of a complex on
$M\otimes_AN:=\oplus_{n\ge 0}M\otimes A[1]^{\otimes n}\otimes N$.
The differential is given by the formula:

$$d(x\otimes a_1\otimes...\otimes a_n\otimes y)=\sum\pm m_i^M(x\otimes a_1\otimes...\otimes a_i)\otimes a_{i+1}\otimes...\otimes a_n\otimes y)+$$
$$\sum\pm x\otimes a_1\otimes...\otimes a_{i-1}\otimes m_k^A(a_i\otimes...\otimes a_{i+k-1})\otimes a_{i+k}\otimes...\otimes a_n\otimes y +$$
$$\sum\pm x\otimes a_1\otimes...\otimes a_{i-1}\otimes m_j^N(a_i\otimes...\otimes a_n\otimes y).$$

We call this complex the {\it derived tensor product} of $M$ and $N$.

For any  $\A$-algebras $A$ and $B$ we define an $A -
B$-bimodule as a $\Z$-graded vector space $M$ together with linear maps
$$ c_{n_1,n_2}^M:A [ 1 ]^{\otimes n_1} \otimes M \otimes B [ 1 ]^{\otimes n_2} \to M [ 1 ] $$
satisfying the natural compatibility conditions (see e.g. [Ke1]).
If $X$ and $Y$ are formal pointed dg-manifolds corresponding to
$A$ and $B$ respectively then an $A-B$-bimodule is the same
as a dg-bundle ${\cal E}$ over $X\otimes Y$ equipped with a homological vector field $d_{\cal E}$
which is a lift of the vector field $d_X\otimes 1+1\otimes d_Y$.

\begin{exa} \label{diagonal bimodule}Let $A=B=M$. We define a structure of diagonal bimodule
on $A$ by setting $c_{n_1,n_2}^A=m_{n_1+n_2+1}^A$.

\end{exa}

\begin{prp}\label{modules and bimodules algebraically}
1) To have a structure of an $\A$-module on the complex $M$ is the same
as to have a
homomorphism of $\A$-algebras $\phi:A \to \underline{End}_{\bf K} ( M )$,
where ${\bf K}$ is a category of complexes of $\K$-vector spaces.

2) To have a structure of an $A-B$-bimodule on a graded
vector space $M$ is the same as to have a structure
of left $A$-module on $M$ and to have a morphism of
$\A$-algebras $\varphi_{A,B}:B^{op}\to Hom_{A-mod}(M,M)$.

\end{prp} Let $A$ be an $\A$-algebra, $M$ be an $A$-module and
$\varphi:=\varphi_{A,A}: A^{op}\to Hom_{A-mod}(M,M)$ be the corresponding
morphism of $\A$-algebras. Then the dg-algebra $Centr(\varphi)$
is isomorphic to the dg-algebra $Hom_{A-mod}(M,M)$.

If $M=_AM_B$ is an $A-B$-bimodule and $N=_BN_C$ is a $B-C$-bimodule
then the complex $_AM_B\otimes_B{_B}N_C$ carries an $A-C$-bimodule
structure. It is called the {\it tensor product} of $M$ and $N$.

Let $f : X \to Y$ be a homomorphism of formal pointed dg-manifolds
corresponding to a homomorphism of $\A$-algebras $A \to B$. Recall that in
Section \ref{centralizer} we constructed the formal neighborhood $U_f$ of $f$ in $Maps ( X, Y
)$ and the
$\A$-algebra $Centr(f)$. On the other hand, we have an $A - mod - B$ bimodule
structure on $B$ induced by $f$. Let us denote this bimodule by $M(f)$. We leave
the proof of the following
result as an exercise to the reader. It will not be used in the paper.

\begin{prp}\label{endomorphisms of bimodule} If $B$ is weakly unital then the dg-algebra $End_{A - mod - B} ( M(f) )$ is quasi-isomorphic to
$Centr(f)$.
\end{prp}

$\A$-bimodules will be used in Part \ref{smoothness and compactness} for the study of homologically smooth
$\A$-algebras. In the case of  $\A$-categories  bimodules give rise to $\A$-functors between the corresponding
categories of modules. Tensor product of bimodules  corresponds
to the composition of $\A$-functors.

\subsection{On the tensor product of $\A$-algebras}\label{tensor product of A-infty algebras}

The tensor product of two dg-algebras $A_1$ and $A_2$ is a dg-algebra.
For $\A$-algebras there is no canonical simple formula
for the $\A$-structure on  $A_1\otimes_\K A_2$ which
generalizes the one in the dg-algebras case.
Some complicated formulas were proposed in [SU]. They are not
symmetric with respect to the permutation $(A_1,A_2)\mapsto (A_2,A_1)$.
One can give a definition of the dg-algebra which is
quasi-isomorphic to the one from [SU] in the case when both $A_1$ and
$A_2$ are weakly unital. Namely, we define the $\A$-tensor product
$$A_1``\otimes" A_2=End_{A_1-mod-A_2}(A_1\otimes A_2).$$
Notice that it is a unital dg-algebra.
One can show that the dg-category $A-mod-B$ is equivalent
(as a dg-category) to $A_1``\otimes" A_2^{op}-mod$.

\section{Yoneda lemma}\label{Yoneda lemma}

\subsection{Explicit formulas for the product and differential on $Centr(f)$}\label{explicit formulas for Centr}

Let $A$ be an $\A$-algebra, and $B = End_{{\bf K}} ( A )$ be
the dg-algebra of endomorphisms of $A$ in the category ${\bf K}$ of
complexes of $\K$-vector spaces. Let $f = f_A : A \to B$ be the natural 
$\A$-morphism
coming from the left action of $A$ on itself.
Notice that $B$ is always a unital
dg-algebra, while $A$ can be non-unital. The aim of this Section
is to discuss the relationship between $A$ and $Centr(f_A)$.
This is the simplest case of the $\A$-version of Yoneda lemma
(the general case easily follows from this one. See also
[Lyu],[LyuOv]).

As a graded vector space $Centr(f_A)$ is isomorphic
to $\prod_{n\ge 0}\ihom(A^{\otimes (n+1)},A)[-n]$.

Let us describe the product on $Centr ( f )$ for $f=f_A$. Let $\phi, \psi$
be two homogeneous elements of $Centr ( f )$. Then
$$ ( \phi \cdot \psi ) ( a_1, a_2, \ldots, a_N ) = \pm \phi ( a_1, \ldots,
   a_{p - 1}, \psi ( a_p, \ldots, a_N ) ) . $$
Here $\psi$ acts on the last group of variables $a_p, \ldots, a_N$, and we use
the Koszul sign convention for $\A$-algebras in order to determine the sign.

Similarly one has the following formula for the differential (see Section \ref{centralizer}):

$$ ( d \phi ) ( a_1, \ldots, a_N ) =
\sum \pm \phi ( a_1, \ldots, a_s,m_i (a_{s+1}, \ldots,
   a_{s+i}),a_{s+i+1} \ldots, a_N ) + $$
 $$\sum \pm m_i ( a_1, \ldots, a_{s - 1}, \phi ( a_s,
   \ldots, a_j ,..., a_N )) . $$

\subsection{Yoneda homomorphism}\label{Yoneda homomorphism}

If $M$ is  an $A-B$-bimodule then one has a homomorphism
of $\A$-algebras $B^{op}\to Centr(\phi_{A,M})$ (see Propositions \ref{modules and bimodules algebraically}, \ref{endomorphisms of bimodule}). We would like
to apply this general observation in the case of the diagonal bimodule
structure on $A$.
Explicitly, we have the $\A$-morphism
$A^{op} \to End_{mod - A} ( A )$ or, equivalently, the
collection of maps $A^{\otimes m} \to Hom ( A^{\otimes n}, A )$. By
conjugation it gives us a collection of maps
$$ A^{\otimes m} \otimes Hom ( A^{\otimes n}, A ) \to Hom ( A^{\otimes ( m + n
   )}, A ) . $$
In this way we get a natural $\A$-morphism $Yo:A^{op} \to Centr ( f_A )$ called the {\it Yoneda homomorphism}.

\begin{prp}\label{weakly unital algebras and Yoneda}
The $\A$-algebra $A$ is weakly unital if and only if the Yoneda homomorphism
is a quasi-isomorphism.

\end{prp}

{\it Proof.} Since $Centr(f_A)$ is weakly unital, then $A$ must be weakly
unital as long as Yoneda morphism is a quasi-isomorphism.

Let us prove the opposite statement. We assume that
$A$ is weakly unital. It suffices to prove that the
cone $Cone(Yo)$ of the Yoneda homomorphism has trivial cohomology.
Thus we need to prove that the cone of the morphism of complexes
$$ ( A^{op}, m_1 ) \to ( \oplus_{n \ge 1} Hom ( A^{\otimes n}, A ), m_1^{Centr( f_A )} ) . $$
is contractible. In order to see this, one considers the extended complex $A \oplus
Centr ( f_A )$. It has natural filtration arising from  tensor powers of
$A$. The corresponding spectral sequence collapses, which gives an explicit
homotopy of the extended complex to the trivial one. This implies the desired
quasi-isomorphism of $H^0 ( A^{op} )$ and $H^0 ( Centr ( f_A ) )$.
$\blacksquare$

\begin{rmk}\label{Centr and unital algebras} It looks like the construction of $Centr(f_A)$ is the first
known canonical construction of a unital dg-algebra quasi-isomorphic to a
given $\A$-algebra (canonical but not functorial).
This is true even in the case of strictly unital $\A$-algebras.
Standard construction via bar and cobar resolutions gives
a non-unital dg-algebra.

\end{rmk}

\part{Smoothness and compactness}\label{smoothness and compactness}

\section{Hochschild cochain and chain complexes of an $\A$-algebra}\label{Hochschild cochain and chain complexes}

\subsection{Hochschild cochain complex}\label{Hochschild cochain complex}

We change the notation for the homological vector field to $Q$, since the letter $d$
will be used for the differential.\footnote{We recall that the super version of the notion of
formal dg-manifold was introduced by A. Schwarz under the name ``$Q$-manifold". Here letter $Q$ refers
to the supercharge notation from Quantum Field Theory.}
Let $( ( X, pt ), Q )$ be a
non-commutative formal pointed dg-manifold corresponding to a non-unital
$\A$-algebra $A$, and $Vect(X)$ the graded Lie algebra of vector fields on $X$ (i.e. continuous derivations of ${\cal O}(X)$).

We denote by $C^{\bullet} ( A, A ) : = C^{\bullet} ( X, X ):=Vect(X)[-1]$ the
\textit{Hochschild cochain complex} of $A$.
As a $\Z$-graded vector space
$$C^{\bullet} ( A, A )=\prod_{n\ge 0}\ihom_{\cal C}(A[1]^{\otimes n},A).$$
The differential on $C^{\bullet} ( A, A )$ is given by $[ Q,
\bullet ]$. Algebraically, $C^{\bullet} ( A, A ) [ 1 ]$ is a DGLA of
derivations of the coalgebra $T ( A [ 1 ] )$ (see Section \ref{A-infty algebras}).

\begin{thm}\label{tangent space}
Let $X$ be a non-commutative formal pointed dg-manifold and $C^{\bullet} ( X,
X )$ be the Hochschild cochain complex. Then one has the following
quasi-isomorphism of complexes
$$ C^{\bullet} ( X, X ) [ 1 ] \simeq T_{id_X} ( Maps ( X, X ) ), $$
where $T_{id_X}$ denotes the tangent complex at the identity map.

\end{thm}
{\it Proof.} Notice that $Maps ( Spec ( k [ \varepsilon ] / (
\varepsilon^2 ) ) \otimes X, X )$ is the non-commutative dg ind-manifold of
vector fields on $X$. The tangent space $T_{id_X}$ from the theorem is the DGLA of vector fields on $X$.
It can be
identified with the set of such $f \in Maps ( Spec ( k [ \varepsilon ] / (
\varepsilon^2 ) ) \otimes X, X )$ that $f |_{ \{ pt \} \otimes X } = id_X$. On the other hand the DGLA $C^{\bullet} ( X, X ) [ 1 ]$ is  the
DGLA of vector fields on $X$. The theorem follows. $\blacksquare$

The Hochschild cochain complex admits a couple of other interpretations.
We leave to the reader to check the equivalence of all of them.
Namely:

a) $C^{\bullet} ( A, A ) \simeq Centr(id_A)$. 

b) For a weakly unital $\A$-algebra
$A$ one has
$C^{\bullet} ( A, A ) \simeq Hom_{A-mod-A}(A,A)$. 

Both a) and b) are quasi-isomorphisms of complexes.

\begin{rmk} \label{relation to Deligne\rq{}s conjecture}Interpretation of $C^{\bullet} ( A, A )[1]=C^{\bullet} ( X, X )[1]$ as vector fields gives a DGLA
structure on this space. It is a Lie algebra of the ``commutative" formal
group in $Vect_\K^{\Z}$, which is an abelianization of the non-commutative
formal group of inner (in the sense of tensor categories) automorphisms $\underline{Aut}(X)\subset Maps(X,X)$.
Because of this non-commutative structure underlying the Hochschild cochain
complex, it is natural to expect that $C^{\bullet} ( A, A )[1]$ carries
more structures than just DGLA. Indeed, Deligne's conjecture (see e.g. [KoSo1]
and the last section of this paper) claims that the DGLA algebra structure
on $C^{\bullet} ( A, A )[1]$ can be extended to a structure of an algebra
over the operad of singular chains of the topological operad of little
discs. Graded Lie algebra structure can be recovered from cells of highest
dimension in the cell decomposition of the topological operad.

\end{rmk}

\subsection{Hochschild chain complex}\label{Hochschild chain complex}

In this subsection we are going to construct a complex of $\K$-vector spaces
which is dual to the Hochschild chain complex of a non-unital $\A$-algebra.

\subsubsection{Cyclic differential forms of order zero}\label{cyclic differential forms of order zero}

Let $( X, pt )$ be a non-commutative formal pointed
manifold over $k$, and ${\cal O} ( X )$ the algebra of functions on
$X$. For simplicity we will assume that $X$ is finite-dimensional,
i.e. $dim_k\,T_{pt}X<\infty$. If $B = B_X$ is a counital
coalgebra corresponding to $X$ (i.e. the coalgebra of distributions on $X$) then
${\cal O} ( X ) \simeq B^{\ast}$. Let us choose affine coordinates $x_1,
x_2, ..., x_n$ at the marked point $pt$. Then we have an isomorphism of
${\cal O} ( X )$ with the topological algebra $\K \langle \langle x_1, ...,
x_n \rangle \rangle$ of formal series in free graded variables $x_1,
..., x_n$.

We define the space of {\it cyclic differential degree zero forms on $X$}
as
$$ \Omega_{cycl}^0 ( X ) = {\cal O} ( X ) / [ {\cal O} ( X ),
   {\cal O} ( X ) ]_{top}, $$
where $[ {\cal O} ( X ), {\cal O} ( X ) ]_{top}$ denotes the topological
commutator (the closure of the algebraic commutator in the adic topology of
the space of non-commutative formal power series).

Equivalently, we can start
with the graded $\K$-vector space $\Omega_{cycl, dual}^0 ( X )$ defined as the
kernel of the composition $B \to B \otimes B \to \bigwedge^2 B$ (first map is
the coproduct $\Delta : B \to B \otimes B$, while the second one is the
natural projection to the skew-symmetric tensors). Then $\Omega_{cycl}^0 ( X )
\simeq ( \Omega_{cycl, dual}^0 ( X ) )^{\ast}$ (dual vector space).

\subsubsection{Higher order cyclic differential forms}\label{higher order cyclic forms}

We start with the definition of the \textit{odd tangent bundle} $T [ 1 ] X$. This is the dg-analog of the total space of the tangent
supervector bundle to a supermanifold with the
changed parity of fibers. It is more convenient to describe this formal
manifold in terms of algebras rather than coalgebras. Namely, the algebra of
functions ${\cal O} ( T [ 1 ] X )$ is a unital topological algebra
isomorphic to the algebra of formal power series $\K \langle \langle x_i, dx_i \rangle \rangle, 1 \le i \le n$,
where $deg \hspace{0.25em} dx_i = deg \hspace{0.25em} x_i + 1$
(we do not impose any commutativity relations between generators).
More invariant
description involves the odd line. Namely, let $t_1 : = Spc ( B_1 )$, where $(
B_1 )^{\ast} = \K \langle \langle\xi\rangle \rangle / ( \xi^2 ), deg \hspace{0.25em} \xi = +
1$. Then we define $T [ 1 ] X$ as the formal neighborhood in $Maps ( t_1, X )$
of the point $p$ which is the composition of $pt$ with the trivial map of
$t_1$ into the point $Spc ( \K )$.

\begin{defn} \label{de Rham and cyclic forms} a) The graded vector space
$$ {\cal O} ( T [ 1 ] X ) = \Omega^{\bullet} ( X ) = \prod_{m \ge 0}
   \Omega^m ( X ) $$
is called the space of de Rham differential forms on $X$.

b) The graded space
$$ \Omega^0_{cycl} ( T [ 1 ] X ) = \prod_{m \ge 0} \Omega^m_{cycl} ( X ) $$
is called the space of cyclic differential forms on $X$.

\end{defn}

In coordinate description the grading is given by the total number of factors $dx_i$.
Clearly each space $\Omega^n_{cycl} ( X ), n \ge 0$ is dual to some vector
space $\Omega^n_{cycl, dual} ( X )$ equipped with the discrete topology (since this is true for $\Omega^0 ( T [ 1
] X )$).

The {\it de Rham differential on $\Omega^{\bullet} ( X )$} corresponds to
the vector field $\partial / \partial \xi$. This can be seen from the  description which uses the odd
line, where the variable $\xi$ has the same meaning. Since $\Omega^0_{cycl}$ is given by the
natural (i.e. functorial) construction, the de Rham differential descends to the
subspace of cyclic differential forms. We will denote the former by $d_{DR}$
and the latter by $d_{cycl}$.

The space of \textit{cyclic 1-forms} $\Omega_{cycl}^1 ( X )$ is a
(topological) span  of expressions $x_1 x_2 ...x_l
\hspace{0.25em} dx_j, x_i \in {\cal O} ( X )$. Equivalently, the space of cyclic 1-forms
consists of expressions $\sum_{1 \le i \le n} f_i ( x_1, ..., x_n )
\hspace{0.25em} dx_i$ where $f_i \in \K \langle \langle x_1, ..., x_n \rangle
\rangle$.

There is a map $\varphi : \Omega_{cycl}^1 ( X ) \to {\cal O} ( X
)_{red} : = {\cal O} ( X ) / \K$, which is defined on $\Omega^1 ( X )$ by
the formula $adb \mapsto [ a, b ]$, and then one checks that the induced map on the cyclic
$1$-forms is well-defined. This map does not have an analog in the
commutative case.\footnote{ V.Ginzburg pointed out that the geometric meaning of the map $\varphi$
as a ``contraction with double derivation" was suggested in Section 5.4 of [Gi3].}

\subsubsection{Non-commutative Cartan calculus}\label{nc calculus}

Let $X$ be a formal graded manifold over a field $\K$.
We denote by $\g:=\g_X$ the graded Lie algebra of continuous linear maps ${\cal O}(T[1]X)\to {\cal O}(T[1]X)$ generated by
de Rham differential $d=d_{dR}$ and contraction maps $i_{\xi}, \xi\in Vect(X)$ which are defined by the formulas $i_{\xi}(f)=0, i_{\xi}(df)=\xi(f)$ for all $f\in {\cal O}(T[1]X)$. Let us define
the Lie derivative $Lie_{\xi}=[d,i_{\xi}]$, where the RHS is the graded commutator.
Then one can easily checks the usual formulas of the Cartan calculus

$$[d,d]=0, Lie_{\xi}=[d,i_{\xi}], [d, Lie_{\xi}]=0,$$
$$[Lie_{\xi},i_{\eta}]=i_{[\xi,\eta]}, [Lie_{\xi},Lie_{\eta}]=Lie_{[\xi,\eta]}, [i_{\xi},i_{\eta}]=0,$$
for any $\xi,\eta\in Vect(X)$.

By naturality, the graded Lie algebra $\g_X$ acts on the space
$\Omega_{cycl}^{\bullet}(X)$ as well as one the dual  space
$(\Omega_{cycl}^{\bullet}(X))^{\ast}$.

\subsubsection{Differential on the Hochschild chain complex}\label{differential on Hochschild chain complex}
Let $Q$ be a homological vector field on $(X,pt)$. Then
$A=T_{pt}X[-1]$ is a non-unital $\A$-algebra.

We define the {\it dual Hochschild chain complex}
$(C_{\bullet} ( A, A ))^{\ast}$ as the graded vector space 
$ \Omega^1_{cycl} ( X )[2] $ endowed with the differential $Lie_Q$.
Our terminology is explained by the observation that $\Omega^1_{cycl} ( X )[2]$
is dual to the conventional Hochschild chain complex
$$ C_{\bullet} ( A, A )= \oplus_{n \ge 0}
( A [ 1 ] )^{\otimes n}\otimes A. $$

Notice that we use the cohomological grading on $C_{\bullet} ( A, A )$,
i.e. chains of
degree $n$ in conventional (homological) grading have degree $- n$ in our grading.
In our grading the differential has degree $+1$.

In coordinates the isomorphism  identifies an
element $f_i ( x_1, ..., x_n ) \otimes x_i \in ( A [ 1 ]^{\otimes n} \otimes A
)^{\ast}$ with the homogeneous element $f_i ( x_1, ..., x_n ) \hspace{0.25em}
dx_i \in \Omega^1_{cycl} ( X )$. Here $x_i \in ( A [ 1 ] )^{\ast}, 1 \le
i \le n$ are ``affine coordinates\rq{}\rq{} on $X$.

The graded Lie algebra $Vect ( X )$ of vector fields of all degrees acts on any
functorially defined space, in particular, on all spaces $\Omega^j ( X ),
\Omega_{cycl}^j ( X )$, etc. Then we have a differential of degree $+ 1$ on $\Omega_{cycl}^j ( X )$ given by $b = Lie_Q$. There is an explicit
formula for the differential $b$ on $C_{\bullet} ( A, A )$ (cf. [T]):
$$ b ( a_0 \otimes ... \otimes a_n ) = \sum \pm a_0 \otimes ... \otimes m_l (
   a_i \otimes ... \otimes a_j ) \otimes ... \otimes a_n $$
$$ + \sum \pm m_l ( a_j \otimes ... \otimes a_n \otimes a_0 \otimes ...
   \otimes a_i ) \otimes a_{i + 1} \otimes ... \otimes a_{j - 1} . $$
It is convenient to depict a cyclic monomial $a_0 \otimes ... \otimes a_n$ in
the following way. We draw a clockwise oriented circle with $n + 1$ points
labeled from $0$ to $n$ such that one point is marked. We assign the elements
$a_0, a_1, ..., a_n$ to the points with the corresponding labels, putting
$a_0$ at the marked point.



\includegraphics{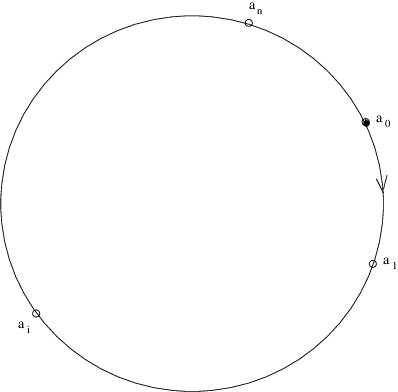}

{\vspace{3mm}}

Then we can write $b = b_1 + b_2$ where $b_1$ is the sum (with appropriate
signs) of the expressions depicted below:

{\vspace{3mm}}

\includegraphics{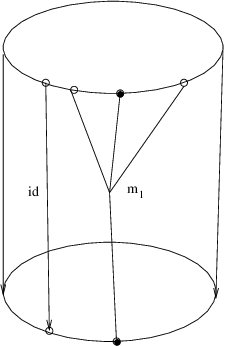}
{\vspace{5mm}}

Similarly, $b_2$ is the sum (with appropriate signs) of the expressions
depicted below:


\includegraphics{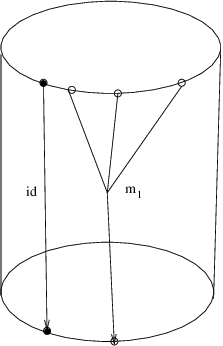}
{\vspace{3mm}}

In both cases maps $m_l$ are applied to a consequitive cyclically ordered
sequence of elements of $A$ assigned to the points on the top circle. The
identity map is applied to the remaining elements. Marked point on the top
circle is the position of the element of $a_0$. Marked point on the bottom
circle depicts the first tensor factor of the corresponding summand of $b$. In
both cases we start cyclic count of tensor factors clockwise from the marked point.

\subsection{The case of strictly unital $\A$-algebras}\label{case of strictly unital algebras}

Let $A$ be a strictly unital $\A$-algebra.
There is a {\it reduced} Hochschild chain complex
$$C_{\bullet}^{red}(A,A)=\oplus_{n\ge 0}A\otimes ((A/\K\cdot 1)[1])^{\otimes n},$$
which is the quotient of $C_{\bullet}(A,A)$. Similarly there is
a reduced Hochschild cochain complex
$$ C^{\bullet}_{red}(A,A)=\prod_{n\ge 0}\ihom_{\cal C}((A/\K\cdot 1)[1]^{\otimes n},A),$$
which is a subcomplex of the Hochschild cochain complex $C^{\bullet}(A,A)$.

Also, $C_{\bullet} ( A, A )$
carries also
the ``Connes's differential " $B$ of degree $- 1$ which is called sometimes ``de Rham differential" and which is
 given by the formula (see [Co],
[T])
$$ B ( a_0 \otimes ... \otimes a_n ) = \sum_i \pm 1 \otimes a_i \otimes ...
   \otimes a_n \otimes a_0 \otimes ... \otimes a_{i - 1}, B^2=0, Bb+bB=0 . $$
Here is the graphical description of $B$, which  will receive an explanation in the
section devoted to generalized Deligne's conjecture:

{\vspace{3mm}}


\includegraphics{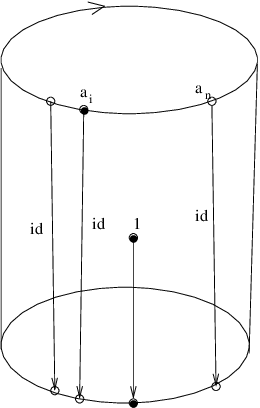}
{\vspace{3mm}}

Let $u$ be an independent variable of degree $+ 2$. It
follows that for a strictly unital $\A$-algebra $A$ one has a
differential $b + uB$ of degree $+ 1$ on the graded vector space $C_{\bullet}
( A, A ) [ [ u ] ]$ which makes the latter into a complex called
\textit{negative cyclic complex} (see [Co], [T]). In fact
$b+uB$ is a
differential on a smaller complex $C_{\bullet} ( A, A ) [ u ]$.
In the non-unital case one can use Cuntz-Quillen complex
instead of the negative cyclic complex (see next subsection).

\subsection{Non-unital case: Cuntz-Quillen complex}\label{non-unital case}

In this subsection we are going to present a
formal dg-version of the mixed
complex introduced by Cuntz and Quillen (see [CQ1]). In the previous
subsection we introduced the Connes differential $B$ in the case of strictly unital
$\A$-algebras. In the non-unital case the construction has to be modified. Let $X=A[1]_{form}$ be the corresponding
non-commutative formal pointed dg-manifold. The algebra of functions ${\cal O} ( X )\simeq
\prod_{n \ge 0} ( A [ 1 ]^{\otimes n} )^{\ast}$ is a complex with the
differential $Lie_Q$.

\begin{prp}\label{cohomology in weakly unital case}
If $A$ is weakly unital then all non-zero cohomology of the complex
${\cal O} ( X )$ are trivial, and $H^0 ( {\cal O} ( X ) )
\simeq \K$.

\end{prp}
{\it Proof.} Let us calculate the cohomology using the spectral sequence
associated with the filtration $\prod_{n\ge n_0}(A[1]^{\otimes n})^{\ast}$.
The term $E_1$ of the spectral sequence is isomorphic to the
complex $\prod_{n\ge 0}((H^{\bullet}(A[1],m_1))^{\otimes n})^{\ast}$
with the differential induced by the multiplication $m_2^A$ on
$H^{\bullet}(A,m_1^A)$. By assumption $H^{\bullet}(A,m_1^A)$ is a unital
algebra, hence all  cohomology groups of the complex vanish except of the zeroth one,
which is isomorphic to $\K$.
This concludes the proof.
$\blacksquare$.

It follows from the above Proposition \ref{cohomology in weakly unital case} that the complex
${\cal O} ( X ) / \K $ is acyclic.
We have the following two morphisms of complexes

$$d_{cycl}:({\cal O}(X)/\K\cdot 1, Lie_Q)\to (\Omega_{cycl}^1(X), Lie_Q)$$
and
$$\varphi: (\Omega_{cycl}^1(X), Lie_Q)\to({\cal O}(X)/k\cdot 1, Lie_Q). $$
Here $d_{cycl}$ and $\varphi$ were introduced in the Section \ref{higher order cyclic forms}.
We have: $deg(d_{cycl})=+1, deg(\varphi)=-1, d_{cycl}\circ \varphi=0,
\varphi\circ d_{cycl}=0.$.

Let us consider a {\it modified} Hochschild chain complex
$$C_{\bullet}^{mod}(A,A):=(\Omega_{cycl}^1(X)[2])^{\ast}\oplus ({\cal O}(X)/\K\cdot 1)^{\ast}$$ with the differential

$b=\left( \begin{array}{cc}
  (Lie_Q)^{\ast} & \varphi^{\ast}\\
  0 & (Lie_Q)^{\ast}
\end{array} \right)$

Let

$B=\left( \begin{array}{cc}
  0 & 0\\
  d_{cycl}^{\ast} & 0
\end{array} \right)$
be an endomorphism of $C_{\bullet}^{mod}(A,A)$ of degree $-1$.
Then $B^2=0$.
Let $u$ be a formal variable of degree $+2$.
We define  {\it modified negative cyclic, periodic cyclic and cyclic chain} complexes
such as follows

$$CC_{\bullet}^{-,mod}(A)=(C_{\bullet}^{mod}(A,A)[[u]],b+uB),$$

$$CP_{\bullet}^{mod}(A)=(C_{\bullet}^{mod}(A,A)((u)),b+uB),$$

$$CC_{\bullet}^{mod}(A)=(CP_{\bullet}^{mod}(A)/CC_{\bullet}^{-,mod}(A))[-2].$$

The corresponding cohomology groups are called (modified) {\it negative cyclic, periodic cyclic and cyclic} homology respectively.

For unital dg-algebras these complexes are quasi-isomorphic to the
standard ones. If $char\,\K=0$ and $A$ is weakly unital then
$CC_{\bullet}^{-,mod}(A)$ is quasi-isomorphic to the complex
$(\Omega_{cycl}^0(X), Lie_Q)^{\ast}$. Notice that the $\K[[u]]$-module 
structure on the cohomology $H^{\bullet}((\Omega_{cycl}^0(X), Lie_Q)^{\ast})$ is not transparent  from its definition.

\section{Homologically smooth and compact $\A$-algebras}\label{homologically smooth and compact algebras}
From now on we will assume that all $\A$-algebras are weakly unital unless we say otherwise.
\subsection{Homological smoothness}\label{homological smoothness}
Let $A$ be an $\A$-algebra over $\K$ and $E_1, E_2, ..., E_n$ be a sequence of
$A$-modules. Let us consider a sequence $(E_{\le i})_{1\le i\le n}$ of $A$-modules together with exact triangles
$$ E_i \to E_{\le i} \to E_{i + 1} \to E_i [ 1 ], $$
such that $E_{\le 1} = E_1$.

We will call $E_{\le n}$ an {\it extension} of the sequence
$E_1,...,E_n$.

Notice  that we can speak in the same way about extensions in the category
of $A-A$-bimodules.

\begin{defn} \label{perfect bimodules}1) A perfect $A$-module is the one which is quasi-isomorphic to a direct
summand of an extension of a sequence of modules each of which is quasi-isomorphic
to $A [ n ], n \in \Z$.

2) A perfect $A - A$-bimodule is the one which is quasi-isomorphic to a direct
summand of an extension of a sequence consisting of bimodules each of which is
quasi-isomorphic to $( A \otimes A )[ n ], n \in \Z$.

\end{defn}

Perfect $A$-modules form a full subcategory $Perf_A$ of the dg-category $A -
mod$. Perfect $A - A$-bimodules form a full subcategory $Perf_{A - mod - A}$
of the category of $A - A$-bimodules.\footnote{Sometimes $Perf_A$ is called a thick triangulated
subcategory of $A-mod$ generated by $A$ and denoted by $\langle A \rangle$. In the case
of $A-A$-bimodules we have a thick triangulated subcategory generated by $A\otimes A$, which is denoted by
$\langle A\otimes A \rangle$.}

\begin{defn} \label{homologically smooth algebras}We say that an $\A$-algebra $A$ is homologically smooth if it is a
perfect $A - A$-bimodule (equivalently, $A$ is a perfect module over the
$\A$-algebra $A ``\otimes" A^{op}$).
\end{defn}

\begin{rmk} \label{saturated categories}An $A-B$-bimodule $M$ gives rise to a dg-functor
$B-mod\to A-mod$ such that $V\mapsto M\otimes_BV$. The diagonal
bimodule $A$ corresponds to the identity functor $Id_{A-mod}: A-mod\to A-mod$.
The notion of homological smoothness can be generalized to the framework
of $\A$-categories. The corresponding notion of saturated $\A$-category
can be spelled out entirely in terms of the identity functor.

\end{rmk}

Let us list few examples of homologically smooth $\A$-algebras.

\begin{exa}\label{examples of homologically smooth algebras} a) Algebra of functions on a smooth affine scheme.

b) $A = \K [ x_1, ..., x_n ]_q$, which is the algebra of polynomials in
variables $x_i, 1 \le i \le n$ subject to the relations $x_i x_j = q_{ij}\,x_j x_i$,
where $q_{ij}\in \K^{\ast}$ satisfy the properties $q_{ii}=1, q_{ij}q_{ji}=1$.
More generally, all quadratic Koszul algebras, which are deformations of
polynomial algebras are homologically smooth.

c) Algebras of regular functions on quantum groups (see [KorSo]).

d) Free algebras $\K \langle x_1, ..., x_n \rangle$.

e) Finite-dimensional associative algebras of finite homological dimension.

f) If $X$ is a smooth scheme over $\K$
then the bounded derived category $D^b(Perf(X))$
of the category of perfect complexes (it is equivalent to $D^b(Coh(X))$)
has a generator $P$ (see [BvB]). Then the dg-algebra $A=End(P)$
(here we understand endomorphisms in the ``derived sense", see [Ke2])
is a homologically smooth algebra.

\end{exa}

Let us introduce an $A - A$-bimodule
$A^! = Hom_{A - mod - A} ( A, A \otimes A)$ (cf. [Gi2]).
The structure of an $A-A$-bimodule is defined similarly to the case
of associative algebras.

\begin{prp}\label{homologically smoothness and perfectness}
If $A$ is homologically smooth then $A^!$ is a perfect $A-A$-bimodule.

\end{prp}

{\it Proof.} We observe that $Hom_{C - mod} ( C, C )$ is a
dg-algebra for any $\A$-algebra $C$.
The Yoneda embedding $C \to Hom_{C - mod} ( C, C )$ is a
quasi-isomorphism of $\A$-algebras.
Let us apply this observation to $C=A\otimes A^{op}$.
Then using the $\A$-algebra $A``\otimes"A^{op}$ (see Section \ref{tensor product of A-infty algebras}) we obtain a quasi-isomorphism of $A-A$-bimodules $Hom_{A - mod -
A} ( A \otimes A, A \otimes A ) \simeq A \otimes A$. By assumption $A$ is
quasi-isomorphic (as an $\A$-bimodule) to a direct summand in an extension of a sequence
$(A \otimes A) [ n_i ]$ for  $n_i \in \Z$. Hence $Hom_{A - mod - A} ( A \otimes
A, A \otimes A )$ is quasi-isomorphic to a direct summand in an extension of a sequence
$(A \otimes A) [ m_i ]$ for $m_i \in \Z$. The result follows. $\blacksquare$

\begin{defn} \label{inverse dualizing bimodule}The bimodule $A^!$ is called the inverse  dualizing
bimodule.

\end{defn}

The terminology is explained by an observation that if $A=End(P)$ where $P$
is a generator of of $Perf(X)$ (see Example \ref{examples of homologically smooth algebras} f)) then the bimodule $A^!$
corresponds to the functor $F\mapsto F\otimes K_X^{-1}[-dim\,X]$, where
$K_X$ is the canonical class of $X$.\footnote{We thank to Amnon Yekutieli for pointing out
that the inverse dualizing module was first mentioned in the paper by Michel van den Bergh ``Existence theorems for dualizing complexes over non-commutative graded and filtered rings", J. Algebra, 195:2, 1997, 662-679.}

\begin{rmk}\label{fibrant dg algebras} In [ToVa] the authors introduced a stronger notion of fibrant dg-algebra.
Informally it corresponds to ``non-commutative homologically smooth affine schemes of finite type".
In the compact case (see the next section) this notion coincides with the one of homologically smooth algebra.

\end{rmk}

\subsection{Compact $\A$-algebras}\label{compact algebras}

\begin{defn} \label{definition of compact algebras}We say that an $\A$-algebra $A$ is compact if the cohomology
$H^{\bullet} ( A, m_1 ) $ is finite-dimensional.
\end{defn}

\begin{exa}\label{examples of compact algebras}
a) If $dim_\K A<\infty$ then $A$ is compact.

b) Let $X/\K$ be a proper scheme of finite type.
According to [BvB] there exists a compact
dg-algebra $A$ such that $Perf_A$ is equivalent
to $D^b(Coh(X))$.

c) If $Y\subset X$ is a proper subscheme (possibly singular) of
a smooth scheme $X$ then the bounded derived category $D^b_Y(Perf(X))$
of the category of perfect complexes on $X$,
which are supported on $Y$ has a generator $P$ such that $A=End(P)$
is compact.  In general it is not homologically smooth for $Y\ne X$.
More generally, one can replace $X$ by a formal smooth scheme containing $Y$, e.g. by the
formal neighborhood of $Y$ in the ambient smooth scheme.
In particular, for $Y=\{pt\}\subset X={\bf A}^1$ and the generator ${\cal O}_Y$ of
$D^b(Perf(X))$ the corresponding graded algebra is isomorphic to $\K\langle \xi\rangle/(\xi^2)$,
where $deg\,\xi=1$.

\end{exa}

\begin{prp}\label{compactness and homological smoothnes imply finiteness of cohomology}

If $A$ is compact and homologically smooth then the Hochschild homology and
cohomology of $A$ are finite-dimensional.

\end{prp}
{\it Proof.} a) Let us start with Hochschild cohomology. We have
an isomorphism of complexes
$C^{\bullet} ( A, A ) \simeq Hom_{A - mod - A} ( A, A )$. Since $A$ is
homologically smooth the latter complex
is quasi-isomorphic to a direct summand of an extension
of the bimodule $Hom_{A - mod - A} ( A
\otimes A, A \otimes A )$. The latter complex in turn is
quasi-isomorphic  to  $A \otimes A$ (see the proof of the
Proposition ref{homologically smoothness and perfectness}). Since $A$ is compact,
the complex $A \otimes A$ has finite-dimensional cohomology.
Therefore any perfect $A-A$-bimodule enjoys the same property.
We conclude that the Hochschild cohomology groups are finite-dimensional
vector spaces.

b) Let us consider the case of Hochschild homology.
With any $A - A$-bimodule $E$ we
associate a complex of vector spaces
$E^{\sharp} = \oplus_{n \ge 0} A [ 1 ]^{\otimes n}\otimes E$ (cf. [Gi2]). The differential on
$E^{\sharp}$ is given by the same formulas as the Hochschild
differential for $C_{\bullet} ( A, A )$ with the
only change: we place an element $e \in E$ instead of an element of $A$
at the marked vertex (see Section \ref{Hochschild cochain and chain complexes}).
Taking $E=A$ with the  structure
of the diagonal $A-A$-bimodule  we obtain $A^{\sharp} =
C_{\bullet} ( A, A )$. On the other hand,
it is easy to see that the complex $(A\otimes A)^{\sharp}$ is quasi-isomorphic
to $(A,m_1)$, since $(A\otimes A)^{\sharp}$ is the quotient of the canonical
free resolution (bar resolution) for $A$ by a subcomplex $A$.
The construction of $E^{\sharp}$ is functorial, hence $A^{\sharp}$ is quasi-isomorphic
to a direct summand of an extension (in the category of complexes)
of a shift of $(A\otimes A)^{\sharp}$, because $A$ is smooth.
Since $A^{\sharp}=C_{\bullet} ( A, A )$ we see that
the  Hochschild homology $H_{\bullet}(A,A)$ is isomorphic to a direct
summand of the cohomology of an extension of a sequence of $k$-modules
$(A[n_i],m_1)$. Since the vector space $H^{\bullet}(A,m_1)$ is finite-dimensional the result follows. $\blacksquare$

\begin{rmk}\label{different qis of Hochschild chain complexes}
For a homologically smooth compact $\A$-algebra $A$ one has
a quasi-isomorphism of complexes
$C_{\bullet} ( A, A ) \simeq Hom_{A - mod - A} ( A^!, A )$
Also, the complex  $Hom_{A - mod - A}( M^!, N )$
is quasi-isomorphic to $( M \otimes_A N )^{\sharp}$ for  two $A -
A$-bimodules $M, N$, such that $M$ is perfect.
Here $M^! : = Hom_{A - mod - A} ( M, A \otimes A )$
Having this in mind one can offer a version of the above proof which
uses the isomorphism
$$Hom_{A - mod - A} ( A^!, A)\simeq
Hom_{A - mod - A} (Hom_{A-mod-A}(A,A\otimes A),A).$$
Indeed, since
$A$ is homologically smooth the bimodule
$Hom_{A - mod - A}(A,A\otimes A)$ is
quasi-isomorphic to a direct summand $P$ of an extension of a shift of
$Hom_{A - mod - A}(A\otimes A,A\otimes A)\simeq A\otimes A$. Similarly,
$Hom_{A - mod - A}(P, A)$ is quasi-isomorphic to a  direct summand of an extension of a shift
of $Hom_{A - mod - A}(A\otimes A,A\otimes A)\simeq A\otimes A$.
Combining the above computations we see that the complex
$C_{\bullet} ( A, A)$ is quasi-isomorphic to a direct
summand of an extension of a shift of the complex  $A\otimes A$. The latter
has finite-dimensional cohomology, since
$A$ enjoys this property.
\end{rmk}

Algebras of finite quivers are homologically smooth and compact. Besides of that there  are two main sources
of homologically smooth compact $\Z$-graded $\A$-algebras.

\begin{exa}\label{examples of smooth proper categories} a) Combining Examples \ref{examples of homologically smooth algebras} f) and \ref{examples of compact algebras} b) we see that if $X/\K$ 
is a smooth proper variety then the
bounded derived category $D^b ( Coh ( X ) )$ is equivalent to the category $Perf_A$
for a homologically smooth compact $\A$-algebra $A/\K$.

b) According to [Se] the derived category $D^b ( F ( X ) )$ of the Fukaya
category $F(X)$ of a K3 surface $X$ is equivalent to $Perf_A$ for a homologically
smooth compact $\A$-algebra $A$. The latter is generated by Lagrangian
spheres, which are vanishing cycles at the critical points for a fibration of
$X$ over ${\bf CP}^1$. This result can be generalized to
other Calabi-Yau manifolds.

\end{exa}
In $\Z/2$-graded case examples of homologically smooth compact
$\A$-algebras come from Landau-Ginzburg models, i.e. from Fukaya-Seidel categories (see [Or], [R])
as well as from Fukaya categories of Fano varieties.

\begin{rmk} \label{moduli stack of nc smooth proper varieties}Formal deformation theory of smooth compact $\A$-algebras
gives a finite-dimensional formal pointed (commutative) dg-manifold.
The global moduli stack can be constructed using methods of [ToVa]).
It can be thought of as a moduli stack of non-commutative
smooth proper varieties.

\end{rmk}

\section{Degeneration Hodge-to-de Rham }\label{Hodge-to-de Rham}

\subsection{Main conjecture}\label{main conjecture}

Let us assume that $char \hspace{0.25em} \K = 0$ and $A$ is a weakly unital
$\A$-algebra, which can be $\Z$-graded or $\Z/2$-graded.

For any $n \ge 0$ we define
the {\it truncated modified  negative cyclic complex}
$C_{\bullet}^{mod,( n )} ( A, A ) = C_{\bullet}^{mod} ( A, A ) \otimes \K [ u ] / ( u^n
)$, where $deg \hspace{0.25em} u = +2$. It is a complex with the differential
$b + uB$. Its cohomology will be denoted by $H^{\bullet} ( C_{\bullet}^{mod,( n )}
( A, A ) )$ and will be called {\it truncated modified  negative cyclic homology}.

\begin{defn} \label{degeneration property}We say that an $\A$-algebra $A$ satisfies the degeneration property
if for any $n \ge 1$ one
has the following: the cohomology $H^{\bullet} ( C_{\bullet}^{mod,( n )} ( A, A ) )$ is a flat $\K [ u ] / ( u^n
)$-module.

\end{defn}

\begin{conj}\label{degeneration conjecture}(Degeneration Hodge-to-de Rham). Let $A$ be a weakly
unital compact homologically smooth $\A$-algebra. Then it satisfies
the degeneration property.

\end{conj}

We call the Conjecture \ref{degeneration conjecture}  the {\it degeneration conjecture}.

\begin{cor} \label{degeneration property and cyclic homology} If  $A$ satisfies the degeneration property then
the modified  negative cyclic homology coincides with $\varprojlim_n H^{\bullet} (
C_{\bullet}^{mod,( n )} ( A, A ) )$, and furthermore it is a flat $\K [ [ u ] ]$-module.
\end{cor}

Assuming the degeneration property for $A$ we see that there is a $\Z$-graded vector
bundle $\xi_A$ over ${\bf A}^1_{form}[-2] = Spf ( \K [ [ u ] ] )$ with the
space of sections isomorphic to
$$ \varprojlim_n H^{\bullet} ( C_{\bullet}^{mod,( n )} ( A, A ) ) =HC_{\bullet}^{-,mod}( A ), $$
which is the  modified negative cyclic homology of $A$. We will skip the word ``modified\rq\rq{} from the notation. The fiber of
$\xi_A$ at $u = 0$ is isomorphic to the Hochschild homology $H_{\bullet}^{mod} ( A,
A ):=H_{\bullet}(C_{\bullet}(A,A))$.

Notice that $\Z$-graded $\K((u))$-module $HP_{\bullet}^{mod}(A)$
of periodic cyclic homology can be described in terms of just one
$\Z/2$-graded vector space
$HP_{even}^{mod}(A)\oplus \Pi HP_{odd}^{mod}(A)$, where
$HP_{even}^{mod}(A)$ (resp. $HP_{odd}^{mod}(A)$) consists
of elements of degree zero (resp. degree $+1$) of
$HP_{\bullet}^{mod}(A)$ and $\Pi$ is the functor of changing the parity.
We can interpret $\xi_A$ in terms of the $\Z/2$-graded geometry (a.k.a. supergeometry) as a
${\bf G}_m$-equivariant supervector bundle over the even formal line ${\bf A}^1_{form}$.
The structure of a ${\bf G}_m$-equivariant supervector bundle $\xi_A$ is equivalent to a filtration $F$ (called Hodge filtration) by
even numbers on $HP_{even}^{mod}(A)$ and by odd numbers on
$HP_{odd}^{mod}(A)$. The associated $\Z$-graded vector space
coincides with $H_{\bullet}(A,A)$.

\begin{rmk} \label{degeneration property over rings}One can speak about
degeneration property (modulo $u^n$) for $\A$-algebras
which are flat over unital commutative $\K$-algebras.
For example, let $R$ be an artinian local $\K$-algebra with the maximal ideal $m$,
and $A$ be a flat $R$-algebra such that $A / m$ is weakly unital, homologically
smooth and compact. Then, assuming the degeneration property for $A / m$, one can easily
see that it holds for $A$ as well. In particular, the negative cyclic homology of
$A$ gives rise to a vector bundle over $Spec ( R ) \times
{\bf A}^1_{form}[-2]$, where the second factor means the formal graded scheme obtained from the formal completion of the affine line with the grading shifted by $2$.
\end{rmk}

We can say few words in support of the degeneration conjecture. One is,
of course, the classical Hodge-to-de Rham degeneration theorem (see Section \ref{relation to Hodge theory}  below). It is an interesting question to express the classical
Hodge theory  algebraically, in terms
of a generator ${\cal E}$ of the derived category of coherent sheaves and the corresponding
$\A$-algebra $A=RHom({\cal E}, {\cal E})$.
The degeneration conjecture also trivially holds for algebras of finite quivers without relations.

In classical algebraic geometry there are basically two approaches
to the proof of degeneration conjecture. One is analytic and uses K\"ahler
metric, Hodge decomposition, etc. Another one is pure algebraic and
uses the technique of reduction to finite characteristic (see [DI]).
Recently Kaledin (see [Kal]) suggested a proof of a version
of the degeneration conjecture based on the reduction to finite
characterstic.

Below we will formulate a conjecture which could lead
to the definition of crystalline cohomology for $\A$-algebras.
Before doing that we remark  that one can define homologically smooth and compact
$\A$-algebras over any commutative ring, in particular, over
the ring of integers $\Z$. We assume that $A$ is a flat $\Z$-module.

\begin{conj} \label{integer version}Suppose that $A$ is a weakly unital $\A$-algebra over $\Z$,
such that it is homologically smooth (but not necessarily compact).
Truncated negative cyclic complexes
$(C_{\bullet}(A,A)\otimes \Z[[u,p]]/(u^n,p^m),b+uB)$ and
$(C_{\bullet}(A,A)\otimes \Z[[u,p]]/(u^n,p^m),b-puB)$
are quasi-isomorphic for all $n,m\ge 1$ and all
prime numbers $p$.

If, in addition, $A$ is compact then the homology of either
of the above complexes is a flat module over $\Z[[u,p]]/(u^n,p^m)$.

\end{conj}

If the Conjecture  \ref{integer version} holds  then the degeneration
conjecture, probably, can be deduced along the lines  of [DI].

One can also make some conjectures about Hochschild complex of an
arbitrary $\A$-algebra, not assuming that it is compact or
homologically smooth. More precisely, let $A$ be a unital
$\A$-algebra over the ring of $p$-adic numbers $\Z_p$. We assume
that $A$ is  topologically free $\Z_p$-module. Let $A_0=A\otimes_{\Z_p} \Z/p$
be the reduction modulo $p$. Then we have the Hochschild complex
$(C_{\bullet}(A_0,A_0),b)$ and the $\Z/2$-graded complex
$(C_{\bullet}(A_0,A_0),b+B)$.
\begin{conj} \label{p-adic version} There is a natural isomorphism of $\Z/2$-graded vector spaces over the field
$\Z/p$:
$$H^{\bullet}(C_{\bullet}(A_0,A_0),b)\simeq H^{\bullet}(C_{\bullet}(A_0,A_0),b+B).$$

There are similar isomorphisms for weakly unital and non-unital $\A$-algebras,
if one replaces $C_{\bullet}(A_0,A_0)$ by $C_{\bullet}^{mod}(A_0,A_0)$.
Also one has similar isomorphisms for $\Z/2$-graded $\A$-algebras.
\end{conj}

The Conjecture \ref{p-adic version} presumably gives an isomorphism used in [DI], but does
not imply the degeneration conjecture.

\begin{rmk} \label{degeneration for categories} One can formulate similar conjectures
for saturated $\A$-categories, which are categorical generalizations of homologically
smooth compact $\A$-algebras. This observation supports the idea
of introducing the category $NCMot$ of non-commutative pure motives.
Objects of the latter are saturated $\A$-categories over a field, and
$Hom_{NCMot}({\cal C}_1,{\cal C}_2)=
K_0(Funct({\cal C}_1,{\cal C}_2))\otimes {\Q}/equiv$. Here ${K}_0$ means the $K_0$-group of the $\A$-category of functors and $equiv$ means
numerical equivalence, i.e. the equivalence relation
generated by the kernel of the Euler form $\langle E, F \rangle:=\chi(RHom(E,F))$, where $\chi$ is the Euler
characteristic. Some of the properties of the category $NCMot$ are discussed in  [Ko4].
One
can formulate non-commutative analogs of Weil and Beilinson conjectures for the category $NCMot$.

\end{rmk}

\subsection{Relationship with the classical Hodge theory}\label{relation to Hodge theory}

Let $X$ be a quasi-projective scheme of finite type over a field $\K$ of
characteristic zero. Then the category $Perf ( X )$ of perfect sheaves on $X$
is equivalent to $H^0 ( A - mod )$, where $A - mod$ is the category of
$\A$-modules over a dg-algebra $A$. Let us recall a construction of $A$.
Consider a complex $E$ of vector bundles which generates the bounded derived
category $D^b ( Perf ( X ) )$ (see [BvB]). Then $A$ is quasi-isomorphic to $RHom ( E, E
)$. More explicitly, let us fix an affine covering $X = \cup_i U_i$. Then the
complex $A:=\oplus_{i_0, i_1, ..., i_n} \Gamma ( U_{i_0} \cap ... \cap U_{i_n},
E^{\ast} \otimes E ) [ -n]$, $n=dim\, X$ computes $RHom ( E, E )$
and carries a structure of dg-algebra. Different choices of $A$ give rise to
equivalent categories $H^0 ( A - mod )$ (derived Morita equivalence).

Properties of $X$ are encoded in the properties of $A$. In particular:

a) $X$ is smooth iff $A$ is homologically smooth;

b) $X$ is compact iff $A$ is compact.

Moreover, if $X$ is smooth then
$$ H^{\bullet} ( A, A ) \simeq Ext^{\bullet}_{D^b(Coh(X\times X))} ( {\cal O}_{\Delta},
   {\cal O}_{\Delta} ) \simeq $$
$$ \oplus_{i, j \ge 0} H^i ( X, \wedge^j T_X )[-(i+j)]] $$
where ${\cal O}_{\Delta}$ is the structure sheaf of the diagonal $\Delta
\subset X \times X$.

Similarly
$$ H_{\bullet} ( A, A ) \simeq \oplus_{i, j \ge 0} H^i ( X, \wedge^j
   T_X^{\ast} ) [ j - i ] . $$
The RHS of the last formula is the Hodge cohomology of $X$. One can consider the
hypercohomology ${\mathbb H}^{\bullet} ( X, \Omega^{\bullet}_X [ [ u ] ] / u^n
\Omega^{\bullet}_X [ [ u ] ] )$ equipped with the differential $ud_{dR}$. Then
the classical Hodge theory ensures degeneration of the corresponding spectral
sequence. In particular the hypercohomology is a flat $\K [ u ] / ( u^n
)$-module for any $n \ge 1$. Usual de Rham cohomology $H^{\bullet}_{dR} ( X )$
is isomorphic to the generic fiber of the corresponding flat vector bundle
over the formal line ${\bf A}^1_{form}[-2]$,
while the fiber at $u = 0$ is isomorphic to the Hodge cohomology
$H^{\bullet}_{Hodge} ( X ) = \oplus_{i, j \ge 0} H^i ( X, \wedge^j T_X^{\ast}
) [ j - i ]$. In order to make a connection with the ``abstract" theory
of the Section \ref{main conjecture} we remark that $H^{\bullet}_{dR} ( X )$ is
isomorphic (as a $\Z/2$-graded vector space) to the periodic cyclic homology $HP_{\bullet} ( A )$ while
$H_{\bullet} ( A, A )$ is isomorphic to $H^{\bullet}_{Hodge} ( X )$.

\section{$\A$-algebras with scalar product}\label{scalar product}

\subsection{Main definitions}\label{definition of algebra with scalar product}
Let $( X, pt, Q )$ be a finite-dimensional formal pointed dg-manifold over a field $\K$ of characteristic zero.

\begin{defn} A symplectic structure of
degree $N\in \Z$ on $X$ is given by a cyclic closed $2$-form $\omega$ of degree $N$
such that its restriction to the tangent space $T_{pt} X$ is non-degenerate.
\end{defn}

One has the following non-commutative analog of the Darboux lemma.

\begin{prp}\label{Darboux lemma}
Symplectic form $\omega$ has constant coefficients in some affine coordinates at the point $pt$.
\end{prp}

{\it Proof.} Let us choose an affine structure at the marked point
and write down $\omega=\omega_0+ \omega_1+\omega_2+....$, where
$\omega_l=\sum_{i,j}c_{ij}(x)dx_i\otimes dx_j$ and $c_{ij}(x)$
is homogeneous of degree $l$
(in particular, $\omega_0$ has constant coefficients).
Next we observe that the following lemma holds.

\begin{lmm}\label{high order Darboux lemma} Let $\omega=\omega_0+r$, where $r=\omega_l+\omega_{l+1}+..., l\ge 1$.
Then there is a change of affine coordinates
$x_i\mapsto x_i+O(x^{l+1})$  which transforms $\omega$ into
$\omega_0+\omega_{l+1}+...$.
\end{lmm}

Lemma implies the Proposition, since we can make an infinite product
of the above changes of variables (it is a well-defined infinite series).
The resulting automorphism of the formal neighborhood of $x_0$ transforms
$\omega$ into $\omega_0$.

{\it Proof of the lemma}. We have $d_{cycl}\omega_j=0$ for all $j\ge l$ (see Section \ref{higher order cyclic forms}  for
the notation).
The change of variables is determined by a vector field $v=(v_1,...,v_n)$
such that $v(x_0)=0$. Namely, $x_i\mapsto x_i-v_i, 1\le i\le n$.
Moreover, we will be looking for a vector field
such that $v_i=O(x^{l+1})$ for all $i$.

We have $Lie_v(\omega)=d(i_v\omega_0)+d(i_vr)$. Since $d\omega_l=0$
we have $\omega_l=d\alpha_{l+1}$ for some form $\alpha_{l+1}=O(x^{l+1})$ in the obvious notation. This follows from the formal Poincare lemma.
Therefore in order to kill the term with $\omega_l$ we need to solve
the equation $d\alpha_{l+1}=d(i_v\omega_0)$. It suffices to solve the
equation $\alpha_{l+1}=i_v\omega_0$. Since $\omega_0$ is non-degenerate,
there exists a unique vector field $v=O(x^{l+1})$ solving last equation.
This proves the lemma. $\blacksquare$

\begin{defn} \label{scalar product}Let $(X,pt,Q,\omega)$ be a non-commutative formal pointed symplectic dg-manifold. A scalar product of degree $N$ on the
$\A$-algebra $A=T_{pt}X[-1]$ is given by a choice of affine coordinates
at $pt$ such that the $\omega$ becomes constant and gives rise
to a non-degenerate bilinear form $A\otimes A\to \K[-N]$.

\end{defn}

\begin{rmk}\label{action functional}
Notice that since $Lie_Q ( \omega ) = 0$ there exists
a cyclic function $S\in \Omega_{cycl}^0(X)$ such that $i_Q \omega = dS$ and $\{ S, S \} = 0$
(here the Poisson bracket corresponds to the symplectic form $\omega$).
It follows  that the deformation
theory of a non-unital $\A$-algebra $A$ with the scalar product is controlled
by the DGLA $\Omega^0_{cycl} ( X )$ equipped with the differential $\{ S,
\bullet \}$.
\end{rmk}

We can restate the above definition in algebraic terms. Let $A$ be a
finite-dimensional  $\A$-algebra, which carries a non-degenerate symmetric bilinear form $(, )$ of
degree $N$. This means that for any two elements $a, b \in A$
such that $deg (a ) + deg ( b ) = N$ we are given a number $( a, b ) \in \K$ such that:

1) for any collection of elements $a_1, ..., a_{n + 1} \in A$ the expression
$( m_n ( a_1, ..., a_n ), a_{n + 1} )$ is cyclically symmetric in the graded
sense (i.e. it satisfies the Koszul rule of signs with respect to the cyclic
permutation of arguments);

2) bilinear form $(\bullet,\bullet)$ is non-degenerate.

In this case we will say that $A$ is an {\it $\A$-algebra with the scalar
product} of degree $N$.

\subsection{Calabi-Yau structure}\label{CY structure}

The above definition requires $A$ to be finite-dimensional.
We can relax this condition
requesting instead that $A$ is compact. As a result we will arrive to a
homological version of the notion of scalar product. More precisely,
assume that $A$ is weakly unital compact $\A$-algebra. Let
$CC_{\bullet}^{mod}(A)=(CC_{\bullet}^{mod}(A,A)[u^{-1}], b+uB)$ be the cyclic complex of $A$.
Let us choose a cohomology class $[\varphi]\in H^{\bullet}(CC_{\bullet}^{mod}(A))^{\ast}$
of degree $N$.
Since the complex $(A,m_1)$ is a subcomplex of
$C_{\bullet}^{mod}(A,A)\subset CC_{\bullet}^{mod}(A)$ we see that $[\varphi]$
defines a linear functional $Tr_{[\varphi]}: H^{\bullet}(A)\to \K[-N]$.

\begin{defn} \label{homologically non-degenerate functionals}We say that $[\varphi]$ is homologically non-degenerate if the bilinear
form of degree $N$ on $H^{\bullet}(A)$ given by $(a,b)\mapsto Tr_{[\varphi]}(ab)$ is
non-degenerate.

\end{defn}

Notice that the above bilinear form defines a symmetric scalar product
of degree $N$ on $H^{\bullet}(A)$ .

\begin{thm} \label{from functionals to scalar products}For a weakly unital compact $\A$-algebra $A$ a
homologically non-degenerate cohomology class $[\varphi]$
gives rise to a class of isomorphisms of non-degenerate scalar products on a minimal model of $A$.
\end{thm}

{\it Proof.} The proof will consist of several steps.
Since $char\,\K=0$ the complex $(CC_{\bullet}^{mod}(A))^{\ast}$
is quasi-isomorphic to $(\Omega^0_{cycl}(X)/\K, Lie_Q)$.

\begin{lmm} Complex $(\Omega^{2,cl}_{cycl}(X), Lie_Q)$ is quasi-isomorphic
to the complex $(\Omega^0_{cycl}(X)/\K, Lie_Q)$.\footnote{See also Proposition 5.5.1 from [Gi3].}

\end{lmm}

{\it Proof.} Notice that as a complex $(\Omega^{2,cl}_{cycl}(X), Lie_Q)$
is isomorphic to the complex $\Omega^1_{cycl}(X)/d_{cycl}\,\Omega_{cycl}^0(X)$.
The latter is quasi-isomorphic to $[{\cal O}(X),{\cal O}(X)]_{top}$ via
$a\,db\mapsto [a,b]$ (recall that $[{\cal O}(X),{\cal O}(X)]_{top}$
denotes the topological closure of the commutator).

By definition $\Omega^0_{cycl}(X)={\cal O}(X)/[{\cal O}(X),{\cal O}(X)]_{top}$.
We know that ${\cal O}(X)/\K$ is acyclic, hence $\Omega^0_{cycl}(X)/\K$
is quasi-isomorphic to $[{\cal O}(X),{\cal O}(X)]_{top}$. Hence the complex
$(\Omega^{2,cl}_{cycl}(X), Lie_Q)$ is quasi-isomorphic to
$(\Omega^0_{cycl}(X)/\K, Lie_Q)$. $\blacksquare$

As a corollary we obtain an isomorphism of cohomology groups
$H^{\bullet}(\Omega^{2,cl}_{cycl}(X))\simeq H^{\bullet}(\Omega^0_{cycl}(X)/\K)$.
Having a non-degenerate cohomology class
$[\varphi]\in H^{\bullet}(CC_{\bullet}^{mod}(A))^{\ast}\simeq H^{\bullet}(\Omega^{2,cl}_{cycl}(X), Lie_Q)$
as above, we can choose
its representative $\omega\in \Omega^{2,cl}_{cycl}(X)$,
$Lie_Q\omega=0$. Let us consider
$\omega(x_0)$. It can be described pure algebraically such as follows.
Notice that there is a natural projection
$H^{\bullet}(\Omega^0_{cycl}(X)/\K)\to (A/[A,A])^{\ast}$ which
corresponds to the taking the first Taylor coefficient of the
cyclic function. Then the above evaluation $\omega(x_0)$
is the image of $\varphi(x_0)$ under the natural map $(A/[A,A])^{\ast}\to (Sym^2(A))^{\ast}$
which assigns to a linear functional  $l$ the bilinear
form $l(ab)$.

We claim that the total map
$H^{\bullet}(\Omega^{2,cl}_{cycl}(X))\to(Sym^2(A))^{\ast}$ is the same
as the evaluation at $x_0$ of the closed cyclic $2$-form. Equivalently,
we claim that $\omega(x_0)(a,b)=Tr_{\varphi}(ab)$. Indeed, if
$f\in \Omega_{cycl}^0(X)/\K$ is the cyclic function corresponding
to $\omega$ then we can write $f=\sum_ia_ix_i+O(x^2)$.
Therefore $Lie_Q(f)=\sum_{l,i,j}a_ic^{ij}_l[x_i,x_j]+O(x^3)$,
where $c^{ij}_l$ are structure constants of ${\cal O}(X)$.
Dualizing we obtain the claim.

\begin{prp} Let $\omega_1$ and $\omega_2$ be two symplectic structures
on the finite-dimensional formal pointed minimal dg-manifold $(X,pt,Q)$ such that
$[\omega_1]=[\omega_2]$ in the cohomology of the complex
$(\Omega_{cycl}^{2,cl}(X),Lie_Q)$ consisting of closed cyclic $2$-forms.
Then there exists a change of coordinates at $x_0$ preserving $Q$ which transforms $\omega_1$ into $\omega_2$.

\end{prp}

\begin{cor}\label{symplectic forms on minimal model} Let $(X,pt,Q)$ be a (possibly infinite-dimensional) formal
pointed dg-manifold endowed with a (possibly degenerate) closed cyclic $2$-form
$\omega$. Assume that the tangent cohomology $H^0(T_{pt}X)$ is finite-dimensional
and $\omega$ induces a non-degenerate pairing on it. Then on the minimal model
of $(X,pt,Q)$ we have a canonical isomorphism class of symplectic forms
modulo the action of the group $Aut(X,pt,Q)$.

\end{cor}

{\it Proof.}
Let $M$ be a (finite-dimensional) minimal model of $A$. Choosing a
cohomology class $[\varphi]$ as above we
obtain a non-degenerate bilinear form
on $M$, which is the restriction $\omega(x_0)$ of a representative
$\omega\in \Omega^{2,cl}(X)$. By construction this scalar product
depends on $\omega$. We would like to show that in fact it depends
on the cohomology class of $\omega$, i.e. on $\varphi$ only.
This is the corollary of the following result.

\begin{lmm} \label{a vector field}Let $\omega_1=\omega+Lie_Q(d\alpha)$. Then there exists
a vector field $v$ such that $v(x_0)=0, [v,Q]=0$ and
$Lie_v(\omega)=Lie_Q(d\alpha)$.

\end{lmm}
{\it Proof.} As in the proof of Darboux lemma we need to find
a vector field $v$, satisfying the condition
$di_v(\omega)=Lie_Q(d\alpha)$. Let $\beta=Lie_Q(\alpha)$.
Then $d\beta=dLie_Q(\alpha)=0$. Since $\omega$ is non-degenerate
we can find $v$ satisfying the conditions of the Proposition
and such that $di_v(\omega)=Lie_Q(d\alpha)$. Using this $v$ we
can change affine coordinates transforming $\omega+Lie_Q(d\alpha)$
back to $\omega$.
This concludes the proof of the Proposition and the Theorem \ref{from functionals to scalar products}.$\blacksquare$

Presumably the above construction is equivalent to the one given in [Kaj].
We will sometimes call the cohomology class $[\varphi]$ a {\it Calabi-Yau structure}
on $A$ (or on the corresponding non-commutative formal pointed dg-manifold $X$). The following example illustrates
the relation to geometry.

\begin{exa} \label{geometric example}Let $X$ be a complex Calabi-Yau manifold of dimension $n$.
Then it carries a nowhere vanishing holomorphic $n$-form $vol$.
Let us fix a holomorphic vector bundle $E$ and consider a dg-algebra $A=\Omega^{0,\ast}(X,End(E))$ of Dolbeault $(0,p)$-forms with values
in $End(E)$. This dg-algebra carries a linear functional
$a\mapsto \int_XTr(a)\wedge vol$. One can check that this is a cyclic
cocycle which defines a non-degenerate pairing on $H^{\bullet}(A)$
in the way described above.

\end{exa}

There is another approach to  Calabi-Yau structures in
the case when $A$ is {\it homologically smooth}.
Namely, we say that $A$ carries a {\it Calabi-Yau structure of dimension
$N$} if $A^!\simeq A[N]$
(recall that $A^!$ is the $A - A$-bimodule $Hom_{A - mod - A} ( A, A \otimes A )$ introduced in Section \ref{homological smoothness}). Then we expect the following
conjecture to be true.

\begin{conj}
If $A$ is a homologically smooth compact finite-dimensional $\A$-algebra  then the existence of  a non-degenerate cohomology class $[\varphi]$
of degree $dim\,A$ is equivalent to the condition $A^! \simeq A [ dim \,A ]$.
\end{conj}

If $A$ is the dg-algebra of endomorphisms
of a generator of $D^b(Coh(X))$ ($X$ is Calabi-Yau) then
the above conjecture holds trivially.

Finally, we would like to illustrate the relationship of the non-commutative symplectic
geometry discussed above with the commutative symplectic geometry
of certain spaces of representations \footnote{It goes back to [Ko2], and since that time
has been discussed in many papers, see e.g. [Gi2].}. More generally we would like
to associate with $X=Spc(T(A[1]))$ a collection of formal
algebraic varieties, so that some ``non-commutative"
geometric structure on $X$
becomes a collection of
compatible ``commutative" structures on
formal manifolds ${\cal M}(X,n):=\widehat{Rep}_0({\cal O}(X), Mat_n(\K))$,
where $Mat_n(\K)$ is the associative
algebra of $n\times n$ matrices over $\K$, ${\cal O}(X)$ is the
algebra of functions on $X$ and $\widehat{Rep}_0(...)$ means the
formal completion at the trivial representation.
In other words, we would like to define a collection
of compatible geometric structure on ``$Mat_n(\K)$-points" of
the formal manifold $X$.
In the case of symplectic structure this philosophy is
illustrated by the following result.

\begin{thm} \label{symplectic structure on representation spaces}Let $X$ be a non-commutative formal symplectic manifold in $Vect_\K$.
Then it defines a collection of symplectic structures on all
manifolds ${\cal M}(X,n), n\ge 1$.

\end{thm}

{\it Proof.} Let ${\cal O}(X)=A, {\cal O}({\cal M}(X,n))=B$.
Then we can choose  isomorphisms
$A\simeq \K\langle\langle x_1,...,x_m \rangle\rangle$ and
$B\simeq
\K\langle\langle x_1^{{\alpha},{\beta}},..., x_m^{{\alpha},{\beta}}\rangle\rangle$,
where $1\le \alpha,\beta \le n$. To any $a\in A$ we can assign
$\widehat{a}\in B\otimes Mat_n(\K)$ such that:

$$\hat{x}_i=\sum_{\alpha,\beta} x_i^{{\alpha},{\beta}}\otimes e_{\alpha,\beta},$$
where $e_{\alpha,\beta}$ is the $n\times n$
matrix with the only non-trivial element equal to $1$ on the intersection
of $\alpha$-th line and $\beta$-th column.
The above formulas define an algebra homomorphism. Composing it with the map
$id_B\otimes Tr_{Mat_n(\K)}$ we get a linear map
${\cal O}_{cycl}(X)\to {\cal O}({\cal M}(X,n))$. Indeed the
closure of the commutator $[A,A]$ is mapped to zero.
Similarly, we have a morphism of complexes
$\Omega_{cycl}^{\bullet}(X)\to \Omega^{\bullet}({\cal M}(X,n))$, such
that
$$dx_i\mapsto \sum_{\alpha,\beta} dx_i^{{\alpha},{\beta}}e_{\alpha,\beta}.$$
Clearly, continuous derivations of $A$ (i.e. vector fields on $X$)
give rise to the vector fields on ${\cal M}(X,n)$.

Finally, one can see that a non-degenerate cyclic $2$-form $\omega$
is mapped to the tensor product of a non-degenerate $2$-form on
${\cal M}(X,n)$ and a nondegenerate $2$-form $Tr(XY)$ on
$Mat_n(\K)$. Therefore a symplectic form on $X$ gives rise
to a symplectic form on ${\cal M}(X,n), n\ge 1$. $\blacksquare$

\section{Hochschild complexes as algebras over operads and  PROPs}\label{Hochschild complexes and operads}

Let $A$ be a strictly unital $\A$-algebra over a field $\K$ of characteristic zero.
In this section we are going to describe a colored dg-operad $P$ such that the pair $(
C^{\bullet} ( A, A ), C_{\bullet} ( A, A ) )$ is an algebra over this operad.
More precisely, we are going to describe $\Z$-graded $\K$-vector spaces $A ( n,
m )$ and $B ( n, m )$, $n, m \ge 0$ which are components of the colored operad
such that $B ( n, m ) \ne 0$ for $m = 1$ only and $A ( n, m ) \ne 0$ for $m =
0$ only together with the colored operad structure and the action

a) $A ( n, 0 ) \otimes ( C^{\bullet} ( A, A ) )^{\otimes n} \to C^{\bullet} (
A, A ),$

b) $B ( n, 1 ) \otimes ( C^{\bullet} ( A, A ) )^{\otimes n} \otimes
C_{\bullet} ( A, A ) \to C_{\bullet} ( A, A ) .$

Then, assuming that $A$ carries a non-degenerate scalar product, we are going to describe
a PROP $R$ associated with moduli spaces of Riemannian surfaces and a structure
of $R$-algebra on $C_{\bullet} ( A, A )$.

\subsection{Configuration spaces of discs}\label{configurations of discs}

We start with the spaces $A ( n, 0 )$. They are chain complexes. The complex
$A ( n, 0 )$ coincides with the complex $M_n$ of the minimal operad $M = ( M_n
)_{n \ge 0}$ described in [KoSo1], Section 5. Without going into details which
can be found in loc. cit. we recall main facts  about the operad $M$.
A basis of $M_n$ as a $k$-vector space is formed by $n$-labeled planar trees
(such trees have internal vertices labeled by the set $\{ 1, ..., n \}$ as
well as other internal vertices which are non-labeled and each has the valency
at least $3$).

We can depict $n$-labeled trees such as follows

\vspace{3mm}
{\includegraphics{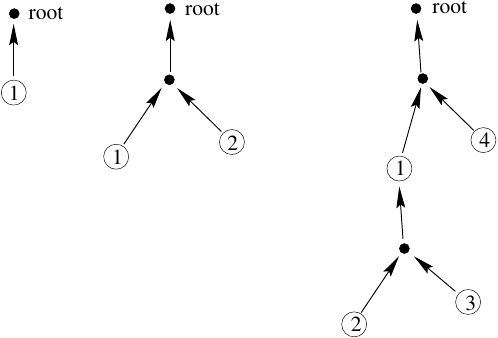}}

\vspace{3mm}

Labeled vertices are depicted as circles with numbers inscribed, non-labeled
vertices are depicted as black vertices. In this way we obtain a graded operad
$M$ with the total degree of the basis element corresponding to a tree $T$
equal to
$$ deg ( T ) = \sum_{v \in V_{lab} ( T )} (1 - |v| ) + \sum_{v \in V_{nonl} ( T
   )} ( 3 - |v| ) $$
where $V_{lab} ( T )$ and $V_{nonl} ( T )$ denote the sets of labeled and
non-labeled vertices respectively , and $|v|$ is the valency of the vertex $v$, i.e. the cardinality of the set of edges attached to $v$.

The notion of an {\it angle} between two edges incoming in a vertex is
illustrated in the following figure (angles are marked by asteriscs).

\vspace{3mm}

\includegraphics{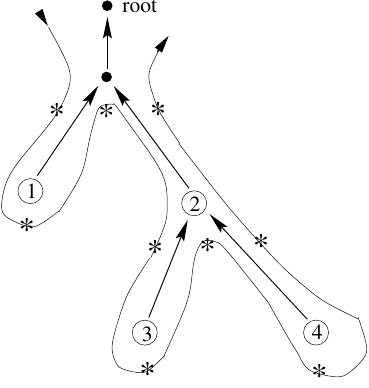}
\vspace{3mm}

Operadic composition and the differential are described in [KoSo1], sections
5.2, 5.3. We borrow from there the following figure which illustrates the
operadic composition of generators corresponding to labeled trees $T_1$ and
$T_2$.

\vspace{3mm}
\includegraphics{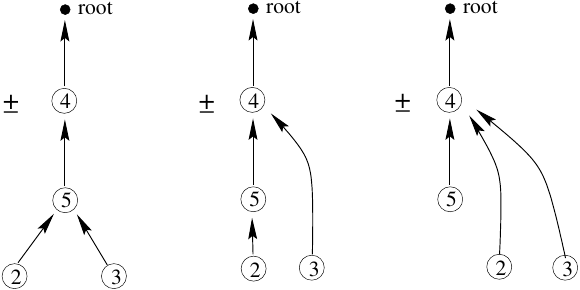}

\vspace{3mm}

{\centerline{{\epsfbox{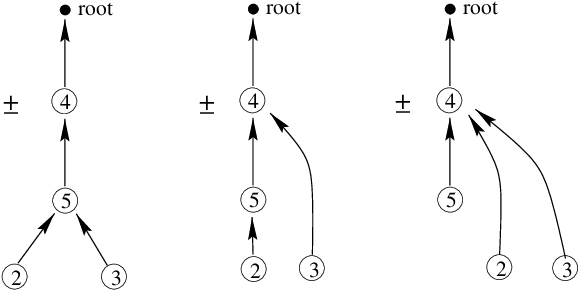}} }}

\vspace{3mm}

Informally speaking, the operadic gluing of $T_2$ to $T_1$ at an internal
vertex $v$ of $T_1$ is obtained by:

a) Removing from $T_1$ the vertex $v$ together with all incoming edges and
vertices.

b) Gluing $T_2$ to $v$ (with the root vertex removed from $T_2$). Then

c) Inserting removed vertices and edges of $T_1$ in all angles between
incoming edges to the new vertex $v_{new}$.

d) Taking the sum (with appropriate signs) over all possible inserting of
edges in c).

The differential $d_M$ is a sum of the ``local" differentials $d_v$,
where $v$ runs through the set of all internal vertices. Each $d_v$ inserts a
new edge into the set of edges attached to $v$. The following figure borrowed
from [KoSo1] illustrates the difference between labeled (white) and
non-labeled (black) vertices.

\vspace{3mm}
\includegraphics{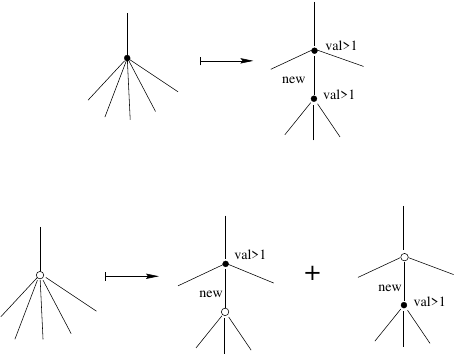}
\vspace{3mm}

In this way we make $M$ into a dg-operad. It was proved in [KoSo1], that $M$
is quasi-isomorphic to the dg-operad $Chains ( FM_2 )$ of singular chains on
the Fulton-Macpherson operad $FM_2$. The latter consists of the compactified
moduli spaces of configurations of points in $\R^2$ (see e.g. [KoSo1], Section
7.2 for a description). It was also proved in [KoSo1] that $C^{\bullet} ( A, A
)$ is an algebra over the operad $M$ (Deligne's conjecture follows from this
fact). The operad $FM_2$ is homotopy equivalent to the famous operad $C_2 = (
C_2 ( n ) )_{n \ge 0}$ of $2$-dimensional discs (little disc operad). Thus
$C^{\bullet} ( A, A )$ is an algebra (in the homotopy sense) over the operad $Chains ( C_2 )$.

\subsection{Configurations of points on the cylinder}\label{configuration of points on cylinder}

Let $\Sigma=S^1\times [0,1]$ denotes the standard cylinder.

Let us denote by $S(n)$ the set of isotopy classes of the following
graphs $\Gamma\subset \Sigma$:

a) every graph $\Gamma$ is a forest
(i.e. disjoint union of finitely many trees
$\Gamma=\sqcup_iT_i$);

b) the set of vertices $V(\Gamma)$ is decomposed into the union
$V_{\partial\Sigma}\sqcup V_{lab}\sqcup V_{nonl}\sqcup V_1$ of four
sets with the following properties:

b1) the  set $V_{\partial\Sigma}$ is the union $\{in\}\cup \{out\}\cup V_{out}$
of three sets of points which belong to the boundary
$\partial\Sigma$ of the cylinder. The set $\{in\}$ consists of one marked point
which belongs to the boundary circle $S^1\times \{1\}$ while
the set $\{out\}$ consists of one marked point which belongs to the boundary
circle $S^1\times \{0\}$. The set $V_{out}$ consists of a finitely many
unlableled points on the boundary circle $S^1\times \{0\}$;

b2) the set $V_{lab}$ consists of $n$ labeled points which belong to the
surface $S^1\times (0,1)$ of the cylinder;

b3) the set $V_{nonl}$ consists of a finitely many non-labeled points
which belong to the
surface $S^1\times (0,1)$ of the cylinder;

b4) the set $V_1$ is either empty or consists of only one element denoted
by ${\bf 1}\in S^1\times (0,1)$ and called  {\it special} vertex;

c) the following conditions on the valencies of vertices are imposed:

c1) the valency of the vertex $out$ is less or equal than $1$;

c2) the valency of each vertex from the set $V_{\partial\Sigma}\setminus V_{out}$
is equal to $1$;

c3) the valency of each
vertex from $V_{lab}$ is at least $1$;

c4) the valency of each vertex
from $V_{nonl}$ is at least $3$;

c5) if the set $V_1$ is non-empty then the valency of the special
vertex is equal to $1$. In this case
the only outcoming edge connects ${\bf 1}$ with the vertex $out$.

d) Every tree $T_i$ from the forest $\Gamma$ has its root vertex in the set
$V_{\partial\Sigma}$.

e) We orient each tree $T_i$ down to its root vertex.










\vspace{3mm}

\includegraphics{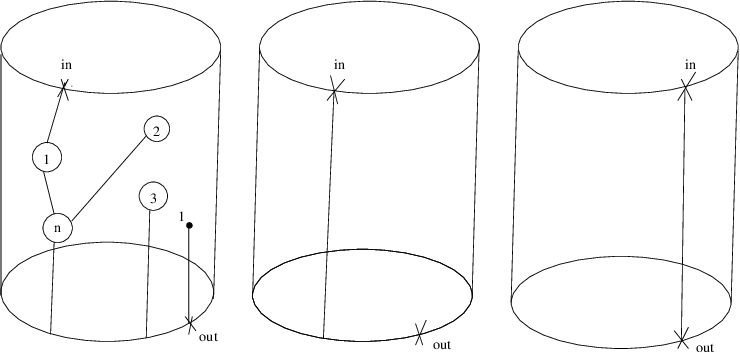}

\vspace{3mm}

\begin{rmk}\label{cell decomposition} Let us consider the configuration space $X_n, n\ge 0$ which
consists of (modulo ${\bf C}^{\ast}$-dilation) equivalence
classes of $n$ points on ${\bf CP}^1\setminus \{0,\infty\}$ together
with two direction lines at the tangent spaces at the points $0$ and
$\infty$. One-point compactification $\widehat{X}_n$ admits
a cell decomposition with cells (except of the point $\widehat{X}_n\setminus X_n$) parametrized by elements
of the set $S(n)$. This can be proved with the help of Strebel
differentials (cf. [KoSo1], Section 5.5).
\end{rmk}

Previous remark is related to the following description
of the sets $S(n)$ (it will be used later in the paper).
Let us
contract both circles of the boundary $\partial\Sigma$ into points. In
this way we obtain a tree on the sphere. Points become vertices of the tree and lines
outcoming from the points become edges. There are two vertices marked by $in$ and
$out$ (placed at the north and south poles respectively). We  orient the
tree towards to the vertex $out$. An additional structure consists of:

a) Marked edge outcoming from $in$ (it corresponds to the edge outcoming from
$in$).

b) Either a marked edge incoming to $out$ (there was an edge incoming to
$out$ which connected it with a vertex not marked by $\bf{1}$) or an angle between
two edges incoming to $out$ (all edges which have one of the endpoint
vertices on the bottom circle become after
contracting it to a point the edges incoming to $out$, and if there was an
edge connecting a point marked by $\bf{1}$ with $out$, we mark the angle
between edges containing this line).

The reader notices that the star of the vertex $out$ can be identified with a
regular $k$-gon, where $k$ is the number of incoming to $out$ edges. For this
$k$-gon we have either a marked point on an edge (case a) above) or a marked
angle with the vertex in $out$ (case b) above).

\subsection{Generalization of Deligne's conjecture}\label{generalized Deligne\rq{}s conjecture}

The definition of the operadic space $B(n,1)$ will be clear from
the description of its action on the Hochschild chain complex.
The space $B(n,1)$ will have a basis parametrized by
elements of the set $S(n)$ described in the previous subsection.
Let us describe the action of a generator of $B ( n, 1 )$ on a pair $(
\gamma_1 \otimes ... \otimes \gamma_n, \beta )$, where $\gamma_1 \otimes ...
\otimes \gamma_n \in C^{\bullet} ( A, A )^{\otimes n}$ and $\beta = a_0
\otimes a_1 \otimes ... \otimes a_l \in C_l ( A, A )$. We attach elements $a_0,
a_1, ..., a_l$ to points on $S^1\times \{1\}$ in a cyclic order,
such that  $a_0$ is attached to the
point $in$. We attach $\gamma_i$ to the $i$th numbered point on the surface $S^1\times (0,1)$.
Then we draw disjoint continuous segments (in all possible ways, considering
pictures up to an isotopy) starting from each point marked by some element $a_i$
and oriented downstairs, with the requirements  as above, with the only
modification that we allow an arbitrary number of points on $S^1\times \{1\}$. We
attach higher multiplications $m_j$ to all non-numbered vertices, so that $j$
is equal to the incoming valency of the vertex. Reading from the top to the
bottom and composing $\gamma_i$ and $m_j$ we obtain (on the bottom circle) an element $
b_0\otimes ...\otimes b_m\in C_{\bullet} ( A, A )$ with $b_0$ attached to the vertex $out$.
If the special vertex ${\bf 1}$ is present then we set $b_0=1$.
This gives the desired action.


\includegraphics{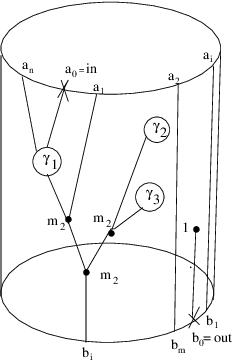}

{\vspace{3mm}}

Composition of the
operations in $B(n,1)$  corresponds to the gluing of the cylinders such that the point
$out$ of the top cylinder is identified with the point $in$ of the bottom
cylinder. \footnote{Technically we should resize vertically the glued cylinder, so that it is identified with $\Sigma=S^1\times [0,1]$. Alternatively
we can introduce an additional parameter $h\in \R_{>0}$ and  consider cylinders $\Sigma_h=S^1\times [0,h]$. Then the composition of the operations depicted 
on $\Sigma_{h_1}$ and $\Sigma_{h_2}$ amounts to the operation depicted on $\Sigma_{h_1+h_2}$. We adopt this convention below when speaking about the colored
operad.  In that case it makes more sense since the case $h=0$ has to be included.}
 If after the gluing there is a line from the point marked ${\bf 1}$
on the top cylinder which does not end at the point $out$ of the bottom
cylinder, we will declare such a composition to be equal to zero.

Let us now consider a topological colored operad $C_2^{col} = ( C_2^{col} ( n, m )
)_{n, m \ge 0}$ with two colors such that $C_2^{col} ( n, m ) \ne \emptyset$
only if $m = 0, 1$, and

a) In the case $m = 0$ it is the little disc operad.

b) In the case $m = 1$ $C_2^{col} ( n, 1 )$ is the moduli space (modulo
rotations) of the configurations of $n \ge 1$ discs on the cyliner $S^1 \times
[ 0, h ]$ $h \ge 0$, and two marked points on the boundary of the cylinder. We
also add the degenerate circle of configurations $n = 0, h = 0$. The topological
space $C_2^{col}(n,1)$ is homotopically equivalent to the configuration
space $X_n$ described in the previous subsection.

Let $Chains ( C_2^{col} )$ be the  colored operad of singular
chains on $C_2^{col}$. Then, similarly to [KoSo1], Section 7, one proves (using the explicit
action of the colored operad $P = ( A ( n, m ), B ( n, m ) )_{n, m \ge 0}$
described above) the following result.

\begin{thm} \label{algebra over colored operad}Let $A$ be a unital $\A$-algebra. Then the pair $( C^{\bullet} ( A, A
), C_{\bullet} ( A, A ) )$ is an algebra over the colored operad $Chains (
C_2^{col} )$ (which is quasi-isomorphic to $P$) such that for $h = 0, n = 0$
and coinciding points $in = out$, the corresponding operation is the
identity.
\end{thm}

\begin{rmk} \label{abstract calculus}The above Theorem \ref{algebra over colored operad} generalizes Deligne's conjecture (see e.g. [KoSo1]). It
is related to the abstract calculus associated with $A$ (see [T], [TaT]). The
reader also notices that for $h = 0, n = 0$ we have the moduli space of two
points on the circle. It is homeomorphic to $S^1$. Thus we have an action of
$S^1$ on $C_{\bullet} ( A, A )$. This action gives rise to the Connes
differential $B$.
\end{rmk}

Similarly to the case of little disc operad, one can prove the following
result.

\begin{prp} \label{colored operad is formal}The colored operad $C_2^{col}$ is formal, i.e. it is quasi-isomorphic
to its homology colored operad.
\end{prp}

If $A$ is non-unital we can consider the direct sum
$A_1=A\oplus \K$ and make it into a unital $\A$-algebra.
The {\it reduced} Hochschild chain complex of $A_1$ is defined
as $C_{\bullet}^{red}(A_1,A_1)=\oplus_{n\ge 0} A_1\otimes ((A_1/\K)[1])^{\otimes n}$ with the same differential as in
the unital case. One defines the reduced Hochschild cochain
complex $C^{\bullet}_{red}(A_1,A_1)$ similarly.
We define the {\it modified} Hochschild chain complex $C_{\bullet}^{mod}(A,A)$ from
the following isomorphism of complexes
$C_{\bullet}^{red}(A_1,A_1)\simeq C_{\bullet}^{mod}(A,A)\oplus \K$.
Similarly, we define the modified
Hochschild cochain complex from the decomposition
$C^{\bullet}_{red}(A_1,A_1)\simeq C^{\bullet}_{mod}(A,A)\oplus \K$.
Then, similarly to the Theorem \ref{algebra over colored operad} one proves the following result.

\begin{prp}\label{algebra over extended colored operad}
The pair $( C_{\bullet}^{mod} ( A, A), C^{\bullet}_{mod} ( A, A ) )$
is an algebra over the colored operad which is an extension of $Chains (C_2^{col} )$
by null-ary operations on Hochschild chain and cochain complexes,
which correspond to the unit in $A$, and
such that for $h = 0, n = 0$
and coinciding points $in = out$, the corresponding operation is the
identity.

\end{prp}

\subsection{Remark about Gauss-Manin connection}\label{GM connection}

Let $R=\K[[t_1,...,t_n]]$ be the algebra of formal series, and $A$ be an $R$-flat
$\A$-algebra. Then the (modified) negative cyclic complex $CC_{\bullet}^{-,mod}(A)=(C_{\bullet}(A,A)[[u]],b+uB)$
is an $R[[u]]$-module. It follows from the existense of the Gauss-Manin
connection (see [Get])  that the cohomology $HC_{\bullet}^{-,mod}(A)$ of this complex  is in fact a module over the ring
$$D_R(A):=\K[[t_1,...,t_n,u]][u\partial/\partial t_1,...,u\partial/\partial t_n].$$
Inedeed, if $\nabla$ is the Gauss-Manin connection from [Get] then $u\partial/\partial t_i$
acts on the cohomology as $u\nabla_{\partial/\partial t_i},1\le i\le n$.

The conjecture below give another interpretation of  this module structure.
Let ${\g}=C^{\bullet}(A,A)[1]$ be the DGLA associated with the Hochschild cochain complex, and
$M:=C_{\bullet}^{-,mod}(A)$. We define a DGLA $\hat{\g}$ which is
the crossproduct $(\g\otimes \K\langle \xi\rangle)\rtimes \K(\partial/\partial \xi)$, where
$deg\,\xi=+1$.

\begin{conj} \label{generalization of Deligne\rq{}s conjecture} There is a structure of an $L_{\infty}$-module on $M$ over $\hat{\g}$ which extends
the natural structure of a $\g$-module and such that $\partial/\partial \xi$ acts as Connes
differential $B$. Moreover this $\hat{\g}$-module structure should follow from the algebra structure over the colored operad $P$
described in Section \ref{generalized Deligne\rq{}s conjecture}.

\end{conj}

It looks plausible that Theorem \ref{algebra over colored operad} gives formulas for the Gauss-Manin connection
from [Get].

\subsection{Flat connections and the colored operad}\label{flat connection and operad}
We start with $\Z$-graded case.
Let us interpret the $\Z$-graded formal scheme $Spf(\K[[u]])$ as even formal line
equipped with the ${\bf G}_m$-action $u\mapsto \lambda^2 u$.
The space $HC_{\bullet}^{-,mod}(A)$ can be interpreted as a space of sections
of a
${\bf G}_m$-equivariant vector bundle $\xi_A$ over the formal line $Spf ( \K [ [ u ] ] )$
corresponding to the $\K [ [ u ] ]$-flat module $\varprojlim_n H^{\bullet} (
C_{\bullet}^{( n )} ( A, A ) )$. The action of ${\bf G}_m$ identifies fibers of this
vector bundle over $u \ne 0$. Thus we have a natural flat connection $\nabla$
on the
restriction of $\xi_A$ to the complement of the point $0$. This connection
has  pole of order one at $u=0$.

We are going to propose a different construction of the connection
$\nabla$ which works also in $\Z/2$-graded case. This connection will
have in general a pole of degree two at $u=0$. In particular we will have the following
result.

\begin{prp} \label{module structure on sections}The space of sections of the vector bundle $\xi_{{A}}$
can be endowed with a structure of a
$\K[[u]] [[u^2\partial/\partial u]]$-module.
\end{prp}

In fact we are going to give an explicit construction of the connection, which
is based on the action of the colored dg-operad $P$ discussed in Section \ref{generalized Deligne\rq{}s conjecture}.
More precisely, we will consider a certain  extension $P^{new}$ of $P$. Before presenting an explicit
formula, we will make few comments.

1. For any $\Z/2$-graded $\A$-algebra $A$ one can define  canonically a
$1$-parameter family of $\A$-algebras $A_{\lambda}, \lambda \in {\bf G}_m$,
such that $A_{\lambda}=A$ as a $\Z/2$-graded vector space and $m_n^{A_{\lambda}}=\lambda
m_n^A$.

2. For simplicity we will assume that $A$ is strictly unital. Otherwise we
will work with
the pair $( C_{\bullet}^{mod} ( A, A ), C^{\bullet}_{mod} ( A, A ) )$
of modified Hochschild complexes.

3. Let us consider an extension $P^{new}$ of the dg-operad $P$ allowing any
non-zero valency for a non-labeled (black) vertex (recall that in the definition of $P$ we
required that such a valency was at least three). All the formulas remain the
same. But the dg-operad $P^{new}$ is no longer formal. It contains a
dg-suboperad generated by trees with all vertices being non-labeled. Action of
this suboperad $P^{new}_{nonl}$ is responsible for the flat connection
discussed below.

4. In addition to the connection along the variable $u$ one has the Gauss-Manin
connection which acts along the fibers of $\xi_A$ (see Section \ref{GM connection}). Probably one
can write down an explicit formula for this connection using the action of the
colored operad $P^{new}$. In what follows are going to describe a connection
which presumably coincides with the Gauss-Manin connection.

Let us now consider a dg-algebra $\K [ B, \gamma_0, \gamma_2 ]$ which is
generated by the following operations of the colored dg-operad $P^{new}$:

a) Connes differential $B$ of degree $- 1$. It can be depicted such as
follows (cf. Section \ref{case of strictly unital algebras}):


\includegraphics{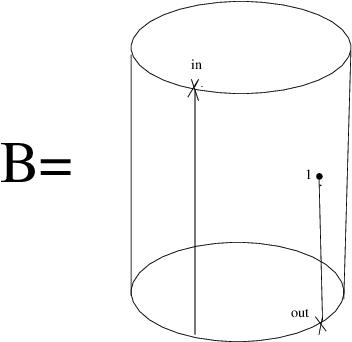}

{\vspace{3mm}}

b) Generator $\gamma_2$ of degree $2$, corresponding to the following figure:

{\vspace{3mm}}


\includegraphics{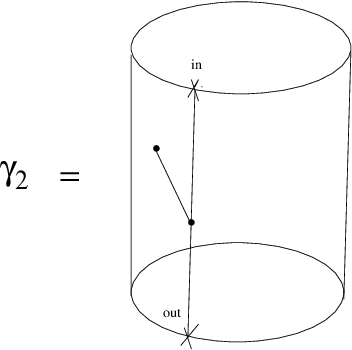}

{\vspace{3mm}}

c) Generator $\gamma_0$ of degree $0$, where $2 \gamma_0$ is depicted below:

{\vspace{3mm}}


\includegraphics{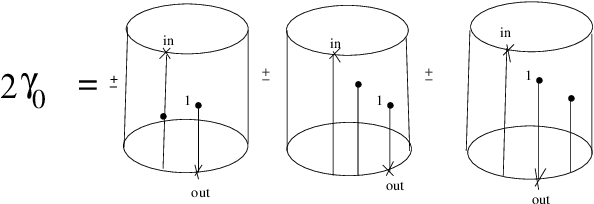}

{\vspace{3mm}}

\begin{prp}\label{identities}
The following identities hold in $P^{new}$:
$$ B^2 = dB = d \gamma_2 = 0, d \gamma_0 = [ B, \gamma_2 ], $$
$$ B \gamma_0 + \gamma_0 B : = [ B, \gamma_0 ]_+ = - B. $$
Here by $d$ we denote the Hochschild chain differential (previously it was
denoted by $b$).

\end{prp}
{\it Proof.} Let us prove that $[ B, \gamma_0 ] = - B$, leaving
the rest as an exercise to the reader. One has the following identities for
the compositions of operations in $P^{new}$: $B \gamma_0 = 0$, $\gamma_0 B =
B$. Let us check, for example, the last identity. Let us denote by $W$ the
first summand on the figure defining $2 \gamma_0$. Then $\gamma_0 B =
\frac{1}{2} WB$. The latter can be depicted in the following way:

{\vspace{3mm}}

\includegraphics{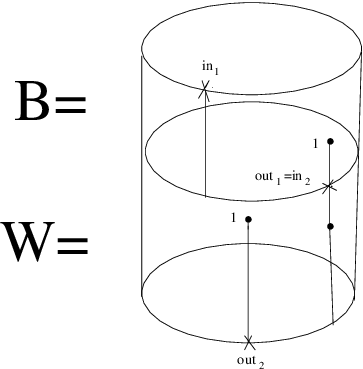}

{\vspace{3mm}}

It is easy  to see that this equals to $2 \cdot \frac{1}{2} B = B$. $\blacksquare$

\begin{cor} \label{Hochschild chains as module}Hochschild chain complex $C_{\bullet} ( A, A )$ is a dg-module over the
dg-algebra $\K [ B, \gamma_0, \gamma_2 ]$.
\end{cor}

Let us consider the truncated negative cyclic complex $( C_{\bullet} ( A, A )
[ [ u ] ] / ( u^n ), d_u = d + uB )$. We introduce a $\K$-linear map $\nabla$
of $C_{\bullet} ( A, A ) [ [ u ] ] / ( u^n )$ into itself such that
$\nabla_{u^2 \partial / \partial u} = u^2 \partial / \partial u - \gamma_2 + u
\gamma_0$. Then we have:

a) $[ \nabla_{u^2 \partial / \partial u}, d_u ] = 0 ;$

b) $[ \nabla_{u^2 \partial / \partial u}, u ] = u^2 .$

Let us denote by $V$ the unital dg-algebra generated by $\nabla_{u^2 \partial
/ \partial u}$ and $u$, subject to the relations a), b) and the relation $u^n
= 0$. From a) and b) one deduces the following result.

\begin{prp} \label{flat connection on cohomology}The complex $( C_{\bullet} ( A, A ) [ [ u ] ] / ( u^n ), d_u = d + uB
)$ is a $V$-module. Moreover, assuming the degeneration conjecture, we see
that the operator $\nabla_{u^2 \partial / \partial u}$ defines a flat
connection on the cohomology bundle corresponding to the module
$$H^{\bullet} ( C_{\bullet} ( A, A ) [ [ u] ] / ( u^n ), d_u )$$
such that the connection has the only singularity at $u = 0$ which is a
pole of second order.
\end{prp}

Taking the inverse limit over $n$ we see that $H^{\bullet} ( C_{\bullet} ( A,
A ) [ [ u ] ], d_u )$ gives rise to a vector bundle over ${\bf A}^1_{form}[-2]$
which carries a flat connection with the second order pole at $u = 0$. 
\begin{rmk}\label{difference between graded cases}
In the above discussion  one should notice a difference between $\Z$-graded and $\Z / 2$-graded
$\A$-algebras. It follows from the explicit formula for the connection
$\nabla$ that the coefficient in front of the second order pole is represented by
multiplication by the cocyle $( m_n )_{n \ge 1}\in C^{\bullet}(A,A)$. Taking the cohomology we see that
the coefficient is trivial in
$\Z$-graded case (because of the invariance with respect to the group action
$m_n \mapsto \lambda \hspace{0.25em} m_n$), but nontrivial in $\Z / 2$-graded
case. Therefore the order of the pole of $\nabla$ is equal to one for
$\Z$-graded $\A$-algebras and is equal to two for $\Z / 2$-graded
$\A$-algebras. We see that in  $\Z$-graded case the connection
along the variable $u$ comes from the
action of the group ${\bf G}_m$ on higher products $m_n$, while
in $\Z/2$-graded case it is more complicated.
\end{rmk}

\subsection{PROP of marked Riemann surfaces and TQFT associated with a Calabi-Yau category}\label{PROP}

Here we will describe a PROP naturally acting on the Hochschild
complexes of a finite-dimensional $\A$-algebra with the scalar
product of degree $N$.

Since we have a quasi-isomorphism of complexes
$$ C^{\bullet} ( A, A ) \simeq ( C_{\bullet} ( A, A ) )^{\ast} [ - N ]  $$
it suffices to consider the chain complex only.

We will assume that $A$ is either $\Z$-graded
(then $N$ is an integer) or $\Z/2$-graded (then $N\in \Z/2$).
We will present the results for non-unital $\A$-algebras. In this case
we will consider the modified Hochschild chain complex
$$C_{\bullet}^{mod}=\oplus_{n\ge 0}A\otimes (A[1])^{\otimes n}\bigoplus
\oplus_{n\ge 1}(A[1])^{\otimes n},$$
equipped with the Hochschild chain differential (see Section \ref{non-unital case}).

Our construction is summarized in i)-ii) below.

i) Let us consider the topological PROP
${\cal M}=({\cal M}(n,m))_{n,m\ge 0}$ consisting of moduli spaces of metrics on compact
oriented surfaces with bondary consisting of $n+m$ circles and some additional marking
(see precise definition below).

ii) Let $Chains({\cal M})$ be the corresponding
PROP of singular chains. Then there is a structure of a
$Chains({\cal M})$-algebra on
$C_{\bullet}^{mod} ( A, A )$, which is encoded in a collection of morphisms of
complexes
$$ Chains( {\cal M}(n,m)) \otimes C_{\bullet}^{mod} ( A, A )^{\otimes n} \to (C_{\bullet}^{mod} ( A, A ))^{\otimes m}. $$

In addition one has the following:

iii) If $A$ is homologically smooth and satisfies the degeneration property
then the structure of $Chains({\cal M})$-algebra extends to a structure of a $Chains(\overline{\cal M})$-algebra, where $\overline{\cal M}$ is the
topological PROP of stable compactifications of ${\cal M}(n,m)$.

\begin{defn} \label{definition of PROP}An element of ${\cal M}(n,m)$ is an isomorphism class of
triples $(\Sigma,h,mark)$ where $\Sigma$ is a compact oriented surface
(not necessarily connected) with metric $h$, and $mark$ is an orientation
preserving isometry between a neighborhood of $\partial\Sigma$ and
the disjoint union of $n+m$ flat semiannuli $\sqcup_{1\le i\le n}(S^1\times [0,\varepsilon))\sqcup
\sqcup_{1\le i\le m}(S^1\times [-\varepsilon,0])$, where $\varepsilon$ is a sufficiently small
positive number. We will call $n$ circles ``inputs" and the rest $m$ circles ``outputs".
We will assume that each connected component of $\Sigma$ has
at least one input, and there are no discs among the connected components.
Also we will add $\Sigma=S^1$ to ${\cal M}(1,1)$ as the identity morphism.
It can be thought of as the limit of cylinders $S^1\times [0,\varepsilon]$ as
$\varepsilon\to 0$.
\end{defn}

The composition is given by the natural gluing of surfaces.

Let us describe a construction of the action of $Chains({\cal M})$ on the Hochschild
chain complex. In fact, instead of $Chains({\cal M})$  we will consider a quasi-isomorphic dg-PROP
$R=R(n,m)_{n,m\ge 0}$
generated by ribbon graphs with additional data.
In what follows we will omit some technical details in the definition of the PROP $R$.
They can be recovered in a more or less straightforward way.

It is well-known (and can be
proved with the help of Strebel differentials) that ${\cal M}(n,m)$ admits a
stratification with strata parametrized by graphs described below.
More precisely, we consider the following class of graphs:

1) Each graph $\Gamma$ is a (not necessarily connected) ribbon graph (i.e. we are given a cyclic order on
the set $Star(v)$ of edges attached to a vertex $v$ of $\Gamma$). It is well-known that
replacing an edge of a ribbon graph by a thin stripe (thus getting
a ``fat graph") and gluing stripes in the cyclic order one
gets a Riemann surface with the boundary.

2) The set
$V( \Gamma )$ of vertices of $\Gamma$ is the union of three sets:
$V (\Gamma ) = V_{in}( \Gamma )
\cup V_{middle}(\Gamma)\cup V_{out}(\Gamma)$. Here $V_{in}(\Gamma)$ consists of
$n$ numbered vertices $in_1,...,in_n$ of the valency $1$ (the corresponding outcoming edges are called tails),
$V_{middle}( \Gamma )$ consists of vertices of the valency greater or equal
than $3$, and $V_{out}(\Gamma)$ consists of $m$ numbered vertices
$out_1,...,out_m$ of  valency
greater or equal than $1$.

3) We assume that the Riemann surface corresponding to $\Gamma$ has  $n$ connected
boundary components each of which has exactly one input vertex.

4) For every vertex $out_j\in V_{out}(\Gamma), 1\le j\le m$ we {\it mark}
either an incoming edge or a pair of adjacent edges
(we call such a pair of edges a {\it corner}).

{\vspace{3mm}}


\includegraphics{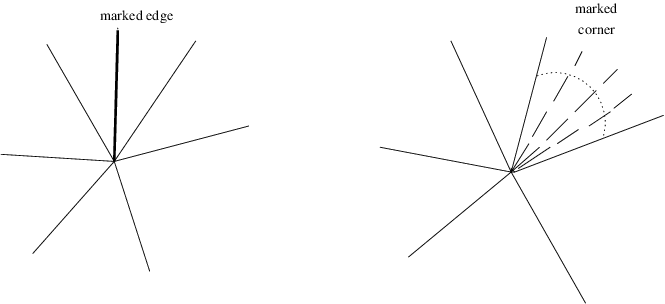}

{\vspace{3mm}}

More pedantically, let $E(\Gamma)$ denotes the set of edges of $\Gamma$ and $E^{or}(\Gamma)$ denotes the set of pairs
$(e,or)$ where $e\in E(\Gamma)$ and $or$ is one of two possible
orientations of $e$. There is an obvious map
$E^{or}(\Gamma)\to V(\Gamma)\times V(\Gamma)$ which assigns
to an oriented edge the pair of its endpoint vertices: source and target. The free involution $\sigma$
acting on $E^{or}(\Gamma)$ (change of orientation) corresponds
to the permutation map on $V(\Gamma)\times V(\Gamma)$.
Cyclic order on each $Star(v)$ means that there is a bijection
$\rho:E^{or}(\Gamma)\to E^{or}(\Gamma)$ such that orbits of
iterations $\rho^n,n \ge 1$ are
elements of $Star(v)$ for some $v\in V(\Gamma)$.
In particular, the corner is given either by a pair of coinciding
edges $(e,e)$ such that $\rho(e)=e$ or by a pair edges $e, e^{\prime}\in Star(v)$ such that $\rho(e)=e^{\prime}$.
Let us define a {\it face} as an orbit of $\rho\circ\sigma$.
Then faces are oriented closed paths. It follows from the condition 2)
that each face contains exactly one edge outcoming from some $in_i$.

We depict below two  graphs in the case $g = 0, n = 2,m=0$.



\includegraphics{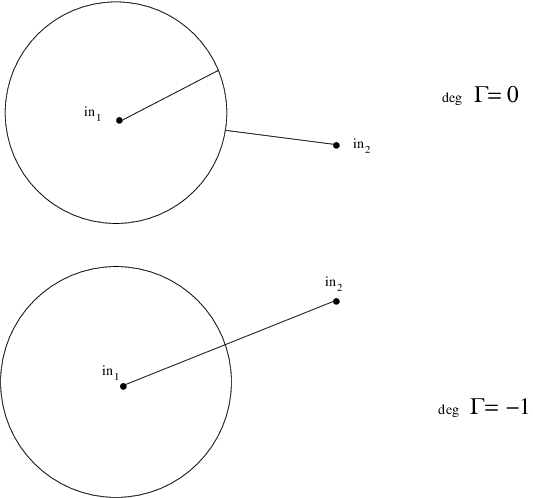}

{\vspace{3mm}}

Here is a picture illustrating the notion of face:

{\vspace{3mm}}


\includegraphics{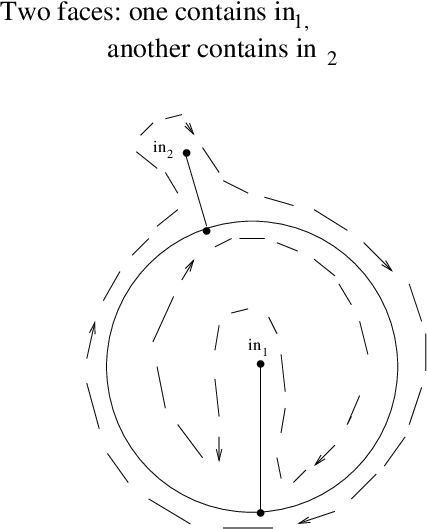}

{\vspace{3mm}}

\begin{rmk} \label{graphs have no automorphisms}The  above data (i.e. a ribbon graph with numerations
of $in$ and $out$ vertices) have no automorphisms. Thus we can
identify $\Gamma$ with its isomorphism class.

\end{rmk}






The functional $( m_n ( a_1, ..., a_n ), a_{n + 1} )$ is depicted such
as follows.

{\vspace{3mm}}


{\vspace{3mm}}

\includegraphics{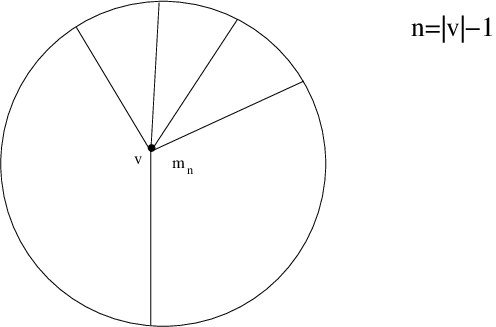}

{\vspace{3mm}}

We define the {\it degree} of $\Gamma$ by the formula

$$deg\,\Gamma=\sum_{v\in V_{middle}(\Gamma)}(3-|v|)+\sum_{v\in V_{out}(\Gamma)}(1-|v|)+\sum_{v\in V_{out}(\Gamma)}\epsilon_v-N\chi(\Gamma),$$
where $\epsilon_v=-1$, if $v$ contains a marked corner and
$\epsilon_v=0$ otherwise. Here $\chi(\Gamma)=|V(\Gamma)|-|E(\Gamma)|$ denotes the Euler characteristic
of $\Gamma$.

\begin{defn} \label{definition of PROP as graded vector space}We define $R(n,m)$ as a graded vector space which
is a direct sum $\oplus_{\Gamma}\psi_{\Gamma}$ of $1$-dimensional
graded vector spaces generated by graphs $\Gamma$ as above, each
summand has degree $deg\,\Gamma$.

\end{defn}

One can see that $\psi_{\Gamma}$ is naturally identified with the tensor
product of $1$-dimensional vector spaces (determinants) corresponding
to vertices of $\Gamma$.

Now, having a graph $\Gamma$ which satisfies conditions 1)-3) above, and
Hochschild chains $\gamma_1, ..., \gamma_n \in C_{\bullet}^{mod} ( A, A )$ we would like to define
an element of $C_{\bullet}^{mod} ( A, A )^{\otimes m}$.
Roughly speaking we are going to assign the above $n$ elements of the Hochschild
complex to $n$ faces corresponding to vertices $in_i, 1\le i\le n$, then
assign tensors corresponding to higher products $m_l$ to
internal vertices $v\in V_{middle}(\Gamma)$, then using the convolution
operation on tensors given by the scalar product on $A$ to read off
the resulting tensor from $out_j, 1\le j\le m$.
More precise algorithm is described below.

a) We decompose the modified Hochschild complex such as follows:

$$C_{\bullet}^{mod}(A,A)=\oplus_{l\ge 0, \varepsilon \in \{0,1\}}C_{l,\varepsilon}^{mod}(A,A),$$
where $C_{l,\varepsilon=0}^{mod}(A,A)=A\otimes (A[1])^{\otimes l}$
and $C_{l,\varepsilon=1}^{mod}(A,A)=\K\otimes (A[1])^{\otimes l}$
according to the definition of modified Hochschild chain complex.
For any choice of $l_i\ge 0,\varepsilon_i\in \{0,1\}, 1\le i\le n$
we are going to construct a linear map of degree zero
$$f_{\Gamma}:\psi_{\Gamma}\otimes C_{l_1,\varepsilon_1}^{mod}(A,A)\otimes...\otimes C_{l_n,\varepsilon_1}^{mod}(A,A)\to (C_{\bullet}^{mod}(A,A))^{\otimes m}.$$

The result will be a sum
$f_{\Gamma}=\sum_{\Gamma^{\prime}}f_{\Gamma^{\prime}}$ of certain
maps. The description of the collection of graphs $\Gamma^{\prime}$
is given below.

b) Each new graph $\Gamma^{\prime}$ is obtained from $\Gamma$ by adding new
edges. More precisely one has
$V(\Gamma^{\prime})=V(\Gamma)$ and for each vertex $in_i\in V_{in}(\Gamma)$
we add $l_i$ new outcoming edges. Then the valency of $in_i$ becomes
$l_i+1$.

\vspace{3mm}


\includegraphics{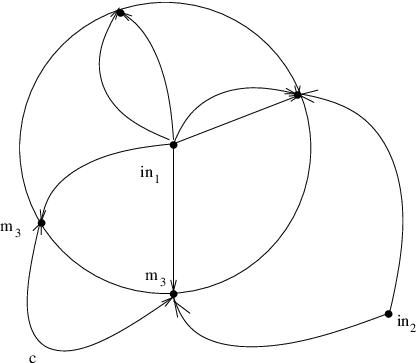}

{\vspace{3mm}}

More pedantically, for every $i, 1\le i\le n$ we have constructed
a map from the set $\{1,...,l_i\}$ to a cyclically ordered
set which is an orbit of $\rho\circ \sigma$ with removed
the tail edge oucoming from $in_i$. Cyclic order on the
edges of $\Gamma^{\prime}$ is induced by the cyclic order
at every vertex and the cyclic order on the path forming the
face corresponding to $in_i$.

\vspace{3mm}


\includegraphics{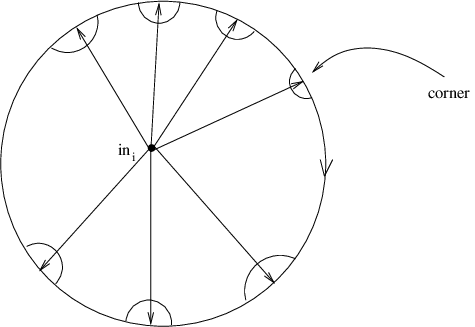}

{\vspace{3mm}}

c) We assign  $\gamma_i\in C_{l_i,\varepsilon_i}$ to $in_i$.
We depict $\gamma_i$ as a ``wheel" representing the Hochschild
cocycle. It is formed by the endpoints of the $l_i+1$ edges outcoming
from $in_i\in V(\Gamma^{\prime})$ and taken in the cyclic order of
the corresponding face. If $\varepsilon_i=1$ then (up to a scalar)
$\gamma_i=1\otimes a_1\otimes...\otimes a_{l_i}$, and we require that the
tensor factor $1$ corresponds to zero in the cyclic order.

\vspace{3mm}


\includegraphics{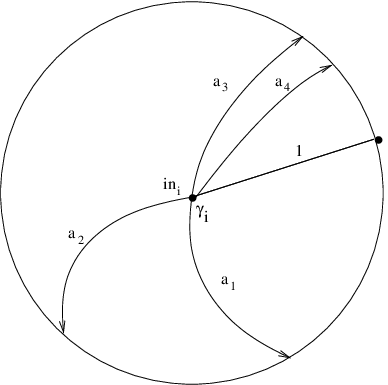}

{\vspace{3mm}}

d) We remove from considerations graphs $\Gamma$ which do not obey the
following property after the step c):

{\it the edge corresponding to the unit $1\in \K$ (see step c))
is of the type $(in_i,v)$ where either $v\in V_{middle}(\Gamma^{\prime})$
and $|v|=3$ or $v=out_j$ for some $1\le j\le m$ and the edge
$(in_i,out_j)$ was the marked edge for $out_j$}.

Let us call {\it unit edge} the one which satisfies one of the above properties.
We define a new graph $\Gamma^{\prime\prime}$ which is obtained
from $\Gamma$ by removing unit edges.

e) Each vertex now has the valency
$|v|\ge 2$. We attach to every such vertex either:

the tensor $c\in A\otimes A$ (inverse to the scalar product),
if $|v|=2$,

or

the tensor $(m_{|v|-1}(a_1,...,a_{|v|-1}),a_{|v|})$ if $|v|\ge 3$.
The latter can be identified with the element of $A^{\otimes |v|}$
(here we use the non-degenerate scalar product on $A$).

Let us illustrate this construction.



\includegraphics{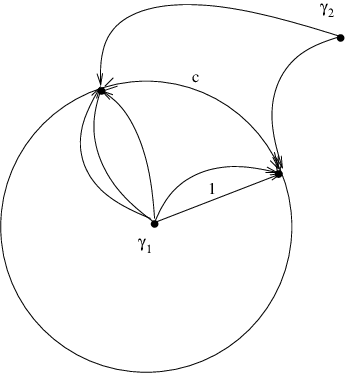}

{\vspace{3mm}}

f) Let us contract indices of tensors corresponding to
$V_{in}(\Gamma^{\prime\prime})\cup V_{middle}(\Gamma^{\prime\prime})$ (see c), e)) along
the edges of $\Gamma^{\prime\prime}$ using the
scalar product on $A$. The result will be an element $a_{out}$ of the tensor product
$\otimes_{1\le j\le m}A^{Star_{\Gamma^{\prime\prime}}(out_j)}$.

g) Last thing we need to do is to interpret the element $a_{out}$
as
an element of $C_{\bullet}^{mod}(A,A)$.
There are three cases.

{\it Case 1.}
When we constructed $\Gamma^{\prime\prime}$ there was a unit
edge incoming to some $out_j$. Then we reconstruct back the removed edge,
attach $1\in \K$ to it, and interpret the resulting
tensor as an element of $C_{|out_j|,\varepsilon_j=1}^{mod}(A,A)$.

{\it Case 2.}
There was no removed unit edge incoming to $out_j$ and we
had a marked edge (not a marked corner) at the vertex $out_j$.
Then we have an honest element of
$C_{|out_j|,\varepsilon_j=0}^{mod}(A,A)$

{\it Case 3.}
Same as in Case 2, but there was a marked corner at $out_j\in V_{out}(\Gamma)$. We have added and removed  new edges when
constructed $\Gamma^{\prime\prime}$. Therefore the marked
corner gives rise to a new set of marked corners at $out_j$ considered
as a vertex of $\Gamma^{\prime\prime}$. Inside every such a corner we
insert a new edge, attach the element $1\in \K$ to it and take
the sum over all the corners. In this way we obtain an element
of $C_{|out_j|,\varepsilon_j=1}^{mod}(A,A)$.
This procedure is depicted below.



\vspace{3mm}

\includegraphics{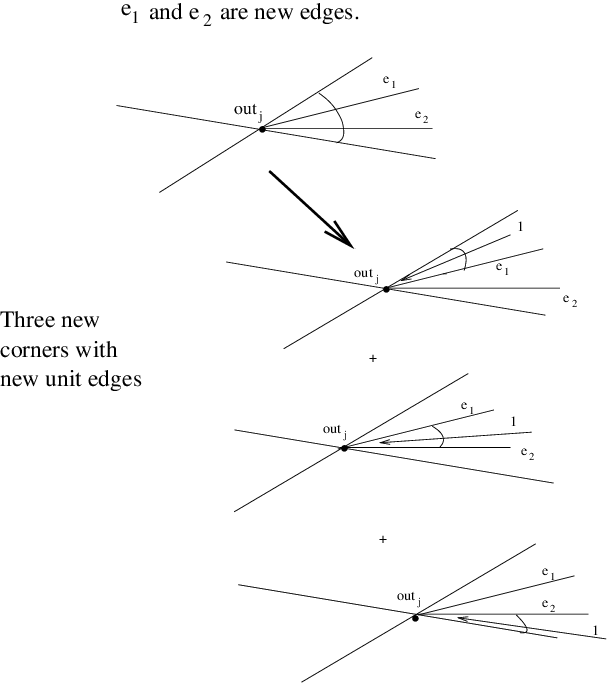}

{\vspace{3mm}}

This concludes the construction of $f_{\Gamma}$.
Notice that
$R$ is a dg-PROP with the differential given by the
insertion of a new
edge between two vertices from $V_{middle} ( \Gamma )$.

Proof of the following Proposition is tedious and technical.  We will omit details here.

\begin{prp} \label{algebra over PROP}The above construction gives rise to a structure of a
$R$-algebra on $C_{\bullet}^{mod}(A,A)$.

\end{prp}

\begin{rmk}\label{homological smoothness was not used}
The above construction did not use homological smoothness of $A$.
\end{rmk}

Finally we would like to say few words about an extension of the $R$-action
to the $Chains(\overline{{\cal M}})$-action. 

If we assume the degeneration property for $A$,
then the action of the PROP $R$ can be extended to
the action of the PROP $Chains(\overline{\cal M})$ of
singular chains of the topological PROP of stable degenerations
of $M_{g,n,m}^{marked}$. In order to see this, one introduces the
PROP $D$ freely generated by $R(2,0)$ and $R(1,1)$, i.e. by singular chains
on the moduli space of cylinders with two inputs and zero outputs (they correspond
to the scalar product on $C_{\bullet}(A,A)$) and by cylinders with one input and one
output (they correspond to morphisms $C_{\bullet}(A,A)\to C_{\bullet}(A,A)$).
In fact the (non-symmetric) bilinear form
$h:H_{\bullet}(A,A)\otimes H_{\bullet}(A,A)\to \K$ does exist
for any compact $\A$-algebra $A$. It is described by the graph
of degree zero on the corresponding figure above.
This is a generalization of the bilinear form
$(a,b)\in A/[A,A]\otimes A/[A,A]\mapsto Tr(axb)\in \K$.

It seems plausible that homological smoothness
implies that $h$ is non-degenerate. This allows us to extend the action
of the dg sub-PROP $D\subset R$ to the action of the dg PROP
$D^{\prime}\subset R$
which contains also $R(0,2)$ (i.e. the inverse to the above bilinear form).
If we assume the degeneration property, then we  can ``shrink"
the action of the homologically non-trivial circle of the cylinders,
since the rotation around this circle corresponds to the differential $B$.
Thus $D^{\prime}$ is quasi-isomorphic to the dg-PROP of chains
on the one-dimensional retracts of the above cylinders, where the retractions
contracts the circles.
Let us denote the dg-PROP generated by singular chains on the retractions by $D^{\prime\prime}$.
Thus, assuming the degeneration property, we see that
the free product dg-PROP $R^{\prime}=R\ast_DD^{\prime\prime}$
acts on $C_{\bullet}^{mod}(A,A)$.
One can show that $R^{\prime}$ is quasi-isomorphic to the dg-PROP
of chains on the topological PROP $\overline{M}_{g,n,m}^{marked}$ of
stable compactifications of the surfaces from $M_{g,n,m}^{marked}$.

\begin{rmk} \label{TFT}a) The above construction is a generalization of
the construction from [Ko3], which assigns  cohomology classes
of ${M}_{g,n}$ to a finite-dimensional $\A$-algebra
with scalar product (trivalent graphs were used in [Ko3]).

b) Different approach to the action of the PROP $R$ was suggested in [Cos].
The above Proposition \ref{algebra over PROP} gives rise to a structure of Topological
Field Theory associated with a non-unital $\A$-algebra with scalar
product (or more generally with a Calabi-Yau category). If the degeneration property holds for $A$ then one can define
a Cohomological Field Theory in the sense of [KoM]

c) Homological smoothness of $A$ is closely related to the existence
of a non-commutative analog of the Chern class of the diagonal
$\Delta\subset X\times X$ of a projective scheme $X$. This Chern
class gives rise to the inverse to the scalar product on $A$.
Generalization of this observation plays an important role in the case of $\A$-categories.
\end{rmk}

\section{Appendix}\label{appendix}

\subsection{Non-commutative schemes and ind-schemes}\label{nc schemes}
Let $\cal{C}$ be an abelian $\K$-linear tensor category. To simplify
formulas we will assume that it is strict (see [McL]). We will also assume
that $\cal{C}$ admits infinite sums. To simplify the exposition we will
assume below as we did in most cases  in the main body of the paper that ${\cal C} =
Vect_\K^{\Z}$.

\begin{defn} \label{categpry of nc schemes}The category of non-commutative affine $\K$-schemes in $\cal{C}$
(notation $NAff_{{\cal C}}$) is defined as the category opposite to the category of
associative unital $\K$-algebras in $\cal{C}$.
\end{defn}

The non-commutative scheme corresponding to the algebra $A$
is denoted by $Spec(A)$. Conversely, if $X$ is a non-commutative
affine scheme then the corresponding algebra  is called the algebra
of regular functions on $X$ and is denoted by ${\cal O}(X)$.
By analogy with commutative case we call a morphism
$f:X\to Y$ a {\it closed embedding} if the corresponding homomorphism
$f^{\ast}:{\cal O}(Y)\to {\cal O}(X)$ is an epimorphism.

Let us recall some terminology related to ind-objects (see for ex. [Gr], [AM],
[KSch]). For a covariant functor $\phi : I \to {\cal A}$ from a small
filtering category $I$ (called filtrant in [KSch]) there is a notion of an
inductive limit $`` \varinjlim  " \phi \in \widehat{{\cal A}}$
and a projective limit $`` \varprojlim  " \phi \in
\widehat{{\cal A}}$, where $\widehat{{\cal A}}$ is the category of functors ${\cal A}\to Sets$.
By definition $`` \varinjlim  " \phi ( X ) =
\varinjlim Hom_{{\cal A}} ( X, \phi ( i ) )$ and $`` \varprojlim
" \phi ( X ) = \varinjlim Hom_{{\cal A}} ( \phi ( i ), X )$. All inductive
limits form a full subcategory $Ind ( {\cal A} ) \subset
\widehat{{\cal A}}$ of ind-objects in $\cal{A}$. Similarly all
projective limits form a full subcategory $Pro ( {\cal A} ) \subset
\widehat{{\cal A}}$ of pro-objects in $\cal{A}$.

\begin{defn} \label{nc ind-schemes}Let $I$ be a small filtering category, and $F : I \to
NAff_{{\cal C}}$ a covariant functor. We say that $`` \varinjlim
"F$ is a non-commutative ind-affine scheme if for a morphism $i \to j$ in $I$
the corresponding morphism $F ( i ) \to F ( j )$ is a closed embedding.
\end{defn}

In other words a non-commutative ind-affine scheme $X$ is an object of $Ind (
NAff_{{\cal C}} )$, corresponding to the projective limit $\varprojlim
\hspace{0.25em} A_{\alpha}, \alpha \in I$, where each $A_{\alpha}$ is a unital
associative algebra in $\cal{C}$, and for a morphism $\alpha \to \beta$ in
$I$ the corresponding homomorphism $A_{\beta} \to A_{\alpha}$ is a surjective
homomorphism of unital algebras (i.e. there exists an exact sequence $0 \to ? \to
A_{\beta} \to A_{\alpha} \to 0$).

\begin{rmk}\label{categorical epimorphisms} Not all categorical epimorphisms of algebras are surjective
homomorphisms (although the converse is true). Nevertheless one can define
closed embeddings of affine schemes for an arbitrary abelian $\K$-linear
category, observing that a surjective homomorphism of algebras $f : A \to B$
is characterized categorically by the condition that $B$ is the cokernel of
the pair of the natural projections $f_{1, 2} : A \times_B A \to A$ defined by
$f$.
\end{rmk}

Morphisms between non-commutative ind-affine schemes are defined as morphisms
between the corresponding projective systems of unital algebras. Thus we have
$$Hom_{NAff_{{\cal C}}} ( \varinjlim_I X_i, \varinjlim_J Y_j ) =
   \varprojlim_I \varinjlim_J Hom_{NAff_{{\cal C}}} ( X_i, Y_j ) . $$
Let us recall that an algebra $M \in Ob ( {\cal C} )$ is called nilpotent
if the natural morphism $M^{\otimes n} \to M$ is zero for all sufficiently
large $n$.

\begin{defn} \label{formal ind-schemes}A non-commutative ind-affine scheme $\hat{X}$ is called formal if it
can be represented as $\hat{X} = \varinjlim Spec ( A_i )$, where $( A_i )_{i
\in I}$ is a projective system of associative unital algebras in $\cal{C}$
such that the homomorphisms $A_i \to A_j$ are surjective and have nilpotent
kernels for all morphisms $j \to i$ in $I$.
\end{defn}

Let us consider few examples in the case when ${\cal C} = Vect_\K\subset Vect_\K^{\Z}$.

\begin{exa} \label{formal line}In order to define the non-commutative formal affine line
$\widehat{{\bf A}}_{NC}^1$ it suffices to define $Hom ( Spec ( A ),
\widehat{{\bf A}}_{NC}^1 )$ for any associative unital algebra $A$. We
define $Hom_{NAff_\K} ( Spec ( A ), \widehat{{\bf A}}_{NC}^1 ) =
\varinjlim \hspace{0.25em} Hom_{Alg_\K} ( \K [ [ t ] ] / ( t^n ), A )$. Then the
set of $A$-points of the non-commutative formal affine line consists of all
nilpotent elements of $A$.
\end{exa}

\begin{exa} \label{formal affine space}For an arbitrary set $I$ the non-commutative formal affine space
$\widehat{{\bf A}}_{NC}^I$ corresponds, by definition, to the topological
free algebra $\K \langle \langle t_i \rangle \rangle_{i \in I}$. If $A$ is a
unital $\K$-algebra then any homomorphism $\K \langle \langle t_i \rangle
\rangle_{i \in I} \to A$ maps almost all $t_i$ to zero, and the remaining
generators are mapped into nilpotent elements of $A$.
In particular, if $I =
{\bf N} = \{ 1, 2, ... \}$ then
$\widehat{{\bf A}}_{NC}^{{\bf N}} =
\varinjlim Spec (\K \langle \langle t_1, ..., t_n \rangle
\rangle / ( t_1, ..., t_n )^m),$ where $( t_1, ..., t_n )$ denotes the
two-sided ideal generated by $t_i, 1 \le i \le n$, and the limit is taken over
all $n, m \to \infty$.
\end{exa}

By definition, a {\it closed subscheme} $Y$ of a scheme $X$ is defined by a
2-sided ideal $J \subset {\cal O} ( X )$. Then ${\cal O} ( Y ) ={\cal O}(X)/J$. If $Y\subset X$ is defined by a
$2$-sided ideal $J \subset {\cal O} ( X )$, then the completion of $X$
along $Y$ is a formal scheme corresponding to the projective limit of algebras
$\varprojlim_n {\cal O} ( X ) / J^n$. This formal scheme will be denoted by
$\hat{X}_Y$ or by $Spf ( {\cal O} ( X ) / J )$.

Non-commutative affine schemes over a given field $\K$ form a symmetric monoidal
category. The tensor structure is given by the {\it ordinary tensor
product} of unital algebras. The corresponding tensor product of
non-commutative affine schemes will be denoted by $X \otimes Y$. It is not a
categorical product, differently from the case of commutative affine schemes,
where the tensor product of algebras corresponds to the Cartesian product $X
\times Y$. For non-commutative affine schemes the analog of the Cartesian
product is the {\it free product} of algebras.

Let $A, B$ be free algebras. Then $Spec ( A )$ and $Spec ( B )$ are
non-commutative manifolds. Since the tensor product $A \otimes B$ in general
is not a smooth algebra, the non-commutative affine scheme $Spec ( A \otimes B
)$ in general is not a manifold.

Let $X$ be a non-commutative ind-affine scheme in $\cal{C}$. A closed
$\K$-point $x \in X$ is by definition a homomorphism of ${\cal O}(X)$ to
the tensor algebra generated by the unit object ${\bf 1}$.
Let
$m_x$ be the kernel of this homomorphism. We define the {\it tangent space} $T_x X$
in the usual way as $( m_x / m_x^2 )^{\ast} \in Ob ( {\cal C} )$. Here
$m_x^2$ is the image of the multiplication map $m_x^{\otimes 2} \to m_x$.

A non-commutative ind-affine scheme with a marked closed $\K$-point will be
called \textit{pointed}. There is a natural generalization of this notion to
the case of many points. Let $Y \subset X$ be a closed subscheme of disjoint
closed $\K$-points (it corresponds to the algebra homomorphism ${\cal O} ( X
) \to {\bf 1} \oplus {\bf 1} \oplus ...$). Then $\hat{X}_Y$ is a formal
manifold. A pair $( \hat{X}_Y, Y )$ (often abbreviated by $\hat{X}_Y$) will be
called (non-commutative) {\it formal manifold with marked points}. If $Y$
consists of one such point then $( \hat{X}_Y, Y )$ will be called a
(non-commutative) {\it formal pointed manifold}.

\subsection {Proof of Theorem \ref{coalgebras as functors}}\label{proof of theorem 2.1.1}

In the category $Alg_{{\cal C}^f}$ every pair of
morphisms has a kernel. Since the functor $F$ is left exact and the category
$Alg_{{\cal C}^f}$ is artinian, it follows from [Gr], Sect. 3.1 that $F$ is
strictly pro-representable. This means that there exists a projective system
of finite-dimensional algebras $( A_i )_{i \in I}$ such that, for any morphism
$i \to j$ the corresponding morphism $A_j \to A_i$ is a categorical
epimorphism, and for any $A \in Ob ( Alg_{\mathcal{C}^f} )$ one has
$$ F ( A ) = \varinjlim_I Hom_{Alg_{{\cal C}^f}} ( A_i, A ) .$$
Equivalently,
$$ F ( A ) = \varinjlim_I Hom_{Coalg_{{\cal C}^f}} ( A_i^{\ast}, A^{\ast}
   ), $$
where $( A_i^{\ast} )_{i \in I}$ is an inductive system of finite-dimensional
coalgebras and for any morphism $i \to j$ in $I$ we have a categorical
monomorphism $g_{ji} : A_i^{\ast} \to A_j^{\ast}$.

All what we need is to replace the projective system of algebras $( A_i )_{i
\in I}$ by another projective system of algebras $( \overline{A}_i )_{i \in
I}$ such that

a) functors $`` \varprojlim "h_{A_i}$ and
$`` \varprojlim "h_{\overline{A}_i}$ are isomorphic (here $h_X$ is the functor
defined by the formula $h_X ( Y ) = Hom ( X, Y )$);

b) for any morphism $i \to j$ the corresponding homomorphism of algebras
$\overline{f}_{ij} : \overline{A}_j \to \overline{A}_i$ is surjective.

Let us define $\overline{A}_i = \bigcap_{i \to j} Im ( f_{ij} )$, where $Im (
f_{ij} )$ is the image of the homomorphism $f_{ij} : A_j \to A_i$
corresponding to the morphism $i \to j$ in $I$. In order to prove a) it
suffices to show that for any unital algebra $B$ in ${\cal C}^f$ the
natural map of sets
$$\varinjlim_I Hom_{{\cal C}^f} ( A_i, B ) \to \varinjlim_I
   Hom_{{\cal C}^f} ( \overline{A}_i, B )$$
(the restriction map) is well-defined and bijective.

The set $\varinjlim_I Hom_{{\cal C}^f} ( A_i, B )$ is isomorphic to $(
\bigsqcup_I Hom_{{\cal C}^f} ( A_i, B ) ) / equiv$, where two maps $f_i :
A_i \to B$ and $f_j : A_j \to B$ such that $i \to j$ are equivalent if $f_i
f_{ij} = f_j$. Since ${\cal C}^f$ is an artinian category, we conclude that
there exists $A_m$ such that $f_{im} ( A_m ) = \overline{A}_i$, $f_{jm} ( A_m
) = \overline{A}_j$. From this observation one easily deduces that $f_{ij} (
\overline{A}_j ) = \overline{A}_i$. It follows that the morphism of functors
in a) is well-defined, and b) holds. The proof that morphisms of functors
biejectively correspond to homomorphisms of coalgebras is similar. This
completes the proof of the theorem. $\blacksquare$

\subsection {Proof of Proposition \ref{coalgebra as union}}\label{proof of proposition 2.1.2}

The result follows from the fact that any $x \in B$
belongs to a finite-dimensional subcoalgebra $B_x \subset B$, and if $B$ was
counital then $B_x$ would be also counital. Let us describe how to construct
$B_x$. Let $\Delta$ be the coproduct in $B$. Then one can write
$$ \Delta ( x ) = \sum_i a_i \otimes b_i, $$
where $a_i$ (resp. $b_i$) are linearly independent elements of $B$.

It follows from the coassociativity of $\Delta$ that
$$ \sum_i \Delta ( a_i ) \otimes b_i = \sum_i a_i \otimes \Delta ( b_i ) . $$
Therefore one can find constants $c_{ij} \in \K$ such that
$$ \Delta ( a_i ) = \sum_j a_j \otimes c_{ij}, $$
and
$$ \Delta ( b_i ) = \sum_j c_{ji} \otimes b_j . $$
Applying $\Delta \otimes id$ to the last equality and using the
coassociativity condition again we get
$$ \Delta ( c_{ji} ) = \sum_n c_{jn} \otimes c_{ni} . $$
Let $B_x$ be the vector space spanned by $x$ and all elements $a_i, b_i,
c_{ij}$. Then $B_x$ is the desired subcoalgebra. $\blacksquare$

\subsection{Formal completion along a subscheme}\label{formal completion along subscheme}

Here we present a construction which generalizes the definition of
a formal neighborhood of a $\K$-point of a non-commutative smooth thin scheme.

Let $X=Spc(B_X)$ be such a scheme and $f:X\to Y=Spc(B_Y)$ be a closed embedding,
i.e. the corresponding homomorphism of coalgebras $B_X\to B_Y$ is injective.
We start with  the category $\mathcal{N}_X$ of nilpotent
extensions of $X$, i.e. homomorphisms $\phi : X \to U$, where $U = Spc ( D )$
is a non-commutative thin scheme, such that the quotient $D / f ( B_X )$ (which
is always a non-counital coalgebra) is locally conilpotent. We recall that the
local conilpotency means that for any $a \in D / f ( B_X )$ there exists $n \ge
2$ such that $\Delta^{( n )} ( a ) = 0$, where $\Delta^{( n )}$ is the $n$-th
iterated coproduct $\Delta$. If $( X, \phi_1, U_1 )$ and $( X, \phi_2, U_2 )$
are two nilpotent extensions of $X$ then a morphism between them is a morphism
of non-commutative thin schemes $t : U_1 \to U_2$, such that $t \phi_1 =
\phi_2$. In particular, $\mathcal{N}_X$ is a subcategory of the naturally
defined category of non-commutative relative thin schemes.

Let us consider the functor $G_f : \mathcal{N}_X^{op} \to Sets$ such that $G (
X, \phi, U )$ is the set of all morphisms $\psi : U \to Y$ such that $\psi
\phi = f$.

\begin{prp}\label{representing triple}
Functor $G_f$ is represented by a triple $( X, \pi, \hat{Y}_X )$ where the
non-commutative thin scheme denoted by $\widehat{Y}_X$ is called the formal
neighborhood of $f ( X )$ in $Y$ (or the completion of $Y$ along $f ( X )$).
\end{prp}

{\it Proof.} Let $B_f \subset B_X$ be the counital subcoalgebra
which is the pre-image of the (non-counital) subcoalgebra in $B_Y / f ( B_X )$
consisting of locally conilpotent elements. Notice that $f ( B_X ) \subset B_f$.
It is easy to see that taking $\widehat{Y}_X : = Spc ( B_f )$ we obtain the triple
which represents the functor $G_f$. $\blacksquare$

Notice that $\widehat{Y}_X \to Y$ is a closed embedding of non-commutative thin
schemes.

\begin{prp}\label{double completion}
If $Y$ is smooth then $\widehat{Y}_X$ is smooth and $\widehat{Y}_X \simeq
\widehat{Y}_{\widehat{Y}_X}$.

\end{prp}

{\it Proof.} Follows immediately from the explicit description
of the coalgebra $B_f$ given in the proof of the  Proposition \ref{representing triple}.
$\blacksquare$

{\vspace{10mm}}

{\bf References}

{\vspace{3mm}}

[AM] M. Artin, B. Mazur, Etale Homotopy, Lect. Notes Math., 100, 1969.

{\vspace{2mm}}

[BD1] A. Beilinson, V. Drinfeld, Chiral algebras, AMS, 2004.

{\vspace{2mm}}

[BD2] A. Beilinson, V. Drinfeld, Quantization of Hitchin's integrable system
and Hecke eigensheaves, preprint, University of Chicago.

{\vspace{2mm}}

[BK] A. Bondal, M. Kapranov, Enhanced triangulated categories, Math. USSR
Sbornik, 70:1, 1991, 93-107.

{\vspace{2mm}} [BvB] A. Bondal, M. Van den Bergh Generators and
representability of functors in commutative and noncommutative geometry,
arXiv:math.AG/0204218. 

{\vspace{2mm}}

[BV] J. Boardman, R. Vogt, Homotopy invariant algebraic structures on
topological spaces, Lect. Notes Math. 347, 1973.

{\vspace{2mm}}

[Co] A. Connes, Non-commutative geometry. Academic Press, 1994.

{\vspace{2mm}}

[Cos] K. Costello, Topological conformal field theories, Calabi-Yau categories
and Hochschild homology, preprint (2004).

{\vspace{2mm}}

[CQ1] J. Cuntz, D. Quillen, Cyclic homology and non-singularity, J. of AMS,
8:2, 1995, 373-442.

{\vspace{2mm}}

[CQ2] J. Cuntz, D. Quillen, Algebra extensions and non-singularity, J. of AMS,
8:2,1995, 251-289.

{\vspace{2mm}}

[CST] J. Cuntz, G. Skandalis, B. Tsygan, Cyclic homology in noncommutative
geometry, Encyclopaedia of Mathematical Sciences, v. 121, Springer Verlag,
p. 74-113.

{\vspace{2mm}}
[DI] P. Deligne, L. Illusie, Relevements modulo $p^2$ et decomposition
du complex de de Rham, Invent. Math. 89, 1987, 247-270.

{\vspace{2mm}}

[DM] P. Deligne, J.S. Milne, Tannakian categories, Lect. Notes Math., 900
(1982), 101-228.

{\vspace{2mm}}

[Dr] V. Drinfeld, DG quotients of DG categories, arXiv:math.KT/0210114.

{\vspace{2mm}}

[FOOO] K. Fukaya, Y.G. Oh, H. Ohta, K. Ono, Lagrangian intersection
Floer theory-anomaly and obstruction. AMS, 2009.

{\vspace{2mm}}

[Get] E. Getzler, Cartan homotopy formulas and Gauss-Manin connection in
cyclic homology, Israel Math. Conf. Proc., 7, 65-78.

{\vspace{2mm}}

[Gi1] V. Ginzburg, Non-commutative symplectic geometry, quiver
varieties and operads, arXiv:math.QA/0005165.

{\vspace{2mm}}

[Gi2] V. Ginzburg, Lectures on Noncommutative
Geometry, arXiv:math.AG/0506603.

{\vspace{2mm}}

[Gi3] V. Ginzburg, Double derivations and cyclic homology, arXiv:math.KT/0505236.

{\vspace{2mm}}

 [Gr] A. Grothendieck, Technique de descente.II. Sem.
Bourbaki, 195, 1959/60).

{\vspace{2mm}}

[KSch] M. Kashiwara, P. Schapira, Ind-sheaves, Asterisque 271, 2001.

{\vspace{2mm}}

[KaKoP], L.Katzarkov, M. Kontsevich, T. Pantev, Hodge theoretic aspects
of mirror symmetry, arXiv:0806.0107.
{\vspace{2mm}}

[Kal] D. Kaledin, Non-commutative Cartier operator and Hodge-to-de Rham
degeneration, arXiv:math.AG/0511665.

{\vspace{2mm}}

[Kau1] R. Kaufmann,Moduli space actions on the Hochschild cochains of a Frobenius algebra I:
Cell Operads, arXiv:math.AT/0606064.

{\vspace{2mm}}

[Kau2] R. Kaufmann,Moduli space actions on the Hochschild cochains of a Frobenius algebra II: Correlators,
arXiv:math.AT/0606065.

{\vspace{2mm}}

[Ke 1] B. Keller, Introduction to A-infinity algebras and
modules, Homology, Homotopy and Applications 3, 2001, 1-35.

{\vspace{2mm}}

[Ke2] B. Keller, On differential graded categories,
arXiv: math.KT/0601185.

{\vspace{2mm}}

[Ko1] M. Kontsevich, Deformation quantization of Poisson manifolds,
arXiv:math.QA/9709040.

{\vspace{2mm}}

[Ko2] M. Kontsevich, Formal non-commutative symplectic geometry. In: Gelfand
Mathematical Seminars, 1990-1992, p. 173-187. Birkhauser Boston, MA, 1993.

{\vspace{2mm}}

[Ko3] M. Kontsevich, Feynman diagrams and low-dimensional topology.
Proc. Europ. Congr. Math. vol. 1, 1992.

{\vspace{2mm}}

[Ko4] M. Kontsevich, Notes on motives in finite characteristic, arXiv:arXiv:math/0702206.

{\vspace{2mm}}

[Ko5] M. Kontsevich, Lecture on triangulated $\A$-categories at Max-Planck Institut f\"ur Mathematik,
2001.

{\vspace{2mm}}

[KoM] M. Kontsevich, Yu.Manin, Gromov-Witten classes,
quantum cohomology and enumerativ geometry,
Comm. Math. Phys., 164:3,1994,525-562.

{\vspace{2mm}}

[KoSo1] M. Kontsevich, Y. Soibelman, Deformations of algebras over operads and
Deligne conjecture, math.QA/0001151, published in Lett. Math. Phys. 21:1
2000, 255-307.

{\vspace{2mm}}

[KoSo2] M. Kontsevich, Y. Soibelman, Deformation theory, book in
preparation, draft available at {https://www.math.ksu.edu/$\sim$ soibel/}.

{\vspace{2mm}} [KorSo] L. Korogodski, Y. Soibelman, Algebras of functions on
quantum groups.I. Math. Surveys and Monographs 56, AMS, 1998.

{\vspace{2mm}}

[Le Br] L. Le Bruyn, Non-commutative geometry an $n$, book in preparation.

{\vspace{2mm}}

[Lyu] V. Lyubashenko, Category of $\A$-categories, arXiv:math.CT/0210047.

{\vspace{2mm}}

[LyuOv] V. Lyubashenko, S. Ovsienko, A construction of $\A$-category,
arXiv:math.CT/0211037.

{\vspace{2mm}}

[McL] S. Mac Lane, Categories for the working mathematician. Springer-Verlag,
1971.

{\vspace{2mm}} 

[Or] D. Orlov,Triangulated categories of singularities and equivalences between Landau-Ginzburg models, arXiv:math.AG/0503630.

{\vspace{2mm}} 

[R]  R. Rouquier,  Dimensions of triangulated categories, arXiv:math.CT/0310134.

{\vspace{2mm}}

[Se] P. Seidel, Homological mirror symmetry for the quartic
surface, arXiv:math.SG/0310414.

{\vspace{2mm}}

[So1] Y. Soibelman, Non-commutative geometry and deformations of $\A$-algebras and $\A$-categories,
Arbeitstagung 2003, Preprint Max-PLanck Institut f\"ur Mathematik, MPIM2003-60h, 2003.

{\vspace{2mm}}

[So2] Y. Soibelman, Mirror symmetry and non-commutative geometry of $\A$-categories,
Journal of Mathematical Physics, vol. 45:10, 2004, 3742-3757.

{\vspace{2mm}}

[SU] S. Sanablidze, R. Umble, A diagonal of the associahedra, arXiv:math. AT/0011065.

{\vspace{2mm}}

[TaT1] D. Tamarkin, B. Tsygan, Non-commutative differential calculus, homotopy
BV-algebras and formality conjectures, arXiv:math.KT/0002116.

{\vspace{2mm}}

[TaT2] D. Tamarkin, B. Tsygan, The ring of differential operators
on forms in non-commutative calculus. Proceedings of Symposia in Pure Math., vol. 73, 2005, 105-131.

{\vspace{2mm}}

[ToVa] B. Toen, M. Vaquie, Moduli of objects in dg-categories, arXiv:math.AG/0503269.

{\vspace{2mm}}

[Ve] J. L. Verdier, Categories derivees, etat 0. Lect. Notes in Math. vol 569,
262-312, 1977.

{\vspace{5mm}}

Addresses:

M.K.: IHES, 35 route de Chartres, F-91440, France,
{email: maxim@ihes.fr}

Y.S.: Department of Mathematics, KSU, Manhattan, KS 66506, USA

{email: soibel@math.ksu.edu}

\end{document}